\documentclass[10pt]{article}
\usepackage{amsfonts}
\usepackage{bbm}
\usepackage{mathrsfs}%
\usepackage{amssymb,amsmath,amsthm, amsfonts}

\def\sqr#1#2{{\vcenter{\vbox{\hrule height.#2pt
              \hbox{\vrule width.#2pt height#1pt \kern#1pt \vrule width.#2pt}
          \hrule height.#2pt}}}}
\def\sqr#1#2{{\vcenter{\vbox{\hrule height.#2pt
              \hbox{\vrule width.#2pt height#1pt \kern#1pt \vrule width.#2pt}
              \hrule height.#2pt}}}}
\def\3n{\negthinspace \negthinspace \negthinspace }
\def\2n{\negthinspace \negthinspace }
\def\1n{\negthinspace }

\def\={\buildrel \triangle \over =}

\def\exp{\mathop{\rm exp}}
\def\sup{\mathop{\rm sup}}
\def\inf{\mathop{\rm inf}}

\def\inf{\hbox{\rm inf$\,$}}

\def\sup{\mathop{\rm sup}}
\def\inf{\mathop{\rm inf}}

\def\|{\Big |}
\def\({\Big (}
\def\){\Big )}
\def\[{\Big[}
\def\]{\Big]}
\def\be{\begin{equation}}
\def\bel{\begin{equation}\label}
\def\ee{\end{equation}}
\def\bt{\begin{theorem}}
\def\bcd{\begin{condition}}
\def\ecd{\end{condition}}
\def\et{\end{theorem}}
\def\bc{\begin{corollary}}
\def\ec{\end{corollary}}
\def\bde{\begin{definition}}
\def\ede{\end{definition}}
\def\bl{\begin{lemma}}
\def\el{\end{lemma}}
\def\bp{\begin{proposition}}
\def\ep{\end{proposition}}
\def\bex{\begin{example}}
\def\eex{\end{example}}
\def\br{\begin{remark}}
\def\er{\end{remark}}
\def\ba{\begin{array}}
\def\ea{\end{array}}
\def\ed{\end{document}}

\def\square#1{\vbox{\hrule\hbox{\vrule height#1%
     \kern#1\vrule}\hrule}}
\def\rectangle#1#2{\vbox{\hrule\hbox{\vrule height#1%
     \kern#2\vrule}\hrule}}

\font\tenbb=msbm10 \font\sevenbb=msbm7 \font\fivebb=msbm5

\newfam\bbfam
\scriptscriptfont\bbfam=\fivebb \textfont\bbfam=\tenbb
\scriptfont\bbfam=\sevenbb

\newtheorem{lemma}{Lemma}[section]
\newtheorem{remark}{Remark}[section]
\newtheorem{example}{Example}[section]
\newtheorem{theorem}{Theorem}[section]
\newtheorem{corollary}{Corollary}[section]

\newtheorem{definition}{Definition}[section]
\newtheorem{proposition}{Proposition}[section]
\newtheorem{condition}{Condition}[section]

\makeatletter
   
   \@addtoreset{equation}{section}
\makeatother

\begin{document}

\title{Path-depending controlled mean-field coupled forward-backward SDEs. The associated stochastic maximum principle\footnotemark[1]
%\footnote{the work}
}
\author{Rainer Buckdahn$^{1,2}$,\,\, Juan Li$^{3,2,\dag}$,\,\, Junsong Li$^{3,\dag}$,\,\, Chuanzhi Xing$^{2,\dag}$\\
{\small $^1$Laboratoire de Math\'{e}matiques de Bretagne Atlantique, Univ Brest, }\\
{\small UMR CNRS 6205, 6 avenue Le Gorgeu, 29200 Brest, France.}\\
{\small $^2$Research Center for Mathematics and Interdisciplinary Sciences,}\\
{\small         Shandong University, Qingdao 266237, P.~R.~China.}\\
%{\small $^3$ School of Mathematics, Shandong University, Jinan 250100, PR China}\\
{\small $^3$School of Mathematics and Statistics, Shandong University, Weihai, Weihai 264209, P.~R.~China.}\\
{\small {\it E-mails: rainer.buckdahn@univ-brest.fr; juanli@sdu.edu.cn; junsongli@mail.sdu.edu.cn; xingchuanzhi@sdu.edu.cn.}}\\
}

\renewcommand{\thefootnote}{\fnsymbol{footnote}}
\footnotetext[1]{Juan Li is supported by the NSF of P.R. China (NOs. 12031009, 11871037), National Key R and D Program of China (NO. 2018YFA0703900), NSFC-RS (No. 11661130148; NA150344).}
\footnotetext[2]{Corresponding authors.}

\date{July 26, 2023}
\maketitle
\textbf{Abstract}.
In the present paper we discuss a new type of mean-field coupled forward-backward stochastic differential equations (MFFBSDEs). The novelty consists in the fact that the coefficients of both the forward as well as the backward SDEs depend not only on the controlled solution processes $(X_t,Y_t,Z_t)$ at the current time $t$, but also on the law of the paths of $(X,Y,u)$ of the solution process and the control process.
The existence of the solution for such a MFFBSDE which is fully coupled through the law of the paths of $(X,Y)$ in the coefficients of both the forward and the backward equations is proved under rather general assumptions. Concerning the law, we just suppose the continuity under the 2-Wasserstein distance of the coefficients with respect to the law of $(X,Y)$. The uniqueness is shown under Lipschitz assumptions and the non anticipativity of the law of $X$ in the forward equation. The main part of the work is devoted to the study of Pontryagin's maximal principle for such a MFFBSDE. The dependence of the coefficients on the law of the paths of the solution processes and their control makes that a completely new and interesting criterion for the optimality of a stochastic control for the MFFBSDE is obtained. In particular, also the Hamiltonian is novel and quite different from that in the existing literature.
Last but not least,  under the assumption of convexity of the Hamiltonian we show that our optimality condition is not only necessary but also sufficient.


\section{Introduction}\label{FBSDEIntroduction}

\qquad $\,\,\,\,$With their seminal paper of 1990 Pardoux and Peng \cite{PP1} were the first to introduce and to investigate nonlinear classical backward stochastic differential equations (BSDEs). They proved the existence and the uniqueness of the solution for nonlinear BSDEs under Lipschitz assumptions. Since their pioneering work the theory of BSDEs has been developing very quickly and dynamically, and this also thanks to its numerous applications, for example, in stochastic control, finance and partial differential equations.

The Hamiltonian systems, derived from Pontryagin's maximum principle for stochastic optimal control problems, have turned out to be forward-backward stochastic differential equations (FBSDEs). Tightly related with the theory of BSDEs and stimulated by their different aplications, namely in stochastic control, also the theory of fully coupled FBSDEs has experienced a rapid and dynamic development and attracted a lot of researchers. Antonelli \cite{A} was the first to prove the existence and the uniqueness of the solution for fully-coupled FBSDEs driven by a Brownian motion. He did it with the help of a fixed point theorem for a sufficiently small time interval. Afterwards
several authors developed different other methods
to study fully-coupled FBSDEs, but now for an arbitrarily fixed time horizon. There are mainly three
methods. One is the so-called ``four-step scheme" (see Ma, Protter, Yong \cite{MPY}) which combines PDE and
probability methods. The authors proved the existence and the uniqueness of the solution for fully-coupled FBSDEs on
an arbitrarily given time interval, but they supposed the diffusion coefficient to be non-degenerate and deterministic. Another method is the purely probabilistic continuation method; we refer to Hu and Peng \cite{HP}, Pardoux and Tang \cite{PT}, and other subsequent works. Last but not least, a third method is inspired by numerical approaches for some linear FBSDEs (see Delarue \cite{D} and Zhang \cite{Z}).

For stochastic control problems both the stochastic maximum principle (SMP) as well as the Bellman dynamic programming principle (DPP) are the most important approaches, in order to get an optimal control strategy. For the use of Bellman's DPP we refer, for instance, to \cite{LW2015} and the references therein.

As concerns the SMP methods, we can distinguish essentiallly two different types. On the one hand, we have Pontryagin's SMP which restricts to the use of the first order derivatives of the coefficients but supposes the convexity of the control state space, and on the other hand, there is Peng's SMP which uses with repsect to the state variables the derivatives up to the second order, but by making use of spike controls it doesn't need the convexity of the state control space anymore.
There have been a lot of recent works on the stochastic maximum principle. For instance, Wu \cite{W3} discussed the SMP for the fully-coupled FBSDEs with a convex control state space. Hu \cite{Hu2017} studied recursive stochastic optimal control problems with a control state space which does not need to be convex, and the driving coefficient of the backward stochastic differential equation contains the $Z$-part of the solution. Zhang, Li and Xiong \cite{ZhangLX2020} investigated  a partially-observed optimal control problem in which the dynamics of the controlled state process is a stochastic differential equation with delay.

Recent works have extended the studies of stochastic control problems and those of the SMP to those involving  mean-field stochastic differential equations (MFSDEs), also known as McKean-Vlasov equations. Mean-field backward stochastic differential equations (mean-field BSDEs) were introduced by \cite{BLP1}; for more properties of mean-field BSDEs we refer to \cite{BDLP2009}. Since then several works have been devoted to the investigation of the SMP involving stochastic equations of mean-field type. So Buckdahn, Li and Ma \cite{BLM} studied the optimal control problem for a class of general mean-field stochastic differential equations in which the coefficients depend nonlinearly on both the state process and its law.
\iffalse
Li \cite{L2018} considered general mean-field FBSDEs with jumps whose generator depends also on the law of the solution. She also studied the associated integral partial differential equation. Without being exhaustive, let us also mention the works by
Hao and Li \cite{HL2014}, Buckdahn, Djehiche, Li~\cite{BDL} and Li~\cite{L}.
\fi

In a lot of the papers which study the SMP in mean-field context, the coefficients of the dynamics of the controlled state process and of the cost functional depend on the expectation of the solution. In the general case studied in the existing literature the coefficients depend on the joint law of the controlled state process and the governing control process both at the current time $t$, or, when the stochastic dynamics are defined through a MFFBSDE, the coefficients at time $t$ depend on the solution $(X_t,Y_t,Z_t)$ as well as on its joint law $P_{(X_t,Y_t,u_t)}$ of the solution and the control process. The novelty of our work is that we don't restrict in our investigations to this ``Markovian" case, we consider the coefficients depending on the solution processes at the current time $t$, but also on the joint law $P_{(X,Y,u)}$ of the trajectories of the solutions processes and the control process. This leads us not only to a new, more general type of MFFBSDEs but also to a new, very interesting criterion for the optimality of a stochastic control process.

To be more precise, in this paper we consider the following fully-coupled mean-field forward-backward stochastic differential equations (MFFBSDEs):
\begin{equation}\label{Eq1.1}
\left\{\begin{aligned}
dX_t&=\sigma(t,X_t,P_{(X,Y)})dB_t+b(t,X_t,P_{(X,Y)})dt,\ t\in[0,T],\\
dY_t&=-f(t,X_t,Y_t,Z_t,P_{(X,Y)})dt+Z_tdB_t,\ t \in [0,T],\\
X_0&=x\in\mathbb{R}^{d},\ Y_{T}=\Phi(X_{T},P_{(X,Y)}),\\
	\end{aligned}\right.
\end{equation}
where the time horizon $T>0$ is arbitrary but fixed and the coefficients
$b:[0,T]\times\mathbb{R}^{d}\times\mathcal{P}_{2}(\mathcal{C}_T^{2})\rightarrow\mathbb{R}^{d},\  \sigma:[0,T]\times\mathbb{R}^{d}\times\mathcal{P}_{2}(\mathcal{C}_T^{2})\rightarrow\mathbb{R}^{d\times d},\
 f:[0,T]\times\mathbb{R}^{d}\times\mathbb{R}\times\mathbb{R}^{d}\times\mathcal{P}_{2}(\mathcal{C}_T^{2})\rightarrow\mathbb{R},\
 \Phi:\mathbb{R}^{d}\times\mathcal{P}_{2}(\mathcal{C}_T^{2})\rightarrow\mathbb{R}$\ are deterministic and satisfy standard assumptions. In particular, for the proof of the existence of a solution we only need for the law variable of the coefficients the continuity in the sense of the 2-Wasserstein distance.
After having stated the existence of a solution of MFBSDE (\ref{Eq1.1}) (Theorem 3.1) we prove the uniqueness (Theorem 3.2) under Lipschitz assumptions on the coefficients and the condition that at any time $t$ the law in the coefficients of the forward equation only takes into account the trajectory of $X$ until $t$, but doesn't depend on $Y$, while in the BSDE at time $t$ the law $P_{(X,Y_{\cdot\vee t})}$ can still depend on the full trajectory of $X$ but concerning $Y$ only on its future trajectory $Y_{\cdot\vee t}$. For simplicity of the presentation we restrict in our paper to the case where the forward SDE in (\ref{Eq1.1}) has only the stochastic integral but no drift. However, let us point out that the case including the drift can be discussed with exactly the same arguments as those which we use in our work.
After having studied the existence and the uniqueness of the solution for our MFBSDE we introduce our stochastic control problem. However, the techniques we use to obtain the existence of a solution for MFFBSDE (\ref{Eq1.1}) involve the Malliavin calculus in order to get a suitable uniform estimate for $Z$. This implies that the control process $u$ cannot appear as a direct variable of the coefficients of our MFFBSDE, unless we restrict to control processes which are smooth in the sense of Malliavin calculus. But as this seems in our eyes to be too restrictive, we prefer to have the dependence on the control only through their joint law with the solution process. This makes that the dynamics of our stochastic control problem takes the following form:
\begin{equation}\label{Eqcontrol_0}
	\left\{\begin{aligned}
		dX^{v}_t&=\sigma(t,X_t,P_{(X^{v}_{\cdot\wedge t},v)})dB_t,\ X^{v}_{0}=x\in\mathbb{R},\\
		dY^{v}_t&=-f(t,X^{v}_t,Y^{v}_t,Z^{v}_t,P_{(X^{v},Y^{v}_{t\vee\cdot},v)})dt+Z^{v}_tdB_t,\\
		Y^{v}_T&=\Phi(X^{v}_T,P_{(X^{v},v)}),
	\end{aligned}\right.
\end{equation}
where $v$ runs the space of admissible stochastic control processes.
The arguments of the Theorems \ref{Thm3.1} and \ref{Theorem_3.2} allow to get the existence and the uniqueness also for (\ref{Eqcontrol_0}).

The main part of our manuscript is devoted to the computations on the SMP for the stochastic control problem whose dynamics are defined through (\ref{Eqcontrol_0}) and whose cost functional is given by

\begin{equation}\label{Eqcontrol_0A}
	\begin{aligned}
		J(v)=E[ & \displaystyle \int_0^T L(t,X^{v}_{t},Y^{v}_{t},Z^{v}_{t},P_{(X^{v},Y^{v},v)})dt+\varphi(X^{v}_{T},P_{(X^{v},Y^{v},v)})].
	\end{aligned}
\end{equation}

The necessary optimality criterion for a control process is stated in Theorem \ref{Thm5.1};
Theorem \ref{Thmsufficient} shows that this condition combined with a convexity assumption is also a sufficient one. However, the associated Hamiltonian $H$ cannot be the classical one

$H(t,x,y,z,p,k,\mu,\nu)=-f(t,x,y,z,\mu)p+\sigma(t,x,\nu)k+L(t,x,y,z,\mu),$

$\qquad(t,x,y,z,p;k,\mu,\nu)\in [0,T]\times \mathbb{R}^4\times\mathcal{P}_2(C([0,T];\mathbb{R}^2)\times L^2([0,T]))\times \mathcal{P}_2(C([0,T];\mathbb{R})\times L^2([0,T]))$,

\noindent anymore. Our new Hamiltonian has to consider two effects:
Let $u$ be an optimal control process and define for the associate dynamics $(X,Y,Z)$ the probability measures $\nu=P_{(X,u)}, \mu= P_{(X,Y,u)}$. The terminal value $Y_{T}=\Phi(X_{T},P_{(X,u)})$ and also $\varphi(X_{T},P_{(X,Y,u)})$ depend on the law of the whole path $(X,u)$ and $(X,Y,u)$, respectively. This has as consequence that they produce their own time-dependent adapted coefficients
$(\partial_{\mu}\Phi)_{j}^{*}(t)[p(T)]$ and $(\partial_{\mu}\varphi)^{*}_{i}(t),\ j=1,2,\ i=1,2,3$, which we have to
take into account in the definition of the Hamiltonian and its derivatives. It adds that, like $(\partial_{\mu}\Phi)_{j}^{*}(t)[p(T)]$, the derivatives with respect to the measure $(\partial_\mu\sigma)_1^*(t)[k],\, (\partial_\mu f)_i^*(t)[p],\, i=1,2,$ are linear functionals but now of the whole solution process $(p,(q,k))$ of the adjoint forward-backward SDE, and don't depend on $p$ and $k$ only in a multiplicative way.
This makes that our Hamiltonian cannot have the classical form. We define the Hamiltonian just as the following vector function:
\begin{equation}\label{Eq_Hamiltonian_00}
	\begin{aligned}
		H(t,x,y,z,\nu,\mu):=(-f(t,x,y,z,\mu),\ \sigma(t,x,\nu),\ L(t,x,y,z,\mu),-\Phi(\cdot,\nu),\varphi(\cdot,\mu)).
	\end{aligned}
\end{equation}
Our paper is organised as follows: In Section \ref{FBSDESec2} we introduce the framework of our studies and some useful notations.
Section \ref{FBSDESec3} is devoted to the proof of the existence of the solution for our general MFFBSDE \eqref{Eq1.1} by using the Malliavin calculus, but also to the proof of the uniquenessof the solution for the MFFBSDE. In Section \ref{FBSDESec4} we introduce the derivatives with respect to the law over suitable function spaces.
In Section \ref{FBSDESec5} we study the stochastic maximum principle for MFFBSDEs with control, and obtain a necessary condition for the optimality of a stochastic control process  (Theorem \ref{Thm5.1}). For this end, we investigate with by far non trivial computations in particular the associated variational equation, which turns out to be on its turn also a MFFBSDE (Lemmas \ref{lem5.1}-\ref{lem5.2}), but also the adjoint FBSDE, which is composed of an affine mean-field forward SDE with delay and an affine mean-field BSDE with anticipation (Lemma \ref{Lemma5.5}). The discussion of these equations leads us to the definition of our Hamiltonian and the computation of its derivatives. Finally, Theorem \ref{theo_variin} with its variational inequality leads to the SMP established in Theorem \ref{Thm5.1}, and Theorem
\ref{Thmsufficient} gives a sufficient optimality condition.

\section{Preliminaries}\label{FBSDESec2}

\qquad $\,\,\,\,$ Let $(\Omega,\mathcal{F},P)$ be a complete probability space endowed with a $d$-dimensional Brownian motion $B$. Let $\mathbb{F}=({\mathcal{F}_t})_{t\geq0}$ be the natural filtration generated by $B$ and augmented by all $P$-null sets (that is, $\mathcal{F}_t=\sigma \{B_s:0\leq s \leq t \}\vee\mathcal{N}_P, \ t\geq 0$, where $\mathcal{N}_P$ is the set of all $P$-null subsets). Moreover, let $T>0$ be a fixed time horizon.

We assume that there is a sub-$\sigma$-field $\mathcal{F}_0\subset\mathcal{F}$ including all $P$-null subsets of $\mathcal{F}$, such that
\begin{itemize}
\item[(i)] the Brownian motion $B$ is independent of $\mathcal{F}_0$;

\item[(ii)] $\mathcal{F}_0$ is ``rich enough'', i.e., $\mathcal{P}_{2}(\mathbb{R}^k)=\{P_\xi, \xi\in L^2(\mathcal{F}_0;\mathbb{R}^k)\},$ $k\geq 1$, where $P_\xi := P\circ \xi^{-1}$ denotes the law of the random variable $\xi$.
\end{itemize}
For given $T>0$ we denote by $\mathcal{C}_{T}:=\mathcal{C}([0,T];\mathbb{R}^{k})$ the space of $\mathbb{R}^{k}$-valued continuous functions on $[0, T ]$ endowed with the supremum norm $|\cdot|_{\mathcal{C}_T}$ where  $|\phi|_{\mathcal{C}_{T}}:=\sup_{t\in[0,T]}|\phi(t)|$, for $\phi\in\mathcal{C}_{T}$, and by $\mathcal{B}(\mathcal{C}_T)$ its topological $\sigma$-field.

 Furthermore, we let $\mathcal{P}(\Omega)$ denote the space of all probability measures on $(\Omega, \mathcal{F})$.
 Consider now the space of all probability measures on $(\mathcal{C}_T, \mathcal{B}(\mathcal{C}_T))$, denoted by $\mathcal{P}(\mathcal{C}_T)$, and for $p \geq 1$ we consider $\mathcal{P}_p(\mathcal{C}_T) (\subseteq \mathcal{P}(\mathcal{C}_T))$ the space of all probabilities with finite $p$-th moment. We recall that the $p$-Wasserstein metric on $\mathcal{P}_p(\mathcal{C}_T)$ is defined as a mapping $W_p: \mathcal{P}_p(\mathcal{C}_T) \times \mathcal{P}_p(\mathcal{C}_T) \mapsto \mathbb{R}_{+}$ such that, for all $\mu, \nu \in \mathcal{P}_p(\mathcal{C}_T)$,
\begin{equation*}
\begin{aligned}
 W_p(\mu, \nu):= \inf\{(\int_{\mathcal{C}_T^{2}}|x-y|_{\mathcal{C}_T}^p \pi(d x, d y))^{\frac{1}{p}}: \pi \in \mathcal{P}_p(\mathcal{C}_T^{2}) \textnormal{ with marginals }\mu\ \textnormal{and}\ \nu\}.
\end{aligned}
\end{equation*}
%Here $|\cdot|_{\mathcal{C}_{T}}$ stands for the supremum norm on $\mathcal{C}_{T}$.
For $p=1$, recall that
due to the Kantorovich-Rubinstein Theorem (refer to \cite{CP2013} or \cite{DEdwards2011}),
\begin{equation*}
  \begin{aligned}
\displaystyle W_{1}(\mu,\mu')
    &=||\mu-\mu'||_{1}\\
    &=\sup \big\{ | \int_{\mathcal{C}_T^{2}} h d\mu-\int_{\mathcal{C}_T^{2}} h d\mu' |,\ h\in Lip_{1}({\mathcal{C}_T^{2}}),\ h(0)=0.\big\},\ \mu,\ \mu'\in\mathcal{P}_{1}(\mathcal{C}^{2}_T).
  \end{aligned}
\end{equation*}
Here $Lip_{1}(\mathcal{C}_{T}^{2})$ denotes the space of all Lipschitz functions $\phi:\ \mathcal{C}_{T}^{2}\rightarrow\mathbb{R}$ with Lipschitz constant 1.
In this paper we shall use the 1-Wasserstein metric $W_1$ and the 2-Wasserstein metric $W_2$, and we abbreviate $(\mathcal{P}_2(\mathcal{C}_T), W_2)$ by $\mathcal{P}_2(\mathcal{C}_T)$. Since $\mathcal{C}_T$ is a separable Banach space, it is known that $\mathcal{P}_2(\mathcal{C}_T)$ is a separable and complete metric space. Furthermore, it is known that (refer to, e.g., \cite{V2003}), for $\mu_n, \mu \in \mathcal{P}_2(\mathcal{C}_T)$,
\begin{equation}
\begin{aligned}
\lim _{n \rightarrow \infty} W_2(\mu_n, \mu)=0 \Longleftrightarrow&\ \mu_n\rightharpoonup\mu\ \text{weakly}\ \text{and,\ as }\ N \rightarrow+\infty,\ \\
&\ \sup _n \int_{\Omega}\bigr|\varphi\bigr|_{\mathcal{C}_T}^2 I\{|\varphi|_{\mathcal{C}_T}\geq N\}\mu_n(d \varphi) \rightarrow 0.
\end{aligned}
\end{equation}
%\stackrel{w}{\rightarrow}

We will use the following spaces of stochastic processes. Recall that $|\cdot|$ denotes the Euclidean norm in the corresponding Euclidean space.
%$\langle\cdot,\cdot\rangle$ be the inner product in the corresponding Hilbert space.
\begin{itemize}
\item $L^{2}(\mathcal{F}; \mathbb{R}^d)$ is the set of $\mathbb{R}^{d}$-valued, $\mathcal{F}$-measurable random variables $\xi$ such that $|\xi|_{L^{2}}=(E[|\xi|^{2}])^{\frac{1}{2}}<\infty$.
\item $S^2_{\mathbb{F}}:= S^2_{\mathbb{F}}(0,T; \mathbb{R}^d)$ is the set of $\mathbb{R}^{d}$-valued, $\mathbb{F}$-adapted continuous processes $Y$ with $\displaystyle E[\sup_{0\leq t\leq T}|Y_t|^2]<\infty.$
\item $M^2_{\mathbb{F}}:= M^2_{\mathbb{F}}(0,T; {\mathbb{R}^{d}})$ is the set of $\mathbb{R}^{d}$-valued $\mathbb{F}$-progressively measurable processes $Z$ on $[0,T]$ with $\displaystyle E[\int_0^T|Z_t|^2dt]<\infty.$
%\item  $\mathcal{P}_{1}(\mathcal{C}_T^{2k})$ is the collection of all probability measures with finite first moments over $\mathcal{C}_T^{2k}=\mathcal{C}([0,T]; \mathbb{R}^{k})\times\mathcal{C}([0,T]; \mathbb{R}^{k})$, endowed with the $1$-Wasserstein metric, where $\mathcal{C}^{k}([0,T]; \mathbb{R}^{k})$ denote the space of $\mathbb{R}^k$-valued continuous functions on $[0, T ]$ and for $\mu,\mu'\in\mathcal{P}_{1}(\mathcal{C}_T^{2k})$,
%\begin{equation*}\label{Eq2.1}
%\begin{aligned}
%\displaystyle W_1(\mu,\mu'):= &\inf\big\{\int_{\mathcal{C}_T^{2k}\times\mathcal{C}_T^{2k}}|x-x'|\rho(dx,dx'):\rho\in\mathcal{P}_{1}(\mathcal{C}_T^{2k}\times\mathcal{C}_T^{2k}), \rho(\cdot\times\mathcal{C}_T^{2k})\\
%&=\mu,\ \rho(\mathcal{C}_T^{2k}\times\cdot)=\mu' \big\}.
%\end{aligned}
%\end{equation*}
\end{itemize}

We now give our assumptions for $\sigma:[0,T]\times\mathbb{R}^{d}\times\mathcal{P}_{2}(\mathcal{C}_T^{2})\rightarrow\mathbb{R}^{d\times d},
 f:[0,T]\times\mathbb{R}^{d}\times\mathbb{R}\times\mathbb{R}^{d}\times\mathcal{P}_{2}(\mathcal{C}_T^{2})\rightarrow\mathbb{R}$ and $\Phi:\mathbb{R}^{d}\times\mathcal{P}_{2}(\mathcal{C}_T^{2})\rightarrow\mathbb{R}:$
\begin{itemize}
\item[\bf{(H1)}] \ $f,\sigma$ are measurable and continuous with respect to $\mu\in \mathcal{P}_{1}(\mathcal{C}_T^{2})$, with a continuity modulus $\rho: \mathbb{R}_{+}\rightarrow\mathbb{R}_{+}$ with $\rho(0+)=0$, uniformly with respect to $t$.
\item[\bf{(H2)}] \ $\sigma,\ \Phi$ and\ $f$ are bounded and have bounded derivatives with respect to $x$ and to $(x,y,z)$, respectively.
\end{itemize}
Notice that from (H2) we get: $f$ and $\sigma$ are Lipschitz in $(x,y,z)$, i.e., there exists a constant $C\in\mathbb{R}_{+}$ such that, $P$-a.s., for all $(x, y,z),\ (x',y',z')\in\mathbb{R}^{d}\times\mathbb{R}^{d}\times\mathbb{R}^{d},\ t\in[0,T],\ \mu\in\mathcal{P}_{1}(\mathcal{C}^{2}_T)$,
\begin{equation*}
 \begin{aligned}
    |f(t,x,y,z,\mu)-f(t,x',y',z',\mu)|&\leq C(|x-x'|+|y-y'|+|z-z'|),\\
    |\sigma(t,x,\mu)-\sigma(t,x',\mu)|&\leq C|x-x'|.
    \end{aligned}
\end{equation*}

At the end of this section, we have still to introduce some notations. By $(\widetilde{\Omega}, \widetilde{\mathcal{F}}, \widetilde{P})$ we denote a copy of the probability space $(\Omega, \mathcal{F}, P)$. For any random variable (of arbitrary dimension) $\vartheta$ over $(\Omega,\ \mathcal{F}, P)$ we denote by $\widetilde{\vartheta}$ a copy over  $(\widetilde{\Omega}, \widetilde{\mathcal{F}}, \widetilde{P})$, i.e., a random variable such that $\widetilde{P}_{\widetilde{\vartheta}}=P_{\vartheta}$. The expectation $\widetilde{E}[\cdot]=\int_{\widetilde{\Omega}}(\cdot) d \widetilde{P}$ acts only over the variables endowed with a tilde. This can be made rigorous by working with the product space $(\Omega, \mathcal{F}, P) \otimes(\widetilde{\Omega}, \widetilde{\mathcal{F}}, \widetilde{P})$ $(=(\Omega, \mathcal{F}, P) \otimes(\Omega, \mathcal{F}, P)$ for the choice $(\widetilde{\Omega}, \widetilde{\mathcal{F}}, \widetilde{P})=(\Omega, \mathcal{F}, P))$, and by extending the random variable $\vartheta$ to $(\Omega, \mathcal{F}, P)\otimes(\widetilde{\Omega}, \widetilde{\mathcal{F}}, \widetilde{P}), \vartheta(\omega,\widetilde{\omega})=\vartheta(\omega)$, and putting $\widetilde{\vartheta}(\omega,\widetilde{\omega})=\widetilde{\vartheta}(\widetilde{\omega})(=\vartheta(\omega)), (\omega, \widetilde{\omega})\in\Omega\otimes\widetilde{\Omega}$.
\begin{remark}\label{rem2.1}
In order not to overcomplicate the already heavy presentation
of this paper, in what follows we shall assume all processes to be 1-dimensional (i.e.,
$d = k = 1$). We emphasize that the higher dimensional cases can be handled in the same way, without substantial difficulties, except for the price of even heavier notations.
\end{remark}

\section{Mean-field FBSDE: Existence and uniqueness }\label{FBSDESec3}
We consider the following mean-field coupled forward-backward stochastic differential equation:
\begin{equation}\label{EqFBSDE3.1}
\left\{ \begin{aligned}
dX_t&=\sigma(t,X_t,P_{(X,Y)})dB_t,\ t \in [0,T],\\
dY_t&=-f(t,X_t,Y_t,Z_t,P_{(X,Y)})dt+Z_tdB_t,\ t \in [0,T],\\
X_0&=x \in \mathbb{R}^d,\ Y_{T}=\Phi(X_{T},P_{(X,Y)}),
\end{aligned}\right. \end{equation}
where $\sigma, f$ and $\Phi$ are introduced above.
%For future research, we first give the existence and uniqueness of solutions of FBSDE $(\ref{FBSDE1})$. In this section, %we give the existence and uniqueness condition of the solution of the equation $(\ref{FBSDE1})$ and give detailed %proof.\\
\begin{theorem}\label{Thm3.1}
We assume \textnormal{(H1)} and \textnormal{(H2)} hold true. Then mean-field FBSDE $(\ref{EqFBSDE3.1})$ has an adapted solution $(X,Y,Z)\in S^2_{\mathbb{F}}\times S^2_{\mathbb{F}}\times M^2_{\mathbb{F}}$.
\end{theorem}

\begin{proof}
 The key for the proof of the existence of a solution is an application of Schauder's fixed point theorem stating that if $V$ is a Hausdorff topological vector space, $\mathcal{K}\subset V$ is a nonempty convex closed subset, and $\mathcal{T}: \mathcal{K}\rightarrow \mathcal{K}$ is a continuous mapping such that $\mathcal{T}(\mathcal{K})\subset \mathcal{K}$ is contained in a compact subset of $\mathcal{K}$, then there exists $\mu\in \mathcal{K}$ such that $\mathcal{T}(\mu)=\mu$.

Let $\mu\in\mathcal{P}_1(\mathcal{C}_T^{2})$, and we consider
\begin{equation}\label{Eq3.2}
 \left\{ \begin{aligned}
 X_t^\mu &= x+\int_0^t\sigma(s, X_s^\mu, \mu)dB_s, \\
  Y_t^\mu &= \Phi(X_T^\mu, \mu)+\int_t^Tf(s, X_s^\mu, Y_s^\mu, Z_s^\mu, \mu)ds-\int_t^TZ_s^\mu dB_s,\ t\in[0, T].\\
\end{aligned}\right.
\end{equation}
Then, under the assumptions (H1) and (H2), \eqref{Eq3.2} has a unique solution $(X^{\mu}, Y^{\mu}, Z^{\mu})\in S^2_{\mathbb{F}}\times S^2_{\mathbb{F}}\times M^{2}_{\mathbb{F}}$.

Let $\theta_s^\mu :=(X_s^\mu, Y_s^\mu, Z_s^\mu),\ s\in [0,T]$. Taking the Malliavin derivatives (refer to Nualart \cite{N2006}), for $r<t$,
%Malliavin导数参考文献
we have, $drdP$-a.e.,
\begin{eqnarray*}
  D_r[X_t^\mu] = \sigma(r, X_r^\mu, \mu)+\int_r^t(\partial_x\sigma)(s, X_s^\mu, \mu)D_r[X_s^\mu]dB_s,\ t\in[r, T].
\end{eqnarray*}
Then, obviously,
\begin{equation}
D_r[X_t^\mu] = \sigma(r, X_r^\mu, \mu)\cdot \exp\Bigr\{\int_r^t(\partial_x\sigma)(s, X_s^\mu, \mu)dB_s-\frac{1}{2}\int_r^t|(\partial_x\sigma)(s, X_s^\mu, \mu)|^2ds\Bigr\},\ 0\leq r\leq t\leq T.
\end{equation}
From our assumptions it follows that, for $p\geq 1$, there is a constant $C_p$ such that,
\begin{equation}\label{Eq3_3.4}
E\Bigr[\big|\sup_{t\in[r,T]}D_r[X_t^\mu]\big|^p\Bigr|\mathcal{F}_r\Bigr]\leq C_p,\ t\in[0, T],\ 0\leq r \leq\ T.
\end{equation}
Moreover, $drdP$-a.e.,
\begin{equation}
\begin{aligned}
D_r[Y_t^\mu]=&\ (\partial_x\Phi)(X_T^\mu, \mu)D_r[X_T^\mu]+\int_t^T\{(\partial_xf)(s, \theta_s^\mu, \mu)D_r[X_s^\mu]+(\partial_yf)(s, \theta_s^\mu, \mu)D_r[Y_s^\mu]\\
&\ +(\partial_zf)(s, \theta_s^\mu, \mu)D_r[Z_s^\mu]\}ds-\int_t^TD_r[Z_s^\mu]dB_s,\ t\in[r, T].
\end{aligned}
\end{equation}
Putting
$$
\alpha_{\nu,t}^\mu:=\exp\Bigr\{\int_\nu^t(\partial_{y}f)(s, \theta_s^\mu, \mu)ds\Bigr\},\ 0\leq \nu\leq t\leq T,
$$
and observing that, for some $C\in\mathbb{R}_{+},\ |\alpha_{\nu,t}^\mu|\leq C,\ 0\leq \nu\leq t\leq T$,\ we have,  $drdP$-a.e.,
\begin{equation}
\begin{aligned}
D_r[Y_t^\mu]=&\ \alpha_{t,T}^\mu\cdot(\partial_x\Phi)(X_T^\mu, \mu)D_r[X_T^\mu]+\int_t^T(\partial_xf)(s, \theta_s^\mu, \mu)\alpha_{t,s}^\mu D_r[X_s^\mu]ds\\
&\ -\int_t^T\alpha_{t,s}^\mu\cdot D_r[Z_s^\mu]\cdot(dB_s-(\partial_zf)(s, \theta_s^\mu, \mu)ds),\ r \leq t \leq T.
\end{aligned}
\end{equation}
Let us introduce
\begin{equation*}
\begin{aligned}
L_t^\mu:=&\exp\Bigr\{\int_0^t(\partial_zf)(s, \theta_s^\mu, \mu)dB_s-\frac{1}{2}\int_0^t|(\partial_zf)(s, \theta_s^\mu, \mu)|^2ds\Bigr\},\ 0\leq t\leq T,\\
B_t^\mu:=& B_t-\int_{0}^{t}(\partial_{z}f)(s, \theta_s^\mu, \mu)ds,\ 0\leq t\leq T,\
P^\mu:= L_T^\mu\cdot P.
\end{aligned}
\end{equation*}
Thanks to Girsanov's theorem we see that
$B^{\mu}=(B^{\mu}_{t})_{t\in[0,T]}$ is an $(\mathbb{F}^B, P^{\mu})$-Brownian motion, and
\begin{eqnarray}
D_r[Y_t^\mu] = E^{\mu}\[\alpha_{t,T}^\mu\cdot(\partial_x\Phi)(X_T^\mu, \mu)D_r[X_T^\mu]+\int_t^T\alpha_{t,s}^\mu\cdot(\partial_xf)(s, \theta_s^\mu, \mu)\cdot D_r[X_s^\mu]ds\bigr|\mathcal{F}_t\],
\end{eqnarray}
$0 \leq r\leq t \leq T$. Hence,
\begin{eqnarray}
|D_r[Y_t^\mu]| \leq CE\[\frac{L_T^\mu}{L_t^\mu}(|D_r[X_T^\mu]|+\int_t^T|D_r[X_s^\mu]|ds)|\mathcal{F}_t\].
\end{eqnarray}
Thus, taking the limit in probability, as $r \leq\ t\downarrow r$, thanks to \eqref{Eq3_3.4} we have
\begin{equation}
\begin{aligned}
|Z_r^\mu| &\leq CE\[\frac{L_T^\mu}{L_r^\mu}(|D_r[X_T^\mu]|+\int_r^T|D_r[X_s^\mu]|ds)\big|\mathcal{F}_r\]\\
&\leq C\(E\[(\frac{L_T^\mu}{L_r^\mu})^2\big|\mathcal{F}_r\]\)^{\frac{1}{2}}\(E\[|D_r[X_T^\mu]|^2+\int_r^T|D_r[X_s^\mu]|^2ds\big|\mathcal{F}_r\]\)^{\frac{1}{2}}\\
&\leq C, \ \ drdP\text{-a.e.},\ \mu\in\mathcal{P}_1(\mathcal{C}_T^{2}).
\end{aligned}
\end{equation}
Consequently,
\begin{equation}\label{Eq3.10Z}
|Z_r^\mu|\leq C,\ drdP\text{-a.e.},\ \text{for any}\ \mu\in\mathcal{P}_1(\mathcal{C}_T^{2}).
\end{equation}
%Now let's verify that the conditions of the fixed point theorem are satisfied.
We also observe that, as $f$ and $\Phi$ are bounded, from a standard BSDE estimate it follows that, for some $C\in\mathbb{R}_{+},\ |Y_{t}^{\mu}| \leq C,\ t\in [0,T], P$-a.s.
Let $\mathcal {M}_1(\mathcal{C}_T^{2})$ denote the normed space
$$\mathcal {M}_1(\mathcal{C}_T^{2}):=\Bigr\{\gamma \rm\ signed\ measure\ over\ (\mathcal{C}_T^{2}, \mathcal {B}(\mathcal{C}_T^{2}))\ \Bigr|\int_{\mathcal{C}_T^{2}}\bigr|\phi\bigr|_{\mathcal{C}_T^{2}}\cdot|\gamma|(d\phi)<+\infty \Bigr\},$$
endowed with the norm:
$$
\parallel \gamma\parallel_1:=\sup\Bigr\{\bigr|\int_{\mathcal{C}_T^2}hd\gamma\bigr|\Bigr|\ h\in Lip_1(\mathcal{C}_T\times\mathcal{C}_T),\ h(0)=0\Bigr\},\,\, \gamma\in \mathcal {M}_1(\mathcal{C}_T^{2}).
$$
Let $\mathscr{C}=(\mathscr{C}_{t})$ be the coordinate process on $\mathcal{C}_T^2: \mathscr{C}_t(\phi)=\phi_t,\  \phi\in\mathcal{C}_T^2,\ t\in[0, T]$. We consider the set
$$\mathcal {K}:=\Bigr\{\mu\in\mathcal{P}_1(\mathcal{C}_T^2)\Bigr|\ \int_{\mathcal{C}_T^2}\sup_{t\in[0,T]}|\mathscr{C}_t|^4d\mu\leq C;\ \int_{\mathcal{C}_T^2}|\mathscr{C}_t-\mathscr{C}_s|^4d\mu\leq C|t-s|^2,\ t,\ s\ \in[0, T]\Bigr\},$$
where $C\in \mathbb{R}_{+}$ will be specified later.\\
Recall that, by the Kantorovich-Rubinstein duality,
\begin{equation}
\begin{aligned}
W_1(\mu, \mu^1) =&\ \sup\Bigr\{\ \bigr|\int_{\mathcal{C}_T^2}hd\mu-\int_{\mathcal{C}_T^2}hd\mu^1\bigr|\Bigr|\ h\in Lip_1(\mathcal{C}_T^2), h(0)=0\Bigr\}\\
=&\ \parallel\mu-\mu^1\parallel_{1},\ \mu,\ \mu^1\in\mathcal{P}_1(\mathcal{C}_T^2).
\end{aligned}
\end{equation}
%from where, in particular, it follows
%$$\mu-\mu^1\in\mathcal {M}_1(\mathcal{C}_T^2).$$

\noindent{\bf{ Step 1}}.\ $\mathcal {K}$ is a convex and compact subset of $(\mathcal{P}_1(\mathcal{C}_T^2), W_1(\cdot, \cdot))$.\\
We know from the definition that $\mathcal {K}\subset \mathcal{P}_1(\mathcal{C}_T^2)$ is convex. Let us show that $\mathcal {K}$ is a compact subset in $(\mathcal{P}_1(\mathcal{C}_T^2), W_1(\cdot, \cdot))$, for some $C\in \mathbb{R}_{+}$ in the definition of $\mathcal {K}$ which is large enough.

Let $(\mu^n)_{n\geq 1}\subset\mathcal{K}$. Then $(\mu^n)_{n\geq 1}$ is tight over $(\mathcal{C}_T^2, \mathcal{B}(\mathcal{C}_T^2))$. Indeed, recall that the Arzel\`{a}-Ascoli Theorem provides a characterization of a compact subset of $\mathcal{C}_T^2$. Moreover, as
\begin{equation}
\mathop{\sup}\limits_{n\geq 1} \int_{\mathcal{C}_T^2}|\mathscr{C}|_{\mathcal{C}_{T}^{2}} I_{\{|\mathscr{C}|_{\mathcal{C}_{T}^{2}}\geq N\}}d\mu^n\leq \frac{1}{N^3}\mathop{\sup}\limits_{n\geq 1}\int_{\mathcal{C}_T^2}|\mathscr{C}|_{\mathcal{C}_{T}^{2}}^4d\mu^n\leq \frac{C}{N^3} \rightarrow 0,\ \text{as}\ N\rightarrow\infty,
\end{equation}
we get the existence of subsequence
$(\mu^{n_k})_{k\geq 1}\subset(\mu^{n})_{n\geq 1},$ and a probability measure $\mu\in\mathcal{P}_1(\mathcal{C}_T^2) \ \mbox{such that}\\ W_1(\mu^{n_k},\mu)\rightarrow 0, k\rightarrow\infty.$
Finally, from Mazur's Lemma we get $\mu\in \mathcal {K}$.  Thus $\mathcal {K}$ is a compact subset.

\noindent {\bf{Step 2}}. $\mathcal{T}$ defined by $\mathcal{T}(\mu):= P_{(X^\mu, Y^\mu)},\ \mu\in\mathcal{P}_{1}(\mathcal{C}_T^2)$, maps $\mathcal{K}$ into $\mathcal{K}$.

Indeed, because of the boundedness assumption on the coefficients we get that for all $p\geq1$ that there exists $C_{p}\in\mathbb{R}_{+}$ such that $E[\sup_{t\in[0, T]}|X_t^\mu|^p]\leq C_p$, and
\begin{equation}\label{eq_SDEt}
      E[|X_t^\mu-X_s^\mu|^p]\leq C_p|t-s|^{\frac{p}{2}},\ t,\ s\in[0, T],\ \mu\in\mathcal{P}_{1}(\mathcal{C}_T^2).
\end{equation}
On the other hand, as $f$ is bounded and $|Z_t^\mu|\leq C$, $dtdP$-a.e.,\ $\mu\in \mathcal{P}_1(\mathcal{C}_T^2)$, we get $E[\sup_{t\in[0, T]}|Y_t^\mu|^p]\leq C_p$, and
\begin{equation}\label{eq_BSDEt}
\ E[|Y_t^\mu-Y_s^\mu|^p]\leq C_p|t-s|^{\frac{p}{2}},\ t,\ s\in[0, T],\ \mu\in\mathcal{P}_{1}(\mathcal{C}_T^2).
\end{equation}

Moreover, from Burkholder-Davis-Gundy's inequality,
\begin{equation*}
E\Bigr[\Bigr|\int_s^tZ_r^\mu dB_r\Bigr|^p\Bigr]\leq C_pE\Bigr[\Bigr(\int_s^t|Z_r^\mu|^2dr\Bigr)^{\frac{p}{2}}\Bigr]\leq C_p|t-s|^{\frac{p}{2}},\  t,\ s\in[0, T],\ \mu\in\mathcal{P}_{1}(\mathcal{C}_T^2).
\end{equation*}
Hence, for $p=4$ (Recall that $\mathscr{C}$ is the coordinate process on $\mathcal{C}_T^2$, and so $\mathscr{C}_{t}(X^\mu,Y^\mu)=(X_t^\mu,Y_t^\mu)$),
\begin{equation*}
\begin{aligned}
\int_{\mathcal{C}_T^2}|\mathscr{C}_t-\mathscr{C}_s|^4dP_{(X^\mu, Y^\mu)}
&=E[|(X_t^\mu, Y_t^\mu)-(X_s^\mu, Y_s^\mu)|^4]\leq CE[|X_t^\mu-X_s^\mu|^4]+CE[|Y_t^\mu-Y_s^\mu|^4]\\
&\leq C |t-s|^{2},\ t, s\in [0,T].
\end{aligned}
\end{equation*}

Then with the above estimates for $(X^{\mu},Y^{\mu}),\ \mu\in\mathcal{P}_1(\mathcal{C}_T^2)$, in particular \eqref{eq_SDEt} and \eqref{eq_BSDEt}, we see $\mathcal{T}: \mathcal{P}_1(\mathcal{C}_T^2)\rightarrow \mathcal {K}$ and, hence, consider $\mathcal {T}: \mathcal {K}\rightarrow\mathcal {K}.$\\
\noindent {\bf{Step 3}}. We have to show that $\mathcal {T}: \mathcal {K}\rightarrow\mathcal {K}$ is continuous.

Let $(\mu^n)_{n\geq 1}\subset\mathcal {K},\, \mu\in\mathcal {K}$, and put $(X^n, Y^n, Z^n):=(X^{\mu^n}, Y^{\mu^n}, Z^{\mu^n})$ and $(X, Y, Z):=(X^{\mu}, Y^{\mu}, Z^{\mu})$. We suppose that $W_1(\mu^n, \mu)\rightarrow 0,\,  n\rightarrow\infty$. We have to show that $W_1(\mathcal {T}(\mu^n), \mathcal{T}(\mu))\rightarrow 0$, as $n\rightarrow\infty$.\\
For this we consider the following associated mean-field FBSDE:
\begin{equation}\label{step3_1}
\left\{
\begin{aligned}
  dX_t^\gamma &= \sigma(t, X_t^\gamma, \mu)dB_t, \ X_0^\gamma=x; \\
  dY_t^\gamma &= -f(t, X_t^\gamma, Y_t^\gamma, Z_t^\gamma, \mu)dt+Z_t^\gamma dB_t, \ Y_T^\gamma=\Phi(X_T^\gamma, \mu),
\end{aligned}
\right.
\end{equation}
for $\gamma=\mu,\ \mu^{n},\ n\geq 1$, respectively. From \eqref{step3_1} we obtain, for all $n\geq1,\, p\ge 2$,
\begin{equation}
E\Bigr[\sup_{t\in[0, T]}|X_t^\mu-X_t^{\mu^n}|^p\Bigr]\leq C\rho(W_1(\mu, \mu^n))^p.
\end{equation}
In fact, as
%\begin{equation*}
%\begin{aligned}
%E[\sup_{t\in[0,T]}|X_t^\mu-X_t^{\mu^n}|^p]=&E[\sup_{t\in[0,T]}|\int_{0}^{t}(\sigma(s, X_s^\mu, \mu)-\sigma(s, %X_s^{\mu^{n}}, \mu))dB_s|^{p}]\\
%\leq& C_{p}E[\sup_{t\in[0,T]}|\int_{0}^{t}|\sigma(s, X_s^\mu, \mu)-\sigma(s, X_s^{\mu^{n}}, \mu)|^{2}ds|^{\frac{p}{2}}]\\
%\leq& C_{p}E[\int_{0}^{T}(|X_t^\mu-X_t^{\mu^{n}}|+\rho(W_1(\mu,\mu^{n})))^{2}dt|^{\frac{p}{2}}]\\
%\leq& C_{p}\int_{0}^{T}E[\sup_{s\in[0,t]}|X_t^\mu-X_t^{\mu^{n}}|^{p}]dt+C_{p}\rho(W_1(\mu,\mu^{n}))^{p}.
%\end{aligned}
%\end{equation*}
\begin{equation*}
\begin{aligned}
E\Bigr[\sup_{t\in[0,s]}|X_t^\mu-X_t^{\mu^n}|^p\Bigr]
\leq&\ C_{p}E\Bigr[\Bigr|\int_{0}^{s}|\sigma(s, X_s^\mu, \mu)-\sigma(s, X_s^{\mu^{n}}, \mu^n)|^{2}ds\Bigr|^{\frac{p}{2}}\Bigr]\\
%\leq& C_{p}E[\int_{0}^{T}(|X_t^\mu-X_t^{\mu^{n}}|+\rho(W_1(\mu,\mu^{n})))^{2}dt|^{\frac{p}{2}}]\\
\leq&\ C_{p}\int_{0}^{s}E\Bigr[\sup_{r\in[0,t]}|X_r^\mu-X_r^{\mu^{n}}|^{p}\Bigr]dt+C_{p}\rho(W_1(\mu,\mu^{n}))^{p},\ s\in[0,T],
\end{aligned}
\end{equation*}
Gronwall's inequality gives the stated result. Also for the BSDE we get by standard estimates
\begin{equation}
\begin{aligned}
E\Big[\sup_{t\in[0, T]}|Y_t^\mu-Y_t^{\mu^n}|^p+\(\int_t^T|Z_s^\mu-Z_s^{\mu^n}|^2ds\)^{\frac{p}{2}}\Big]\leq C_p\rho(W_1(\mu, \mu^n))^p.
\end{aligned}
\end{equation}
Thus,
\begin{equation}
\begin{aligned}
&W_{1}(\mathcal {T}(\mu), \mathcal{T}(\mu^n))=W_{1}(P_{(X^{\mu}, Y^{\mu})}, P_{(X^{\mu^n}, Y^{\mu^n})})\\
\leq\ &E\[\sup_{t\in[0, T]}|(X_t^\mu, Y_t^{\mu})-(X_t^{\mu^n}, Y_t^{\mu^n})|\]\\
\leq\ &C\rho(W_1(\mu, \mu^n))\rightarrow0,\ \text{as}\ n\rightarrow\infty.
\end{aligned}
\end{equation}
\noindent {\bf{Step 4}}. Application of Schauder's fixed point theorem.

$\mathcal {T}(\mathcal {K})$ is compact in $\mathcal {K}$, since $\mathcal {K}$ is compact and $\mathcal {T}: \mathcal {K}\rightarrow\mathcal {K}$ is continuous. Embedding $\mathcal {K}$ into the separable linear Hausdorff space ${M}_1(\mathcal{C}_T^2)$ allows to apply Schauder's fixed point theorem (for details we refer ro Li and Min \cite{LH2016}),
from where we conclude the existence of $\mu\in \mathcal {K}(\subset\  \mathcal{P}_1(\mathcal{C}_T^2))$ such that
$\mu=\mathcal {T}(\mu)=P_{(X^\mu, Y^\mu)}$, where $(X^{\mu},Y^{\mu},Z^{\mu})$ is the solution of FBSDE (\ref{Eq3.2}).
\end{proof}

\begin{remark}\label{Remark3.1}
For Theorem \ref{Thm3.1} we have assumed with \textnormal{(H1)} that
\begin{equation}\label{remqrk3.1_eq1}
|f(t, x, y, z, \mu)-f(t, x, y, z, \mu^{\prime})|+|\sigma(t, x, \mu)-\sigma(t, x, \mu^{\prime})| \leq \rho(W_1(\mu, \mu^{\prime})),
\end{equation}
where $(t, x, y, z) \in[0, T] \times \mathbb{R}^d \times \mathbb{R}^d \times \mathbb{R}^d,\, \mu, \mu' \in \mathcal{P}_1(\mathcal{C}_T^2)$, and for some continuity modulus $\rho$ with $\rho: \mathbb{R}_{+}\rightarrow \mathbb{R}_{+},\, \rho(0+)=0$.
However, as we have seen in the proof of Theorem \ref{Thm3.1}, we need to work only with $\mu, \mu' \in \mathcal{K}$, and so we can replace equivalently \eqref{remqrk3.1_eq1} by
\begin{equation}\label{remqrk3.1_eq2}
|f(t, x, y, z, \mu)-f(t, x, y, z, \mu')|+|\sigma(t, x, \mu)-\sigma(t, x, \mu')| \leq \rho(W_2(\mu, \mu^{\prime})),
\end{equation}
for all $(t, x, y, z) \in[0, T] \times \mathbb{R}^d \times \mathbb{R}^d \times \mathbb{R}^d, \mu, \mu' \in \mathcal{K}$, and for some continuity modulus $\rho$ with $\rho(0+)=0$. \\
Indeed, without loss of generality we can suppose that $\rho: \mathbb{R}_{+} \rightarrow \mathbb{R}_{+}$ is increasing. Then, for $C$ introduced in the definition of $\mathcal{K}$,
\begin{equation}\label{remqrk3.1_eq3}
W_1(\mu, \mu^{\prime}) \leq W_2(\mu, \mu') \leq (2C^{\frac{1}{4}})^{\frac{3}{4}}W_1(\mu, \mu')^{\frac{1}{4}},\ \mu,\ \mu' \in \mathcal{K},
\end{equation}
and so
\begin{equation}\label{remqrk3.1_eq4}
\rho(W_1(\mu, \mu^{\prime})) \leq \rho(W_2(\mu, \mu')) \leq \rho'(W_1(\mu, \mu')),\ \mu,\ \mu' \in \mathcal{K},
\end{equation}
for $\rho'(r):=\rho((2 C^{\frac{1}{4}})^{\frac{3}{4}} r^{\frac{1}{4}}), r \geqslant 0$.\\
To see the second inequality in \eqref{remqrk3.1_eq3}, choose two $\mathcal{C}_{T}^2$-valued random variables $\xi$ and $\eta$ over $(\Omega, \mathcal{F}, P)$ such that $P_{\xi}=\mu, P_{\eta}=\mu'$ and $W_1(\mu, \mu')=E[|\xi-\eta|_{\mathcal{C}_{T}^{2}}]$, and note that by the H\"{o}lder inequality
\begin{equation*}
\begin{aligned}
W_2(\mu, \mu^{\prime}) & \leq (E[|\xi-\eta|^{2}_{\mathcal{C}_{T}^{2}}])^{\frac{1}{2}} \leq(E[|\xi-\eta|_{\mathcal{C}_{T}^{2}}])^{\frac{1}{4}}(E[|\xi-\eta|^{4}_{\mathcal{C}_{T}^{2}}])^{\frac{3}{16}} \\
& \leq((E[|\xi|^{4}_{\mathcal{C}_{T}^{2}}])^{\frac{1}{4}}+(E[|\eta|_{\mathcal{C}_{T}^{2}}^4])^{\frac{1}{4}})^{\frac{3}{4}}
(E[|\xi-\eta|_{\mathcal{C}_{T}^{2}}])^{\frac{1}{4}}
 \leq(2 C^{\frac{1}{4}})^{\frac{3}{4}} (W_1(\mu, \mu^{\prime}))^{\frac{1}{4}}.
\end{aligned}
\end{equation*}
\end{remark}

In what follows we consider the following mean field FBSDE,
\begin{equation}\label{eq_FBuni}
\left\{
\begin{aligned}
dX_t &= \sigma(t, X_t, P_{X_{.\wedge t}})dB_t,\ t\in[0,T], \\
dY_t &= -f(t, X_t, Y_t, Z_t, P_{(X, Y_{.\vee t})})dt+Z_tdB_t,\ t\in[0,T],\\
X_{0}&= x\in\mathbb{R},\ Y_T=\Phi(X_T, P_{X}).
\end{aligned}
\right.
\end{equation}
It is a special case of the equation \eqref{EqFBSDE3.1}, and by the above Theorem \ref{Thm3.1} and Remark \ref{Remark3.1}, the existence of the solution is proved. Let us discuss the uniqueness of the solution.
We recall $\mathcal{P}_{2}(\mathcal{C}_T^{2})$ is the set of probability measure $\mu$ on $\mathcal{C}_T^{2}$ with finite second order moment equipped with the Wasserstein distance $W_{2}$,
\begin{equation*}%\label{eq_w2}
\begin{aligned}
W_2(\mu,\mu'):= \inf\Big\{\Bigr\{\int_{\mathcal{C}_T^{2}\times\mathcal{C}_T^{2}}|x-x'|^{2}\rho(dx,dx')\Bigr\}^{\frac{1}{2}}:\ \rho\in\mathcal{P}_{2}(\mathcal{C}_T^{2}\times\mathcal{C}_T^{2}),\ \rho(\cdot\times\mathcal{C}_T^{2})
  =\mu,\ \rho(\mathcal{C}_T^{2}\times\cdot)=\mu' \Big\}.
\end{aligned}
\end{equation*}
Also note that $(\mathcal{C}_T^{2},|\cdot|)$ and $(\mathcal{P}_{2}(\mathcal{C}_T^{2}),W_{2})$ are Polish spaces, they are namely complete and separable.
Let us replace \textnormal{(H1)} by the following assumption:
\begin{itemize}
\item[\bf{(H3)}]\ The coefficients $\sigma, f$ and $\Phi$ are bounded and $\sigma, f$ are Lipschitz with respect to $(x,\nu)$ and to $(x,y,z,\mu)$, respectively. i.e., there exists a constant $C\in\mathbb{R}_{+}$ such that, $P$-a.s., for all $t\in[0,T], (x, y,z),(x',y',z')\in\mathbb{R}^{d}\times\mathbb{R}^{d}\times\mathbb{R}^{d},\mu, \mu'\in\mathcal{P}_{1}(\mathcal{C}_T^{2}), \nu, \nu'\in\mathcal{P}_{1}(\mathcal{C}_T)$,
\begin{equation*}
\begin{aligned}
 |f(t,x,y,z,\mu)-f(t,x',y',z',\mu')|&\leq C(|x-x'|+|y-y'|+|z-z'|+W_{2}(\mu,\mu')),\\
 |\sigma(t,x,\nu)-\sigma(t,x',\nu')|&\leq C(|x-x'|+W_{2}(\nu,\nu')).
\end{aligned}
\end{equation*}
\end{itemize}

\begin{theorem}\label{Theorem_3.2}
Under assumptions \textnormal{(H2)} and \textnormal{(H3)}, equation \eqref{eq_FBuni} has a unique solution $(X,Y,Z)\in S^2_{\mathbb{F}}\times S^2_{\mathbb{F}}\times M^2_{\mathbb{F}}$.
\end{theorem}
\begin{proof}
The existence of the solution follows from Theorem \ref{Thm3.1} and Remark \ref{Remark3.1}.
Let $(X,Y,Z)$ and $(X',Y',Z')$ be two solutions of equation \eqref{eq_FBuni}. Set $\overline{X}=X-X',\ \overline{Y}=Y-Y',\ \overline{Z}=Z-Z'$, then we have
\begin{equation}
\left\{
\begin{aligned}
d\overline{X}_t& =(\sigma(t, X_t, P_{X_{.\wedge t}})-\sigma(t, X'_t, P_{X'_{.\wedge t}}))dB_t, \\
d\overline{Y}_t &=-(f(t, X_t, Y_t, Z_t, P_{(X, Y_{.\vee t})})-f(t, X'_t, Y'_t, Z'_t, P_{(X', Y'_{.\vee t})}))dt+\overline{Z}_tdB_t,\\
\overline{X}_{0}&=0,\ \overline{Y}_T=\Phi(X_T, P_{X_{.\wedge T}})-\Phi(X'_T, P_{X'_{.\wedge T}}).
\end{aligned}
\right.
\end{equation}
Thus, by standard estimates
\begin{equation*}
\begin{aligned}
&\ E[\sup_{t\in[0,s]}|\overline{X}_t |^{2}]
%=&E[\sup_{t\in[0,T]}|\int_{0}^{t}(\sigma(s, X_s, P_{X_{.\wedge s}})-\sigma(s, X'_s, P_{X'_{.\wedge s}}))dB_{s}|^{2}]\\
\leq C'E\[\int_{0}^{s}|\sigma(t, X_t, P_{X_{.\wedge t}})-\sigma(t, X'_t, P_{X'_{.\wedge t}})|^{2}dt\]\\
\leq&\  C'E\[\int_{0}^{s}(|\overline{X}_{t}|+W_{2}(P_{X_{.\wedge t}},P_{X'_{.\wedge t}}))^{2}dt\]
\leq C'E\[\int_{0}^{s}\sup_{r\in[0,t]}|\overline{X}_{r}|^{2}dt\],\ s\in[0,T],
\end{aligned}
\end{equation*}
and from  Gronwall's inequality we obtain  $\displaystyle E[\sup_{t\in[0,T]}|\overline{X}_t |^{2}]= 0$. Consequently, $X=X'$. The uniqueness of the solution of the SDE is proved.
Next we consider the BSDE
\begin{equation}\label{eq_w2}
\left\{\begin{aligned}
		d\overline{Y}_{\!t}& =-(f(t, X_t, Y_t, Z_t, P_{(X, Y_{.\vee t})})-f(t, X'_t, Y'_t, Z'_t, P_{(X', Y'_{.\vee t})}))dt+\overline{Z}_tdB_t,\\
\overline{Y}_{T\!}&=\Phi(X_T, P_{X_{.\wedge T}})-\Phi(X'_T, P_{X'_{.\wedge T}}).
	\end{aligned}\right.
\end{equation}
For some suitable $\beta>0$ we apply It\^{o}'s formula to $e^{\beta t}|\overline{Y}_{t}|^{2}$. This yields
\begin{equation}\label{eq_etY^{2}Ito}
\begin{aligned}
&\ e^{\beta t}|\overline{Y}_{t}|^{2}+\int_{t}^{T}e^{\beta s}|\overline{Z}_{s}|^{2}ds+\beta\int_{t}^{T}e^{\beta s}|\overline{Y}_{s}|^{2}ds+2\int_{t}^{T}e^{\beta s}\overline{Y}_{s}\overline{Z}_{s}dB_{s}\\
=&\ e^{\beta T}|\overline{Y}_{T}|^{2}+2\int_{t}^{T}e^{\beta s}\overline{Y}_{s}\(f(s, X_s, Y_s, Z_s, P_{(X, Y_{.\vee s})})-f(t, X'_s, Y'_s, Z'_s, P_{(X', Y'_{.\vee s})})\)ds,\ t\in [0,T].
\end{aligned}
\end{equation}
From standard estimates we have $\displaystyle E[(\int_{0}^{T}|\overline{Y}_{s}\overline{Z}_{s}|^{2}ds)^{\frac{1}{2}}]<+\infty$. Consequently, $\displaystyle\int_{0}^{t}e^{\beta s}\overline{Y}_{s}\overline{Z}_{s}dB_{s}$ is a uniformly integrable martingale. By taking expectation and using the polarization formula $2ab\leq\alpha a^{2}+\frac{1}{\alpha}b^{2}$, we obtain
\begin{equation*}
\begin{aligned}
&\ E\Bigr[e^{\beta t}|\overline{Y}_{\!t}|^{2}+\int_{t}^{T}e^{\beta s}|\overline{Z}_{s}|^{2}ds+\beta\int_{t}^{T}e^{\beta s}|\overline{Y}_{\!s}|^{2}ds\Bigr]\\
\leq &\ (4C^{2}+\frac{C}{\beta})E\Bigr[e^{\beta T}\sup_{t\in[0,T]}|\overline{X}_{t}|^{2}\Bigr]+CE\Bigr[\int_{t}^{T}e^{\beta s}|\overline{X}_{s}|^{2}ds\Bigr]
 +(4C+2C^{2})E\Bigr[\int_{t}^{T}e^{\beta s}|\overline{Y}_{\!s}|^{2}ds\Bigr]\\
 &\ +\frac{1}{2}E\Bigr[\int_{t}^{T}e^{\beta s}|\overline{Z}_{s}|^{2}ds\Bigr]+Ce^{\beta T}\int_{t}^{T}E[\sup_{r\in[s,T]}e^{\beta r}|\overline{Y}_{\!r}|^{2}]ds,\ t\in[0,T] .\\
\end{aligned}
\end{equation*}
As $\overline{X}_{t}=0,\ t\in[0,T]$, letting $\beta=(4C+2C^{2}+\frac{1}{2})$, we get
\begin{equation}\label{eqnew0}
\begin{aligned}
 E\Bigr[\int_{t}^{T}e^{\beta s}|\overline{Z}_{s}|^{2}ds+\int_{t}^{T}e^{\beta s}|\overline{Y}_{\!s}|^{2}ds\Bigr]
\leq  2Ce^{\beta T }\int_{t}^{T}E[\sup_{r\in[s,T]}e^{\beta r}|\overline{Y}_{\!r}|^{2}]ds,\ \ t\in[0,T] .
\end{aligned}
\end{equation}
On the other hand, again from \eqref{eq_etY^{2}Ito}, Burkholder-Davis-Gundy's inequality and the polarization formula, we deduce that for $t\in[0,T]$,
\begin{equation}\label{eqnew1}
\begin{aligned}
E[\sup_{s\in[t,T]}& e^{\beta s}|\overline{Y}_{s}|^{2}]
%\leq&E\[ e^{\beta T}|\overline{Y}_{T}|^{2}+2\int_{t}^{T}e^{\beta s}|\overline{Y}_{s}(f(s, X_s, Y_s, Z_s, P_{(X, Y_{.\vee s})})-f(t, X'_s, Y'_s, Z'_s, P_{(X', Y'_{.\vee s})}))|ds\]\\
%&+2\Bigr|\sup_{s\in[t,T]}\int_{t}^{T}e^{\beta s}\overline{Y}_{s}\overline{Z}_{s}dB_{s}\Bigr|,\\
\leq (4C^{2}+\frac{C}{\beta})E\Bigr[e^{\beta T}\sup_{t\in[0,T]}|\overline{X}_{t}|^{2}\Bigr]+CE\Bigr[\int_{t}^{T}e^{\beta s}|\overline{X}_{s}|^{2}ds\Bigr]
 +\frac{1}{2}E\Bigr[\int_{t}^{T}e^{\beta s}|\overline{Z}_{s}|^{2}ds\Bigr]\\
 &+\!(4C\!+\!2C^{2})E\Bigr[\!\!\int_{t}^{T}\!\!e^{\beta s}|\overline{Y}_{s}|^{2}ds\Bigr]\!+\!Ce^{\beta T}\!\int_{t}^{T}\!\!\!E[\sup_{r\in[s,T]}e^{\beta r}|\overline{Y}_{\!r}|^{2}]ds
\!+\!6E\[\(\!\int_{t}^{T}e^{2\beta s}|\overline{Y}_{\!s}\overline{Z}_{s}|^{2}ds\)^{\frac{1}{2}}\!\].
\end{aligned}
\end{equation}
Notice that
\begin{equation}\label{eqnew2}
\begin{aligned}
 6E\big[(\int_{t}^{T}e^{2\beta s}|\overline{Y}_{\!s}\overline{Z}_{s}|^{2}ds)^{\frac{1}{2}}\big]
\leq&\ 6 E\big[(\sup_{s\in[t,T]}e^{\frac{\beta s}{2}}|\overline{Y}_{\!s}|)\cdot
(\int_{t}^{T}e^{\beta s}|\overline{Z}_{s}|^{2}ds)^{\frac{1}{2}}\big]\\
\leq&\ \frac{1}{2}E[\sup_{s\in[t,T]}e^{\beta s}|\overline{Y}_{\!s}|^{2}]+18 E\big[\int_{t}^{T}e^{\beta s}|\overline{Z}_{s}|^{2}ds\big].
 \end{aligned}
\end{equation}
As $\overline{X}_{t}=0,\ t\in[0,T]$, by combining with \eqref{eqnew0}-\eqref{eqnew2},  we have, for some constant $C'$,
\begin{equation*}
\begin{aligned}
&E[\sup_{s\in[t,T]} e^{\beta s}|\overline{Y}_{\!s}|^{2}]\leq C'\int_{t}^{T}E[\sup_{r\in[s,T]}e^{\beta r}|\overline{Y}_{\!r}|^{2}]ds,\ t\in[0,T] .
\end{aligned}
\end{equation*}
By applying Gronwall's inequality, we obtain $\displaystyle E[\sup_{s\in[t,T]} e^{\beta s}|\overline{Y}_{s}|^{2}]=0$, for all $t\in[0,T]$, i.e.,
$Y=Y'$.  Then from \eqref{eqnew0}, $Z=Z'$. We conclude that equation  \eqref{eq_FBuni}  has a unique solution.
\end{proof}

\begin{example}
In order to illustrate our results, let us give the following two examples:

\textnormal{i)} $\sigma(t, x, P_{(X, Y)}):=\widetilde{\sigma}(t, x, \int_0^TE[h(X_s, Y_s)]ds)$,
where $h$ is a deterministic function;

\textnormal{ii)} $\sigma(t, x, P_{(X, Y)}):=\widetilde{\sigma}(t, x, E[g(X_{\phi(t)}, Y_{\psi(t)})]),\ \phi,\ \psi: [0, T]\rightarrow[0, T] \rm\ \mbox{are measureable functions}$. \\
In addition to \textnormal{i),\ ii)}, let us suppose that
\begin{itemize}
\item[\textnormal{(A1)}] \
 $\widetilde{\sigma}(\cdot, \cdot, \cdot)$ is measurable, bounded, Lipschitz in $x$ and has continuity modulus $\rho_{1} : \mathbb{R}_{+}\rightarrow\mathbb{R}_{+}$ continuous and increasing, such that $$|\widetilde{\sigma}(t,x,\nu)-\widetilde{\sigma}(t,x,\nu')|^{2}\leq\rho_{1}(|\nu-\nu'|^{2}),\ (t,x)\in[0,T]\times\mathbb{R},\nu,\nu'\in\mathbb{R};\ \text{and}$$
\item[$\textnormal{(A2)}$] $h, g: \mathbb{R}^{2}\rightarrow\mathbb{R}$ are Lipschitz with respect to $(x, y)$.
\end{itemize}
Then, $\sigma$ defined in \textnormal{i)} and in \textnormal{ii)} satisfies \textnormal{(H1)} and \textnormal{(H2)}. Indeed, for example, for \textnormal{ii)}, we have\\
\begin{equation*}
|\sigma(t,x,P_{(X^{1},Y^{1})})-\sigma(t,x,P_{(X^{2},Y^{2})})|^{2}
\leq\rho_{1}(C^{2}(E[(|\widetilde{X}^{1}_{\phi(t)}-\widetilde{X}^{2}_{\phi(t)}|^{2}+|\widetilde{Y}^{1}_{\psi(t)}-\widetilde{Y}^{2}_{\psi(t)}|^{2})^{\frac{1}{2}}])),
\end{equation*}
 for all $(\widetilde{X}^{i},\widetilde{Y}^{i}), i=1,2$, with $P_{(\widetilde{X}^{i}, \widetilde{Y}^{i})}=P_{(X^{i}, Y^{i})}$. Thus,
\begin{equation*}
\begin{aligned}
|\sigma(t,x,P_{(X^{1},Y^{1})})-\sigma(t,x,P_{(X^{2},Y^{2})})|^{2}\leq \rho_{1}(C^{2}E[\sup_{t\in[0,T]}(|\widetilde{X}^{1}_{t}-\widetilde{X}^{2}_{t}|^{2}+|\widetilde{Y}^{1}_{t}-\widetilde{Y}^{2}_{t}|^{2})^{\frac{1}{2}}]),
\end{aligned}
\end{equation*}
from where we deduce that
\begin{equation*}
\begin{aligned}
|\sigma(t,x,P_{(X^{1},Y^{1})})-\sigma(t,x,P_{(X^{2},Y^{2})})|^{2}
\leq \rho_{1}(C^{2}W_{1}^{2}(P_{(X^{1},Y^{1})},P_{(X^{2},Y^{2})}))=\rho(W_{1}^{2}(P_{(X^{1},Y^{1})},P_{(X^{2},Y^{2})})),
\end{aligned}
\end{equation*}
for $\rho(x):=\rho_{1}(C^{2}x),\ x\in\mathbb{R}_{+}$.
\end{example}

\section{Derivative with respect to a measure over a Banach space}\label{FBSDESec4}

Let $\mathcal{K}$ be a real separable Banach space endowed with the norm $|\cdot|_{\mathcal{K}}$. We denote by
$\mathcal{K}'=\{l:\mathcal{K} \rightarrow\mathbb{R}\ |\ l\ \text{continuous linear functional}\}$ the dual space of $\mathcal{K}$ and we endow $\mathcal{K}'\times\mathcal{K}$ with the duality product $\langle l, x\rangle_{\mathcal{K}'\times\mathcal{K}}=l(x),\ l\in\mathcal{K}',\ x\in\mathcal{K}$. We are here mainly motivated by the case, where
\begin{equation}\label{Eq_4.1}
\begin{aligned}
\mathcal{K}=\mathcal{C}([0,T];\mathbb{R}^{n})\times L^{2}([0,T])\ \text{and}\ \mathcal{K}'=BV([0,T];\mathbb{R}^{n})\times L^{2}([0,T]).
\end{aligned}
\end{equation}
Here $BV([0,T];\mathbb{R}^{n})$ is the space of all $\mathbb{R}^{n}$-valued bounded variation c\`{a}dl\`{a}g functions defined over $[0,T]$. Recall that $\mathcal{C}([0,T];\mathbb{R}^{n})'=BV([0,T];\mathbb{R}^{n})$, and
$\langle h, \varphi\rangle_{BV^{n}_{T}\times\mathcal{C}^{n}_{T}}= \int_{[0,T]}\varphi(t) h(dt), \, \varphi\in\mathcal{C}^{n}_{T}:=\mathcal{C}([0,T];\mathbb{R}^{n}),\ h\in BV^{n}_{T}:=BV([0,T];\mathbb{R}^{n})$.

Then, obviously, for \eqref{Eq_4.1},
\begin{equation}\label{Eq_4.2}
\begin{aligned}
\langle (h,y), (\varphi,x)\rangle_{\mathcal{K}'\times\mathcal{K}}=&\int_{[0,T]}\varphi(t) h(dt)+\int_0^Ty(t)x(t)dt,\\
&(\varphi,x)\in\mathcal{C}^{n}_{T}\times L^{2}([0,T];\mathbb{R}^{d}),\ (h,y)\in BV^{n}_{T}\times L^{2}([0,T];\mathbb{R}^{d}).
\end{aligned}
\end{equation}
Let us come back to a general real separable Banach space $(\mathcal{K}, |\cdot|_{\mathcal{K}})$, and recall that $$\mathcal{P}_{2}(\mathcal{K})=\{m\in\mathcal{P}(\mathcal{K})\ \text{probability measure over}\ (\mathcal{K},\ \mathcal{P}(\mathcal{K})): \int_{\mathcal{K}}|x|^{2}_{\mathcal{K}} m(dx)<\infty\}.$$
We now recall the notion of differentiability of a function $f: \mathcal{P}_{2}(\mathcal{K})\rightarrow \mathbb{R}$.
For more details for the case where $\mathcal{K}=\mathbb{R}^{d}$, the reader can refer to Carmona and  Delarue \cite{CD2018I}.\\
\begin{definition}
We say that $u: \mathcal{P}_{2}(\mathcal{K})\rightarrow\mathbb{R}$ has the linear functional derivative $\mathcal{D}_{m}u: \mathcal{P}_{2}(\mathcal{K})\times\mathcal{K}\rightarrow\mathbb{R}$, if $\mathcal{D}_{m}u$ is a continuous function over $\mathcal{P}_{2}(\mathcal{K})\times\mathcal{K}$ with at most quadratic growth such that, for all $m,\, m'\in\mathcal{P}_{2}(\mathcal{K})$,
\begin{equation}\label{Eq_Derivative1}
\begin{aligned}
u(m')-u(m)=\int_{0}^{1}\int_{\mathcal{K}}\mathcal{D}_{m}u(tm'+(1-t)m,x)(m'(dx)-m(dx))dt.
\end{aligned}
\end{equation}
\end{definition}
Let us suppose that for $u:\mathcal{P}_{2}(\mathcal{K})\rightarrow\mathbb{R}$ the derivative $\mathcal{D}_{m}u:\mathcal{P}_{2}(\mathcal{K})\times\mathcal{K}\rightarrow\mathbb{R}$ exists, is continuous and of at most quadratic growth, and that, for all $m\in\mathcal{P}_{2}(\mathcal{K}), \mathcal{D}_{m}u(m,\cdot): \mathcal{K}\rightarrow\mathbb{R}$ is differentiable, i.e., there exists $\partial_{x}(\mathcal{D}_{m}u)(m,\cdot):\mathcal{K}\rightarrow\mathcal{K}'$ such that, for all $x\in\mathcal{K}$,
\begin{equation}\label{Eq_Derivative2}
\begin{aligned}
\mathcal{D}_{m}u(m,y)-\mathcal{D}_{m}u(m,x)
=\langle\partial_{x}(\mathcal{D}_{m}u)(m,x),y-x\rangle_{\mathcal{K}'\times\mathcal{K}}+o(|y-x|_{\mathcal{K}}),\  \text{as}\ \mathcal{K}\ni y\rightarrow x.
\end{aligned}
\end{equation}

Let us make the following assumptions for $\mathcal{D}_{m}u:\mathcal{P}_{2}(\mathcal{K})\times\mathcal{K}\rightarrow\mathbb{R}$:
\begin{itemize}
\item[\bf{(H4)}] The derivative $\partial_{x}(\mathcal{D}_{m}u):\mathcal{P}_{2}(\mathcal{K})\times\mathcal{K}\rightarrow\mathcal{K}'$ exists, is continuous and of at most linear growth.
\end{itemize}
However, to simplify our arguments, we suppose also {\bf{(H4)'}}:
\begin{itemize}
\item[i)] $\partial_{x}(\mathcal{D}_{m}u):\mathcal{P}_{2}(\mathcal{K})\times\mathcal{K}\rightarrow\mathcal{K}'$ is bounded, and
\item[ii)] There exists a continuity modulus $\rho_{u}:\mathbb{R}_{+}\rightarrow \mathbb{R}_{+}$ continuous and increasing, with $\rho_{u}(0)=0$ and $\rho^{2}_{u}(\cdot)$ is concave such that
\begin{equation*}
\begin{aligned}
|\partial_{x}(\mathcal{D}_{m}u)(m',x')-\partial_{x}(\mathcal{D}_{m}u)(m,x)|_{\mathcal{K}'}\leq \rho_{u}(W_{2}(m,m')+|x-x'|_{\mathcal{K}}),\ (m',x'),\ (m,x)\in \mathcal{P}_{2}(\mathcal{K})\times\mathcal{K}.
\end{aligned}
\end{equation*}
\end{itemize}
Given any $m,m'\in \mathcal{P}_{2}(\mathcal{K})$ we put $$\Pi(m,m')=\{ \Pi\in\mathcal{P}_{2}(\mathcal{K}\times\mathcal{K}): \Pi(\cdot\times\mathcal{K})=m,\ \Pi(\mathcal{K}\times\cdot)=m'\},$$
and we denote by $\Pi_{0}(m,m')$ the set of $\Pi\in\Pi(m,m')$ for which
\begin{equation*}
\begin{aligned}
\int_{\mathcal{K}\times\mathcal{K}}|x-y|^{2}\Pi(dxdy)=W_{2}^{2}(m,m').
\end{aligned}
\end{equation*}

It is well-known that $\Pi_{0}(m,m')\neq\emptyset$.
\begin{proposition}\label{Prop_3.1_taylor}
Under our assumption (H4)' on $u$ and $\mathcal{D}_{m}u$ we have the following first order Taylor expansion for $u:\mathcal{P}_{2}(\mathcal{K})\rightarrow \mathbb{R}$ at $m\in\mathcal{P}_{2}(\mathcal{K})$ holds true:
\begin{equation}\label{Taylor_expansion1}
\begin{aligned}
u(m')=u(m)+\int_{\mathcal{K}} \mathcal{D}_{m}u(m,x)(m'(dx)-m(dx))+o(W_{2}(m,m')),\ \text{as}\ W_{2}(m,m')\rightarrow 0\ (m'\in\mathcal{P}_{2}(\mathcal{K})).
\end{aligned}
\end{equation}
\end{proposition}

The proof of this proposition is a slight extension of the corresponding result in Carmona and Delarue \cite{CD2018I}. Moreover, also under our above assumptions the following result concerning Lion's L-derivative of $u:\mathcal{P}_{2}(\mathcal{K})\rightarrow\mathbb{R}$ extends easily from $\mathcal{K}=\mathbb{R}^{n}$ (see Carmona and Delarue \cite{CD2018I})
 to the general separable Banach space $\mathcal{K}$.
\begin{proposition}\label{Prop3.2L-derivatives}
Given any $\xi\in L^{2}(\Omega,\mathcal{F},P;\mathcal{K})$, we have
\begin{equation}\label{Eq_L_derivative}
\begin{aligned}
u(P_{\xi+\eta})=u(P_{\xi})+E[\langle\partial_{x}(\mathcal{D}_{m}u)(P_{\xi},\xi),\eta\rangle_{\mathcal{K}'\times\mathcal{K}}]
+R(\xi,\eta),\ \eta\in L^{2}(\Omega,\mathcal{F},P;\mathcal{K}),
\end{aligned}
\end{equation}
where $|R(\xi,\eta)|=o((E[|\eta|^{2}_{\mathcal{K}}])^{\frac{1}{2}})$, as $E[|\eta|^{2}_{\mathcal{K}}])^{\frac{1}{2}}\rightarrow0$.
\end{proposition}
\begin{definition}
We denote $L$-derivative by $(\partial_{\mu}u(P_{\xi},x))$:
\begin{equation}
\begin{aligned}
\partial_{\mu}u(m,x)=\partial_{x}(\mathcal{D}_{m}u)(m,x),\ (m,x)\in\mathcal{P}_{2}(\mathcal{K})\times\mathcal{K}.
\end{aligned}
\end{equation}
Moveover, when we speak about the differentiability of $u:\ \mathcal{P}_{2}(\mathcal{K})\rightarrow \mathbb{R}$,
we mean the L-differentiability.
\end{definition}

As we have already mentioned, for our studies we are mainly interested in the case $\mathcal{K}=\mathcal{C}_{T}^{k}\times L^{2}([0,T];\mathbb{R}^{k})$, where $\mathcal{C}_{T}^{k}:=\mathcal{C}([0,T];\mathbb{R}^{k})$.

And so, if, for instance, $\sigma:\mathcal{P}_{2}(\mathcal{C}_{T}^{k}\times {L}^{2}([0,T];\mathbb{R}^{k}))\rightarrow\mathbb{R}$ is differentiable we have $\partial_{\mu}\sigma:\mathcal{P}_{2}(\mathcal{C}_{T}^{k}\times L^{2}([0,T];\mathbb{R}^{k}))\times\mathcal{C}_{T}^{k}\times{L}^{2}([0,T];\mathbb{R}^{k})\rightarrow BV_{T}^{k}\times{L}^{2}([0,T];\mathbb{R}^{k})$, and as the elements of $BV_{T}^{k}\times{L}^{2}([0,T];\mathbb{R}^{k})$ take their values in $\mathbb{R}^{k+1}$, we denote the $i$-th component of $(\partial_{\mu}\sigma)(m,\varphi)$ by $(\partial_{\mu}\sigma)_{i}(m,\varphi)$,
$ (\partial_{\mu}\sigma)(m,\varphi)=((\partial_{\mu}\sigma)_{i}(m,\varphi))_{1\leq i\leq k+1},\ (m,\varphi)\in \mathcal{P}_{2}(\mathcal{C}_{T}^{k}\times{L}^{2}([0,T];\mathbb{R}^{k}))\times\mathcal{C}_{T}^{k}\times{L}^{2}([0,T];\mathbb{R}^{k}).$

Moreover, we will write $\langle\cdot,\cdot\rangle_{k}=\langle\cdot,\cdot\rangle_{\mathcal{K}'\times\mathcal{K}}$ for the duality product.
\begin{definition}
For $l=(l_{1},\cdot\cdot\cdot,l_{k})\in BV_{T}^{k},\ \widehat{v}\in{L}^{2}([0,T];\mathbb{R}^{k}), (f(=(f_{1},\cdot\cdot\cdot,f_{k})),v)\in \mathcal{C}^{k}_{T}\times{L}^{2}([0,T];\mathbb{R}^{k}),$ the inner product $\langle\cdot,\cdot\rangle_{k}$ is defined as \\
$$\langle(l,\widehat{v}),(f,v)\rangle_{k}:=\sum\limits_{i=1}^{k}\int_{0}^{T}f_{i}(t)l_{i}(dt)+\int_{0}^{T}\widehat{v}_{s}v_{s}ds.$$
\end{definition}

\section{Maximum principle for the controlled mean-field FBSDEs}\label{FBSDESec5}
Let $\mathcal{U}_{T}:=L^{2}([0,T]; \mathbb{R}^{k})$ and  $U$ be a nonempty convex subset of $\mathbb{R}^k$. We denote by $L_{\mathbb{F}}^{\infty-}(\Omega,L^{2}([0,T];\mathbb{R}^{k}))$ the space of all $\mathbb{F}$-adapted process with values in $\mathbb{R}^k$ such that $\displaystyle E[(\int_{0}^{T}|v(t)|^{2}dt)^{\frac{p}{2}}]<+\infty,$ for all $p\geq2$.  An element $v$ of the space
$$\mathcal{U}_{ad}=\{v(\cdot)\in L_{\mathbb{F}}^{\infty-}(\Omega,L^{2}([0,T];\mathbb{R}^{k}))|\ v_{t}\in U,\ 0\leq t \leq T,\ a.s.a.e.\}$$
 is called an admissible control.

We consider the following controlled mean-field coupled forward-backward SDE:
\begin{equation}\label{Eqcontrol}
\left\{\begin{aligned}
dX^{v}_t&=\sigma(t,X_t,P_{(X^{v}_{\cdot\wedge t},v)})dB_t,\ t\in[0,T],\\
dY^{v}_t&=-f(t,X^{v}_t,Y^{v}_t,Z^{v}_t,P_{(X^{v},Y^{v}_{t\vee\cdot},v)})dt+Z^{v}_tdB_t,\ t\in[0,T],\\
X^{v}_{0}&=x\in\mathbb{R},\ Y^{v}_T=\Phi(X^{v}_T,P_{(X^{v},v)}),
\end{aligned}\right.
\end{equation}
where $v\in\mathcal{U}_{ad}$ and
\begin{equation*}
  \begin{aligned}
  &f:\ [0,T]\times\mathbb{R}\times\mathbb{R}\times\mathbb{R}\times\mathcal{P}_{2}(\mathcal{C}_T^{2}\times\mathcal{U}_{T})\rightarrow\mathbb{R},\\
  &\sigma:\ [0,T]\times\mathbb{R}\times\mathcal{P}_{2}(\mathcal{C}_T\times\mathcal{U}_{T})\rightarrow\mathbb{R},\\
  &\Phi:\ \mathbb{R}\times\mathcal{P}_{2}(\mathcal{C}_T\times\mathcal{U}_{T})\rightarrow\mathbb{R},\\
  \end{aligned}
\end{equation*}
and we define the following cost functional:

\begin{equation*}\label{Eqcost}
\begin{aligned}
 J(v)=E[ & \displaystyle \int_0^T L(t,X^{v}_{t},Y^{v}_{t},Z^{v}_{t},P_{(X^{v},Y^{v},v)})dt+\varphi(X^{v}_{T},P_{(X^{v},Y^{v},v)})],
\end{aligned}
\end{equation*}
where
$$\begin{array}{lll}
& &L:[0,T]\times\mathbb{R}\times\mathbb{R}\times\mathbb{R}\times\mathcal{P}_{2}(\mathcal{C}_T^{2}\times\mathcal{U}_{T})\rightarrow\mathbb{R}, \, \, \varphi:\mathbb{R}\times\mathcal{P}_{2}(\mathcal{C}^{2}_T\times\mathcal{U}_{T})\rightarrow\mathbb{R}.\\
\end{array}
$$
The optimal control problem consists in minimizing the cost functional $J(v)$ over all admissible controls. An admissible control $u\in\mathcal{U}_{ad}$ is called optimal, if the cost functional $J(v)$ attains its minimum at $u$.
Equation (\ref{Eqcontrol}) is called the state equation, the solution associated with $u(\cdot)$ is denoted by $(X,Y,Z):=(X^{u},Y^{u},Z^{u})$.
Recall the definition of partial derivatives with respect to the measure that we have introduced in Section 4,  $\partial_{\mu}\sigma=((\partial_{\mu}\sigma)_{1},(\partial_{\mu}\sigma)_{2}),\ \partial_{\mu}f=((\partial_{\mu}f)_{1}, (\partial_{\mu}f)_{3}, (\partial_{\mu}f)_{3})$.
We have the following standard assumptions:
\begin{itemize}
\item[\bf{(H5)}] The functions $ f,\ \sigma,\ L,\ \Phi,\ \mbox{and}\ \varphi$ are bounded and continuously differentiable to $(x,\ y,\ z,\ \mu)$ and the derivatives are bounded.
\end{itemize}
Recall that $(\partial_{\mu}\sigma)_{1}(t,x,\mu;(\phi,\widehat{v}))\in BV_{T}$. Identifying the elements of $BV_{T}$ with the measures on $([0,T];\mathcal{B}([0,T]))$ they generate, we suppose:
\begin{itemize}
\item[\bf{(H6)}] (i) For fixed $(t,x,\mu;(\phi,\widehat{v}))\in [0,T]\times\mathbb{R}\times\mathcal{P}_{2}(\mathcal{C}_{T}\times \mathcal{U}_{T})\times(\mathcal{C}_{T}\times \mathcal{U}_{T})$, $(\partial_{\mu}\sigma)_{1}(t,x,\mu;(\phi,\widehat{v}))(dr) =(\partial_{\mu}\sigma)_{1}(t,x,\mu;(\phi,\widehat{v}))(r)dr$, where the function at the right-hand side
$(\partial_{\mu}\sigma)_{1}:[0,T]\times\mathbb{R}\times\mathcal{P}_{2}(\mathcal{C}_{T}\times\mathcal{U}_{T})\times (\mathcal{C}_{T}\times\mathcal{U}_{T})\rightarrow \mathbb{R}$ is supposed to be Borel measurable and bounded as well as Lipschitz with respect to $(x,\mu,\phi,\widehat{v})$, uniformly with respect to $t, r\in[0,T]$.
Also for $(\partial_{\mu}\Phi)_{1}$ and for $(\partial_{\mu}f)_{i}, (\partial_{\mu}L)_{i}, (\partial_{\mu}\varphi)_{i}, i=1,2$, we make the same identification as for $(\partial_{\mu}\sigma)_{1}$, and we suppose
\end{itemize}
\begin{equation*}
\begin{aligned}
(\partial_{\mu}\Phi)_{1}(t,x,\mu;(\phi,\widehat{v}))(dr)
=&\ (\partial_{\mu}\Phi)_{1}(t,x,\mu;(\phi,\widehat{v}))(r)dr,\\
(\partial_{\mu}f)_{i}(t,x,y,z,\mu;(\phi_{1},\phi_{2},\widehat{v}))(dr) =&\ (\partial_{\mu}f)_{i}(t,x,y,z,\mu;(\phi_{1},\phi_{2},\widehat{v}))(r)dr, \\
 (\partial_{\mu}L)_{i}(t,x,y,z,\mu;\mu;(\phi_{1},\phi_{2},\widehat{v}))(dr)
 =&\ (\partial_{\mu}L)_{i}(t,x,y,z,\mu;(\phi_{1},\phi_{2},\widehat{v}))(r)dr,\\
(\partial_{\mu}\varphi)_{i}(t,x,y,z,\mu;\mu;(\phi_{1},\phi_{2},\widehat{v}))(dr)
=&\ (\partial_{\mu}\varphi)_{i}(t,x,y,z,\mu;(\phi_{1},\phi_{2},\widehat{v}))(r)dr,
 \end{aligned}
\end{equation*}
where for the functions at the right-hand side we make the following assumptions:
\begin{itemize}
\item[\ ]
(ii) $(\partial_{\mu}f)_{i}, i=1,2$, is Borel measurable, bounded and Lipschitz with respect to $(x, y, z, \mu,  \phi_{1}, \phi_{2}, \widehat{v})$, uniformly with respect to $t, r\in[0,T]$;
 $(\partial_{\mu}\Phi)_{1}$ is Borel measurable, bounded and Lipschitz with respect to $(x, \mu, \phi, \widehat{v})$, uniformly with respect to $t, r\in[0,T]$;\\
(iii) $(\partial_{\mu}L)_{i}, i=1,2,$ is Borel measurable, bounded and Lipschitz with respect to $(x,y,z,\mu, \phi_{1},\phi_{2},\widehat{v})$, uniformly with respect to $t, r\in[0,T]$, in the same sense we suppose for $(\partial_{\mu}\varphi)_{i},i=1, 2,$ is Borel measurable, bounded and Lipschitz with respect to $(t,x,y,z,\mu,\phi_{1},\phi_{2},\widehat{v})$, uniformly with respect to $t, r\in[0,T]$.
\end{itemize}

Under the above assumptions, the existence and the uniqueness of the solution of \eqref{Eqcontrol} can be shown by the arguments used for the proof of the Theorems \ref{Thm3.1} and \ref{Theorem_3.2}.
\begin{proposition}
Suppose {\rm(H5)} and {\rm(H6)}. Then the controlled MFFBSDE \eqref{Eqcontrol} has a unique solution.
\end{proposition}
\begin{proof}
For any arbitrary but fixed admissible control $v\in\mathcal{U}_{T}$, let $$\mathcal{H}_{1}=\{P_{(X,Y,v)}\in\mathcal{P}_{1}(\mathcal{C}^{2}_{T}\times \mathcal{U}_{T}) |\ X, Y \text{continuous adapted process over}\ [0,T] \},$$
replace the space $\mathcal{P}_{1}(\mathcal{C}^{2}_{T})$ in the proof of Theorem \ref{Thm3.1}. The proof of \eqref{Eq3.10Z} is not concerned by the change of the space of measures $\mu$. In the argument for the proof of Theorem \ref{Thm3.1} the normed space $\mathcal{M}_{1}(\mathcal{C}_{T}^{2})$ is replaced by $\mathcal{M}_{1}(\mathcal{C}_{T}^{2}\times \mathcal{U}_{T})$ defined in an obvious way. Let $(\varphi,\psi)=(\varphi_{t},\psi_{t})_{t\in[0,T]}$ be the coordinate process on $\mathcal{C}^{2}_{T}\times \mathcal{U}_{T}$, $(\varphi_{t},\psi_{t})(\phi)=\phi_{t}, \phi\in \mathcal{C}^{2}_{T}\times \mathcal{U}_{T}, t\in[0,T]$.
Then we replace $\mathcal{K}$ in the proof of Theorem \ref{Thm3.1} by
$$\mathcal{K}_{1}=\Bigr\{\mu\in\mathcal{H}_{1}\Bigr| \int_{\mathcal{C}^{2}_{T}\times \mathcal{U}_{T}}|\varphi|^{4}\mu(d\varphi d\psi)\leq C,  \int_{\mathcal{C}^{2}_{T}\times \mathcal{U}_{T}}|\varphi_{t}-\varphi_{s}|^{4}\mu(d\varphi d\psi)\leq C|t-s|^{2},t,s\in[0,T]\Bigr\}.$$
Taking into account the arguments given in the proof of Theorem \ref{Thm3.1}, it is straight-forward to show that
$\mathcal{K}_{1}\subset\mathcal{H}_{1}\subset \mathcal{P}_{1}(\mathcal{C}^{2}_{T}\times L^{2}([0,T]))$ is convex and compact.
Let us define now the map $\mathcal{T}$ by $\mathcal{T}(\mu):=P_{(X^{\mu},Y^{\mu},v)}, \mu\in\mathcal{H}_{1}$.
Here $(X^{\mu},Y^{\mu})$ is defined by \eqref{Eq3.2}, but now with $\mu\in\mathcal{H}_{1}$. Using the argument of the proof of Theorem \ref{Thm3.1} we get that  $\mathcal{T}$ maps $\mathcal{K}_{1}$ into $\mathcal{K}_{1}$ and is continuous on $\mathcal{K}_{1}$. This allows also here the application of Schauder's fixed point theorem. The fixed point
$\mu=\mathcal{T}(\mu)=P_{(X^{\mu},Y^{\mu},v)}$ thus obtained has as consequence that $(X^{\mu},Y^{\mu},Z^{\mu})$ is a solution of MFFBSDE \eqref{Eqcontrol}. Finally, the proof of the uniqueness is a straight-forward extension of the approach used for the proof of Theorem \ref{Theorem_3.2}.
\end{proof}

Let $u$ be an optimal control and $(X, Y, Z):=(X^{u}, Y^{u}, Z^{u})$ be the corresponding optimal solution.
Let $v\in\mathcal{U}_{T}$ be such that $u+v\in\mathcal{U}_{ad}$. Since $U$ is convex, then also for any
  $0\leq\rho\leq1,\ u^{\rho}=u+\rho v$ is in $\mathcal{U}_{ad}$.

We introduce the following linear mean-field FBSDE:
\begin{equation}\label{eq_vari}
\left\{
\begin{aligned}
dX^1_{t}=&\{\sigma_x(t,X_{t}, P_{(X_{\cdot\wedge t},u)})X^1_{t}\\
        &+\widetilde{E}[\langle(\partial_{\mu}\sigma)(t,X_{t},P_{(X_{\cdot\wedge t},u)}; (\widetilde{X}_{\cdot\wedge t}, \widetilde{u})),\ (\widetilde{X}_{\cdot\wedge t}^{1}, \widetilde{v})\rangle_{1}] \}dB_t,\ t\in[0,T],\\
dY^1_{t}=&-\{f_x(\theta_{t})X^1_{t}+f_y(\theta_{t})Y^1_{t}+f_z(\theta_{t})Z^1_{t}\\%+f_v(\theta_{t},u_t)v_{t}
         &+\widetilde{E}[\langle(\partial_{\mu}f)(\theta_{t}; (\widetilde{X}, \widetilde{Y}_{\cdot\vee t},\widetilde{u})),\ (\widetilde{X}^{1}, \widetilde{Y}_{\cdot\vee t}^{1},\widetilde{v})\rangle_{2}]\}dt +Z^1_{t}dB_t,\ t\in[0,T],\\
X^1_{0}=&0,\ Y^1_{T}=\Phi_{x}(X_{T},P_{(X,u)})X^1_{T}+\widetilde{E}[\langle(\partial_{\mu}\Phi)(X_{T},P_{(X,u)}; (\widetilde{X},\widetilde{u})),\ (\widetilde{X}^{1},\widetilde{v})\rangle_{1}].
\end{aligned}
\right.
\end{equation}
where $\theta_{t}:=(t,X_{t},Y_{t},Z_{t},P_{(X,Y_{\cdot\vee t,u})})$.
Recall that, for example, $\langle(\partial_{\mu}f)(\theta(t);(\widetilde{X},\widetilde{Y}_{\cdot\vee t},\widetilde{u})),(\widetilde{X}^{1},\widetilde{Y}^{1}_{\cdot\vee t},\widetilde{v})\rangle_{2}=\int_{0}^{T}(\partial_{\mu}f)(\theta(t);(\widetilde{X},\widetilde{Y}_{\cdot\vee t},\widetilde{u}))(r)\cdot(\widetilde{X}_{r}^{1},\widetilde{Y}^{1}_{r\vee t},\widetilde{v}_{r})^{T}dr$.
%%%%%%%%%%%%%%%%%%
Equation (\ref{eq_vari}) is called the variational equation. First of all, we prove the existence and the uniqueness of the solution $(X^1, Y^1, Z^1)$ of mean-field FBSDE (\ref{eq_vari}).
\begin{lemma}\label{lem5.1}
Let the assumptions \textnormal{(H5)} and \textnormal{(H6)} be satisfied. Then the above linear mean-field FBSDE (\ref{eq_vari}) has a unique solution $(X^{1},Y^{1},Z^{1})\in S_{\mathbb{F}}^{2}\times S_{\mathbb{F}}^{2}\times M_{\mathbb{F}}^{2}$.
\end{lemma}
\begin{proof}
We prove first the existence and the uniqueness of the solution to the forward equation in the variational system (\ref{eq_vari}):
\begin{equation}\label{eq_variSDE}
\left\{
\begin{aligned}
dX^1_{t}=&\ \{\sigma_x(t,X_{t}, P_{(X_{\cdot\wedge t},u)})X^1_{t}
%+\sigma_{v}(t,X_{t},P_{(X_{\cdot\wedge t},u)},u_t)v_{t}
         +\widetilde{E}[\langle(\partial_{\mu}\sigma)(t,X_{t},P_{(X_{\cdot\wedge t})}; (\widetilde{X}_{\cdot\wedge t}, \widetilde{u})),\ (\widetilde{X}_{\cdot\wedge t}^{1}, \widetilde{v})\rangle_{1}] \}dB_t,\\
        X_{0}^{1}=&0.
\end{aligned}
\right.
\end{equation}
Let us begin with the proof of the uniqueness of the solution. Let $X^{1},X^{2}$ be two solutions in ${S}_{\mathbb{F}}^{2}.$ We put $X^{3}:=X^{1}-X^{2}$. Then, $X^{3}_{0}=0$, and from the linearity of equation \eqref{eq_variSDE} in $(\widetilde{X}^{i},\widetilde{v})$ we have
\begin{equation*}
\begin{aligned}
dX_{t}^{3}=\{\sigma_x(t,X_{t}, P_{(X_{\cdot\wedge t},u)})X^3_{t}+\widetilde{E}[\langle(\partial_{\mu}\sigma)(t,X_{t},P_{(X_{\cdot\wedge t},u)}; (\widetilde{X}_{\cdot\wedge t}, \widetilde{u})),\ (\widetilde{X}_{\cdot\wedge t}^{3}, 0)\rangle_{1}]\}dB_{t},\ X_{0}^{3}&=0.
\end{aligned}
\end{equation*}
Taking into account that $|\sigma_{x}|\leq C$\ and $|(\partial_{\mu}\sigma)_{1}|\leq C$, a standard estimate yields
\begin{equation*}
\begin{aligned}
&E[\sup_{0\leq s\leq t}|X_{s}^{3}|^{2}]\\
\leq&\ 4E\[\int_{0}^{t} |\sigma_x(s,X_{s}, P_{(X_{\cdot\wedge s},u)})X^3_{s}+\widetilde{E}[\langle(\partial_{\mu}\sigma)(s,X_{s},P_{(X_{\cdot\wedge s},u)}; (\widetilde{X}_{\cdot\wedge s}, \widetilde{u})),\ (\widetilde{X}_{\cdot\wedge s}^{3}, 0)\rangle_{1}]|^{2}ds\]\\
=&\ 4E\[\int_{0}^{t}|\sigma_x(s,X_{s}, P_{(X_{\cdot\wedge s},u)})X^3_{s}+\widetilde{E}\[\int_{0}^{T} \widetilde{X}_{r \wedge s}^{3}(\partial_{\mu}\sigma)_{1}(s,X_{s},P_{(X_{\cdot\wedge s},u)}; (\widetilde{X}_{\cdot\wedge s},\widetilde{u})) (r)dr\]|^{2}ds\]\\
\leq&\ 16C^{2}\int_{0}^{t}E[\sup_{0\leq r\leq s}|X_{r}^{3}|^{2}]ds,\ t\in[0,T].
\end{aligned}
\end{equation*}

From Gronwall's inequality we obtain $X^{3}=0$, i.e., $X^{1}=X^{2}$.\\
Now we construct the solution by Picard iteration. Let $X^{1,0}_{t}:=\eta\in S^2_{\mathbb{F}},\ 0\leq t\leq T$, and for $n\geq 1$, define
\begin{equation}\label{eq_5.4}
\begin{aligned}
X_{t}^{1,n}=&\int_{0}^{t}
\{\sigma_x(s,X_{s}, P_{(X_{\cdot\wedge s},u)})X^{1,n-1}_{s}%+\sigma_{v}(s,X_{s},P_{(X_{\cdot\wedge s},u)},u_s)v_{s}\\
+\widetilde{E}[\langle(\partial_{\mu}\sigma)(s,X_{s},P_{(X_{\cdot\wedge\cdot s},u)}; (\widetilde{X}_{\cdot\wedge s}, \widetilde{u})),\ (\widetilde{X}_{\cdot\wedge s}^{1,n-1}, \widetilde{v})\rangle_{1}]\}dB_{s},\ t\in[0,T].
\end{aligned}
\end{equation}
First of all, we prove by induction that $X^{1,n}_{t}\in S^2_{\mathbb{F}}$. Indeed, $X^{1,0}=\eta\in S^2_{\mathbb{F}}$, and assuming $X^{1,n-1}_{t}\in S^2_{\mathbb{F}}$, thanks to the assumption (H3), we have that $X^{1,n}$ is adapted and
\begin{equation*}
\begin{aligned}
E[\sup_{t\in[0,T]}|X_{t}^{1,n}|^{2}]
%=&E\[\sup_{0\in[0,T]}\|\int_{0}^{t}
%\{\sigma_x(t,X_{s}, P_{(X_{\cdot\wedge t},u)})X^{1,n-1}_{s}\\%+\sigma_{v}(s,X_{s},P_{(X_{\cdot\wedge s},u)},u_s)v_{s}\\
%&+\widetilde{E}[\langle(\partial_{\mu}\sigma)(s,X_{s},P_{(X_{s\wedge\cdot},u)}; (\widetilde{X}_{s\wedge\cdot}, \widetilde{u})),\ (\widetilde{X}_{s\wedge\cdot}^{1,n-1}, \widetilde{v})\rangle_{1}]\}dB_{s}\|^{2}\]\\
\leq& CE\[\int_{0}^{T}\Bigr\{\sigma_x(t,X_{t}, P_{(X_{\cdot\wedge t},u)})X^{1,n-1}_{t}\\%+\sigma_{v}(t,X_{t},P_{(X_{\cdot\wedge t},u)},u_t)v_{t}\\
&\ +\widetilde{E}\[\langle(\partial_{\mu}\sigma)(t,X_{t},P_{(X_{\cdot\wedge t},u)}; (\widetilde{X}_{\cdot\wedge t}, \widetilde{u})),\ (\widetilde{X}_{\cdot\wedge t}^{1,n-1}, \widetilde{v})\rangle_{1}\]\Bigr\}^{2}dt\],
\end{aligned}
\end{equation*}
and as $\sigma_{x}, \sigma_{v}$ and $\partial_{\mu}\sigma$ are bounded,
\begin{equation*}
\begin{aligned}
E[\sup_{t\in[0,T]}|X_{t}^{1,n}|^{2}]
\leq&\  CE\[\int_{0}^{T}\{\sigma_x(t,X_{t}, P_{(X_{\cdot\wedge t},u)})X^{1,n-1}_{t}\}^{2}dt\]\\%+\sigma_{v}(t,X_{t},P_{(X_{\cdot\wedge t},u)},u_t)v_{t}\}^{2}dt]\\
&\ +CE\[\int_{0}^{T}|\widetilde{E}[\langle(\partial_{\mu}\sigma)(t,X_{t},P_{(X_{\cdot\wedge t},u)}; (\widetilde{X}_{\cdot\wedge t}, \widetilde{u})),\ (\widetilde{X}_{\cdot\wedge t}^{1,n-1}, \widetilde{v})\rangle_{1}]|^{2})]dt\]\\
%\leq& CE[\int_{0}^{T}|X^{1,n-1}_{t}|^{2}dt]+CE[\int_{0}^{T}|v_{t}|^{2}dt]\\
%&+CE[\int_{0}^{T}\widetilde{E}[|\int_{0}^{T}(\partial_{\mu}\sigma)_{1}(t,X_{t},P_{(X_{t\wedge\cdot})},u_t; (\widetilde{X}_{t\wedge\cdot}, \widetilde{u}))(r)\widetilde{X}_{r\wedge t}^{1,n-1}dr|^{2}\\
%&+\int_{0}^{T}|(\partial_{\mu}\sigma)_{2}(t,X_{t},P_{(X_{t\wedge\cdot},u)},u_t; (\widetilde{X}_{t\wedge\cdot}, \widetilde{u}))(r)\widetilde{v}_{r}dr|^{2})]dt]\\
\leq&CE\[\sup_{t\in[0,T]}|X^{1,n-1}_{t}|^{2}\]%+CE\[\int_{0}^{T}|v_{t}|^{2}
%%+CE[\sup_{0\in[0,T]}|X^{1,n-1}_{t}|^{2}]\times \widetilde{E}\[\{\int_{0}^{T}|(\partial_{\mu}\sigma)_{1}(t,X_{t},P_{(X_{t\wedge\cdot},u)}; (\widetilde{X}_{t\wedge\cdot}, \widetilde{u}))(dr)|\}^{2}\]\\
%%&+\widetilde{E}\[\int_{0}^{T}|(\partial_{\mu}\sigma)_{2}(t,X_{t},P_{(X_{t\wedge\cdot})}; (\widetilde{X}_{t\wedge\cdot}, \widetilde{u}))(dr)|^{2}\]\times \widetilde{E}\[\int_{0}^{T}|\widetilde{v}_{r}|^{2}drdt\]
+C\widetilde{E}\[\int_{0}^{T}|\widetilde{v}_{r}|^{2}dr\]
< +\infty.
\end{aligned}
\end{equation*}
Consequently, $X^{1,n}\in S^2_{\mathbb{F}},\ n\geq 1$.
Now we put $\Delta X^{1,n}=X^{1,n}_{t}-X^{1,n-1}_{t},\ n\geq1$, and again from the linearity of \eqref{eq_variSDE} in $(\widetilde{X}^{i},\widetilde{v})$, for $t\in[0,T]$,
\begin{equation}\label{eq_Picard2}
\begin{aligned}
\Delta X_{t}^{1,n}=
\int_{0}^{t}\{\sigma_x(s,X_{s}, P_{(X_{s\wedge\cdot},u)})\Delta X^{1,n-1}_{s}+\widetilde{E}[\langle(\partial_{\mu}\sigma)(s,X_{s},P_{(X_{s\wedge\cdot},u)}; (\widetilde{X}_{s\wedge\cdot}, \widetilde{u})),\ (\Delta\widetilde{X}_{s\wedge\cdot}^{1,n-1}, 0)\rangle_{1}]\}dB_{s}.
\end{aligned}
\end{equation}
Then, by Burkholder-Davis-Gundy inequality, as $\sigma_{x}$ and $(\partial_{\mu}\sigma)_{1}$ are bounded,
\begin{equation*}
\begin{aligned}
&E[\sup_{s\in[0,t]}|\Delta X_{s}^{1,n}|^{2}]\\
\leq&\ C E\[\int_{0}^{t}\{|\sigma_x(s,X_{s}, P_{(X_{s\wedge\cdot},u)})\Delta X^{1,n-1}_{s}|^{2}+|\widetilde{E}[\langle(\partial_{\mu}\sigma)(s,X_{s},P_{(X_{s\wedge\cdot},u)}; (\widetilde{X}_{s\wedge\cdot}, \widetilde{u})),\ (\Delta\widetilde{X}_{s\wedge\cdot}^{1,n-1}, 0)\rangle_{1}]|^{2}\}ds\]\\
\leq&\ C\int_{0}^{t}E[\sup_{r\in[0,s]}|\Delta X^{1,n-1}_{r}|^{2}]ds\\
&\ +C\int_{0}^{t}E\[\widetilde{E}\[\sup_{r\in[0,s]}|\widetilde{X}^{1,n-1}_{r}|^{2} \Bigr\{\int_{0}^{T}|(\partial_{\mu}\sigma)_{1}(s,X_{s},P_{(X_{s\wedge\cdot},u)}; (\widetilde{X}_{s\wedge\cdot}, \widetilde{u}))(r)|dr\Bigr\}^{2}\]ds\]\\
\leq&\ C\int_{0}^{t}E[\sup_{s\in[0,t]}|\Delta X^{1,n-1}_{s}|^{2}]ds,\ t\in[0,T].
\end{aligned}
\end{equation*}
% last inequality is due to the definition of the inner product and H\"{o}lder inequality, meanwhile we noticed that $C\in\mathbb{R}_{+}$ is a constant which do not depends on $T$ and $n$. We choose $T\leq\delta:=\frac{1}{4C}$, then we obtain,
%\begin{equation*}
%\begin{aligned}
%&E[\sup_{0\in[0,T]}|\Delta X_{t}^{1,n}|^{2}]\leq \frac{1}{4^{n-1}} E[\sup_{0\in[0,T]}|\Delta X_{t}^{1,1}|^{2}]\leq \frac{C}{4^{n}}.
%\end{aligned}
%\end{equation*}
%Therefore, for $m\geq n$,
%\begin{equation*}
%\begin{aligned}
%&E[\sup_{0\in[0,T]}|X_{t}^{1,m}-X_{t}^{1,n}|^{2}]\leq \sum_{i=n+1}^{m}E[\sup_{0\in[0,T]}\Delta| X_{t}^{1,i+1}|]\leq\sum_{i=n+1}^{m}\frac{C}{2^{i}}\leq\frac{C}{2^{n}}\rightarrow  0,\ \mbox{as}\ n\rightarrow\infty.
%\end{aligned}
%\end{equation*}
From here it is standard to show that there exists $X\in S^{2}_{\mathbb{F}}$ on $[0,T]$ such that
\begin{equation}\label{eq_5.6}
\begin{aligned}
&E[\sup_{t\in[0,T]}|X_{t}^{1,n}-X_{t}|^{2}]\rightarrow  0,\ \mbox{as}\ n\rightarrow\infty.
\end{aligned}
\end{equation}
Passing to the limit $n\rightarrow\infty$ in \eqref{eq_5.4} we see that $X$ solves \eqref{eq_variSDE}.
%Now the remaining problem is to prove the existence of a solution for any positive constant $T$. In fact, after we prove %the existence of the solution on $[0,\delta]$, we consider the solution on the interval $[\delta,2\delta]$, repeat the %proof steps from \eqref{eq_5.4} to \eqref{eq_5.6}, we will get the result there exists a solution on  %$[\delta,2\delta]$. For any $T$, we always find $T\leq l\delta$ with $l$ is large enough. Therefore we prove for any %positive constant $T$ there exists a solution $X$ solve the equation \eqref{eq_variSDE}. It is worth noting that the %initial value of the equation on $[k\delta,(k+1)\delta]$ is the terminal value of the equation on $[(k-1)\delta, %k\delta]$. We omit the details of the proof.

Now, after proving the existence and the uniqueness of the solution for the forward equation \eqref{eq_variSDE}, we consider the backward equation in the variational system (\ref{eq_vari}):
\begin{equation}\label{eq_variBSDE}
\left\{
\begin{aligned}
dY^1_{t}=&\ -\{f_x(\theta_{t})X^1_{t}+f_y(\theta_{t})Y^1_{t}+f_z(\theta_{t})Z^1_{t}\\%+f_v(\theta_{t},u_t)v_{t}\\
         &\ +\widetilde{E}[\langle(\partial_{\mu}f)(\theta_{t}; (\widetilde{X}, \widetilde{Y}_{\cdot\vee t},\widetilde{u})),\ (\widetilde{X}^{1}, \widetilde{Y}_{\cdot\vee t}^{1},\widetilde{v})\rangle_{2}]\}dt +Z^1_{t}dB_t,\\
Y^1_{T}=&\ \Phi_{x}(X_{T},P_{(X,u)})X^1_{T}+\widetilde{E}[\langle(\partial_{\mu}\Phi)(X_{T},P_{(X,u)}; (\widetilde{X},\widetilde{u})),\ (\widetilde{X}^{1},\widetilde{v})\rangle_{1}].
\end{aligned}
\right.
\end{equation}
%%%%%%%%%%%%%%%%%%%%%%%%%%%%%%%%%%%%%%%%%%%%%%%%%%%%%%%%%%%%%%%%%%
To prove the existence and the uniqueness of the solution for \eqref{eq_variBSDE}, we let $(\xi,\eta)\in S^{2}_{\mathbb{F}}\times M^{2}_{\mathbb{F}}$ and denote by $(\overline{Y},\overline{Z})\in S^{2}_{\mathbb{F}}\times M^{2}_{\mathbb{F}}$ the solution of

\begin{equation}
\left\{
\begin{aligned}
d\overline{Y}_{t}=&-\{f_x(\theta_{t})X^1_{t}+f_y(\theta_{t})\xi_{t}+f_z(\theta_{t})\eta_{t}\\%+f_v(\theta_{t},u_t)v_{t}\\
         &+\widetilde{E}[\langle(\partial_{\mu}f)(\theta_{t}; (\widetilde{X}, \widetilde{Y}_{\cdot\vee t},\widetilde{u})),\ (\widetilde{X}^{1}, \widetilde{\xi}_{\cdot\vee t},\widetilde{v})\rangle_{2}]\}dt +\overline{Z}_{t}dB_t,\\
\overline{Y}_{T}=&\Phi_{x}(X_{T},P_{(X,u)})X^1_{T}+\widetilde{E}[\langle(\partial_{\mu}\Phi)(X_{T},P_{(X,u)}; (\widetilde{X},\widetilde{u})),\ (\widetilde{X}^{1},\widetilde{v})\rangle_{1}].
\end{aligned}
\right.
\end{equation}
We define $\Phi:{S}^{2}_{\mathbb{F}}\times{M}^{2}_{\mathbb{F}}\rightarrow{S}^{2}_{\mathbb{F}}\times{M}^{2}_{\mathbb{F}}$ by putting $\Phi(\xi,\eta)=(\overline{Y},\overline{Z})$. Given any $(\xi^{1},\eta^{1}),\ (\xi^{2},\eta^{2})\in{S}^{2}_{\mathbb{F}}\times{M}^{2}_{\mathbb{F}}$, let $(\overline{Y}^{i},\overline{Z}^{i})=\Phi(\xi^{i},\eta^{i}),\ i=1,2$, and put $(\overline{\xi},\overline{\eta}):=(\xi^{1}-\xi^{2},\eta^{1}-\eta^{2}),\ (\widehat{Y},\widehat{Z}):=(\overline{Y}^{1}-\overline{Y}^{2},\overline{Z}^{1}-\overline{Z}^{2})$.
Then
\begin{equation}
\left\{
\begin{aligned}
d\widehat{Y}_{t}=& -\{f_y(\theta_{t})\overline{\xi}_{t}+f_z(\theta_{t})\overline{\eta}_{t}
         +\widetilde{E}[\langle(\partial_{\mu}f)(\theta_{t}; (\widetilde{X}, \widetilde{Y}_{\cdot\vee t},\widetilde{u})),\ (0, \widetilde{\overline{\xi}}_{\cdot\vee t},0)\rangle_{2}]\}dt +\widehat{Z}_{t}dB_t,\\
\widehat{Y}_{T}=&\ 0.
\end{aligned}
\right.
\end{equation}
Applying It\^{o}'s formula to $e^{\beta t}|\widehat{Y}_t|^{2}$ ($\beta>0$ will be specified later) yields
\begin{equation}\label{Eq_Ito_Y^2}
\begin{aligned}
&\ e^{\beta t}|\widehat{Y}_{s}|^{2}+\int_{t}^{T}e^{\beta s}|\widehat{Z}_{s}|^{2}ds+\beta\int_{t}^{T}e^{\beta s}|\widehat{Y}_{s}|^{2}ds\\
=&\ 2\int_{t}^{T}e^{\beta s}\widehat{Y}_{s}\(f_{y}(\theta_{s})\overline{\xi}_{s}+f_z(\theta_{s})\overline{\eta}_{s}
+\widetilde{E}[\langle(\partial_{\mu}f)(\theta_{s}; (\widetilde{X}, \widetilde{Y}_{\cdot\vee s},\widetilde{u})),\ (0, \widetilde{\overline{\xi}}_{\cdot\vee s},0)\rangle_{2}]\)ds\\
&\ -\int_{t}^{T}e^{\beta s}\widehat{Y}_{s}\widehat{Z}_{s}dB_s,\ t\in[0,T].
\end{aligned}
\end{equation}
As $(\widehat{Y},\widehat{Z})\in S_{\mathbb{F}}^{2}\times M_{\mathbb{F}}^{2}$, it is standard to show the process $\int_{0}^{\cdot}e^{\beta s}\widehat{Y}_{s}\widehat{Z}_{s}dB_{s}$ is a uniformly integrable martingale. Consequently, with $C\in\mathbb{R}_{+}$ as bound of the derivatives of $f$, we have
\begin{equation}\label{eq_Ito_estimate1}
\begin{aligned}
&\ \beta E\[\int_{0}^{T}e^{\beta s}|\widehat{Y}_{s}|^{2}ds\]+E\[\int_{0}^{T}e^{\beta s}|\widehat{Z}_{s}|^{2}ds\]\\
\leq&\ 2CE\[\int_{0}^{T}e^{\beta s}|\widehat{Y}_{s}|(|\overline{\xi}_{s}|+|\overline{\eta}_{s}|)ds\]
         +2CE\[\int_{0}^{T}e^{\frac{\beta s}{2}}|\widehat{Y}_{s}|\widetilde{E}\[\int_{0}^{T}e^{\frac{\beta s}{2}}|\widetilde{\overline{\xi}}_{r\vee s}|dr\]ds\].
\end{aligned}
\end{equation}
But as $\beta>0$,
\begin{equation}\label{eq_Ito_estimate2}
\begin{aligned}
\widetilde{E}\[\int_{0}^{T}e^{\frac{\beta s}{2}}|\widetilde{\overline{\xi}}_{r\vee s}|dr\]
\leq E\[\int_{0}^{T}e^{\frac{\beta (r\vee s)}{2}}|\overline{\xi}_{r\vee s}|dr\]
\leq TE\[\sup_{r\in[0,T]}(e^{\frac{\beta r}{2}}|\overline{\xi}_{r}|)\],
\end{aligned}
\end{equation}
we get from \eqref{eq_Ito_estimate1}, for some $\delta>0$ (to specify later),
\begin{equation}\label{eq_Ito_estimate3}
\begin{aligned}
&\beta E\[\int_{0}^{T}e^{\beta s}|\widehat{Y}_{s}|^{2}ds\]+E\[\int_{0}^{T}e^{\beta s}|\widehat{Z}_{s}|^{2}ds\]\\
\leq&\ C_{\delta}E\[\int_{0}^{T}e^{\beta s}|\widehat{Y}_{s}|^{2}ds\]
+\delta \(E\[\int_{0}^{T}e^{\beta s}|\overline{\xi}_{s}|^{2}ds\]+ E\[\int_{0}^{T}e^{\beta s}|\overline{\eta}_{s}|^{2}ds\]+ T^{2}E\[\sup_{r\in[0,T]}(e^{\frac{\beta r}{2}}|\overline{\xi}_{r}|)^{2}\]\)\\
\leq& C_{\delta}E\[\int_{0}^{T}e^{\beta s}|\widehat{Y}_{s}|^{2}ds\]+\delta T(1+T)E\[\sup_{r\in[0,T]}(e^{\beta r}|\overline{\xi}_{r}|^{2})\]+\delta E\[\int_{0}^{T}e^{\beta s}|\overline{\eta}_{s}|^{2}ds\].
        \end{aligned}
\end{equation}
On the other hand, with the help of the Burkholder-Davis-Gundy inequality, we deduce from \eqref{Eq_Ito_Y^2}
that

\begin{equation*}%\label{eq_Ito_estimate4}
\begin{aligned}
& E\[\sup_{t\in[0,T]}e^{\beta t}|\widehat{Y}_{t}|^{2}\]\\
\leq&\  2CE\[\int_{0}^{T}e^{\frac{\beta s}{2}}|\widehat{Y}_{s}|(e^{\frac{\beta s}{2}}|\overline{\xi}_{s}|+e^{\frac{\beta s}{2}}|\overline{\eta}_{s}|+ TE[\sup_{r\in[0,T]}e^{\frac{\beta r}{2}}|\overline{\xi}_{r}|])ds\]
+C_{1}E\[\(\int_{0}^{T}(e^{\beta s}|\widehat{Y}_{s}|^{2})(e^{\beta s}|\widehat{Z}_{s}|^{2})ds\)^{\frac{1}{2}}\]\\
\leq&\ C_{\delta}E\[\int_{0}^{T}e^{\beta s}|\widehat{Y}_{s}|ds\]
+\delta \( E\[\int_{0}^{T}e^{\beta s}|\overline{\xi}_{s}|^{2}ds\]+ E\[\int_{0}^{T}e^{\beta s}|\overline{\eta}_{s}|^{2}ds\]
+T^{2} E[\sup_{r\in[0,T]}(e^{\beta r}|\overline{\xi}_{r}|^{2})]\)\\
&+C_{1}^{2}E\[\int_{0}^{T}e^{\beta s}|\widehat{Z}_{s}|^{2}ds\]+\frac{1}{2} E\[\sup_{t\in[0,T]}e^{\beta t}|\widehat{Y}_{t}|^{2}\],
\end{aligned}
\end{equation*}
then we have
\begin{equation}\label{eq_Ito_estimate5}
\begin{aligned}
 E\[\sup_{t\in[0,T]}e^{\beta t}|\widehat{Y}_{t}|^{2}\]
\leq&\ 2C_{\delta}E\[\int_{0}^{T}e^{\beta s}|\widehat{Y}_{s}|ds\]
+2C_{1}^{2}E\[\int_{0}^{T}e^{\beta s}|\widehat{Z}_{s}|^{2}ds\]\\
&\ +2\delta \(T(1+T) E[\sup_{r\in[0,T]}(e^{\beta r}|\overline{\xi}_{r}|^{2})]+ E\[\int_{0}^{T}e^{\beta s}|\overline{\eta}_{s}|^{2}ds\]\).
\end{aligned}
\end{equation}
Consequently, from  \eqref{eq_Ito_estimate3} and \eqref{eq_Ito_estimate5} for some small $\rho>0$ (to be specified),
\begin{equation*}%\label{eq_Ito_estimate5}
\begin{aligned}
&\beta E\[\int_{0}^{T}e^{\beta s}|\widehat{Y}_{s}|^{2}ds\]+\rho E\[\sup_{t\in[0,T]}e^{\beta t}|\widehat{Y}_{t}|^{2}\]
+E\[\int_{0}^{T}e^{\beta s}|\widehat{Z}_{s}|ds\]\\
\leq& C_{\delta}(1+2\rho)E\[\int_{0}^{T}e^{\beta s}|\widehat{Y}_{s}|^{2}ds\]+2\rho C^{2}_{1}E\[\int_{0}^{T}e^{\beta s}|\widehat{Z}_{s}|ds\]+\delta T(1+T)(1+2\rho)E[\sup_{r\in[0,T]}(e^{\beta r}|\overline{\xi}_{r}|^{2})]\\
&+\delta (1+2\rho)E\[\int_{0}^{T}e^{\beta s}|\overline{\eta}|^{2}_{s}ds\].
\end{aligned}
\end{equation*}
Choose now $\rho=\frac{1}{4C_{1}^{2}}(\leq 1),\ \delta=\frac{\rho}{12}(\frac{1}{T(1+T)}\wedge 1)$ and $\beta=3C_{\delta}$.
Then
\begin{equation}\label{eq_Ito_estimate6}
\begin{aligned}
\rho E\[\sup_{t\in[0,T]}e^{\beta t}|\widehat{Y}_{t}|^{2}\]+\frac{1}{2}E\[\int_{0}^{T}e^{\beta s}|\widehat{Z}_{s}|ds\]
\leq \frac{\rho}{2}E[\sup_{t\in[0,T]}(e^{\beta t}|\overline{\xi}_{t}|^{2})]+\frac{1}{4}E\[\int_{0}^{T}e^{\beta s}|\overline{\eta}_{s}|^{2}ds\].
\end{aligned}
\end{equation}
This proves that for ${S}^{2}_{\mathbb{F}}\times{M}^{2}_{\mathbb{F}}$ endowed with the norm
$$||(\xi,\eta)||_{\rho}:=\(\rho E[\sup_{t\in[0,T]}(e^{\beta t}|\overline{\xi}_{t}|^{2})]+\frac{1}{2}E\[\int_{0}^{T}e^{\beta s}|\overline{\eta}_{s}|^{2}ds\]\)^{\frac{1}{2}},$$
$(\xi,\eta)\in{S}^{2}_{\mathbb{F}}\times{M}^{2}_{\mathbb{F}}$, the mapping
$\Phi:{S}^{2}_{\mathbb{F}}\times{M}^{2}_{\mathbb{F}}\rightarrow{S}^{2}_{\mathbb{F}}\times{M}^{2}_{\mathbb{F}}$
is a contraction, and its unique fixed point is the unique solution of \eqref{eq_variBSDE}.
%%%%%%%%%%%%%%%%%%%%%%%%%%%%%%%%%%%%%%%%%%%%%%%%%%%%%%%%%%%%%%%%%%
\end{proof}

%\begin{remark}\label{rem5.1}
%When $l= f,\,\sigma,\,L$, respectively, $l_x$ is the partial derivative of $l(t,x,y,z,\mu,v)$ with respect to $x$; $l_{\mu}$ is the partial derivative of $l(t,x,y,z,\mu,v)$ with respect to $\mu$. Similar to $l_y,\,l_z,\,l_v$.
%\end{remark}

We denote by $(X^{\rho},\ Y^{\rho},\ Z^{\rho}):=(X^{u^\rho},\ Y^{u^\rho},\ Z^{u^\rho})$\ the trajectory associated with the control $u^{\rho}$. Then we have the following estimates.

\begin{lemma}\label{lem5.1.1}
We assume \textnormal{(H5)} and \textnormal{(H6)} hold true. Then for $v\in \mathcal{U}_{ad}$ and $p\geq2$, there is a constant $C_{p}\in\mathbb{R}_{+}$ such that
\begin{equation}\label{SDEestimates1_1}
E\Bigr[\sup_{t\in[0,T]}\Bigr|\frac{1}{\rho}(X^{\rho}_{t}-X_{t})\Bigr|^{p}\Bigr]\leq C_{p}E\Bigr[\Bigr(\int_{0}^{T}|v_t|^{2}dt\Bigr)^{\frac{p}{2}}\Bigr]
\end{equation}
and
\begin{equation}\label{Eq21-2-2_1}
E\Bigr[\sup_{s\in[0,T]}\Bigr(\frac{1}{\rho}|Y_s^\rho-Y_{s}|\Bigr)^p\Bigr]
+E\Bigr[\Bigr(\int_0^T\Bigr(\frac{1}{\rho}|Z_s^\rho-Z_{s}|\Bigr)^2ds\Bigr)^{\frac{p}{2}}\Bigr]\leq C_pE\Bigr[\Bigr(\int_{0}^{T}|v_t|^{2}dt\Bigr)^{\frac{p}{2}}\Bigr],
\end{equation}
for all $\rho\in(0,1)$.
\end{lemma}
\begin{proof}
We first prove estimate \eqref{SDEestimates1_1}.
Recall that
\begin{equation}\label{Eqcontrol1}
\left\{\begin{aligned}
dX^{\rho}_t&=\sigma(t,X^{\rho}_t,P_{(X^{\rho}_{\cdot\wedge t},u^{\rho})})dB_t,\ X^{\rho}_{0}=x\in\mathbb{R},\\
dY^{\rho}_t&=-f(t,X^{\rho}_t,Y^{\rho}_t,Z^{\rho}_t,P_{(X^{\rho},Y^{\rho}_{\cdot\vee t},u^{\rho})})dt+Z^{\rho}_tdB_t,\ t\in[0,T),\\
Y^{\rho}_T&=\Phi(X^{\rho}_T,P_{(X^{\rho},u^{\rho})}),\\
\end{aligned}\right.
\end{equation}
and $(X,Y,Z)$ is the solution related with the optimal control $u$. We have
\begin{equation}\label{eq_5.23_1}
\begin{aligned}
\frac{1}{\rho}(X_{t}^{\rho}-X_{t})&=\frac{1}{\rho}\int_{0}^{t}(\sigma(s,X^{\rho}_{s},P_{(X^{\rho}_{\cdot\wedge s},u^{\rho})})-\sigma(t,X_{s},P_{(X_{\cdot\wedge s},u)}))dB_{s},\ t\in[0,T].
\end{aligned}
\end{equation}
Let us use the following notations:
\begin{equation*}
\begin{aligned}
\vartheta_{t}^{\lambda,\rho}:=(X_{\cdot\wedge t}+\lambda(X^{\rho}_{\cdot\wedge t}-X_{\cdot\wedge t}), u+\lambda\rho v),\
X_{t}^{\lambda,\rho}:=X_{t}+\lambda(X_{t}^{\rho}-X_{t}),\ \lambda\in[0,1].
%u_{t}^{\lambda,\rho}:=u_{t}+\lambda\rho v_{t}.
\end{aligned}
\end{equation*}
For the coefficient $\sigma$ we have
\begin{equation}\label{eq_Chain1}
\begin{aligned}
&\frac{1}{\rho}\[\sigma(t,X^{\rho}_{t},P_{(X^{\rho}_{\cdot\wedge t},u^{\rho}) } )-\sigma(t,X_{t},P_{(X_{\cdot\wedge t},u)})\]
= \frac{1}{\rho}\int_{0}^{1}\partial_{\lambda}[\sigma(t,X^{\lambda,\rho}_{t},P_{\vartheta_{t}^{\lambda,\rho}})] d\lambda\\
=&\ \!\!\int_{0}^{1}\!\!\sigma_{x}(t,X^{\lambda,\rho}_{t},P_{\vartheta_{t}^{\lambda,\rho}})\frac{1}{\rho}(X^{\rho}_{t}-X_{t})d\lambda
%+\int_{0}^{1}\sigma_{v}(t,X^{\lambda,\rho}_{t},P_{\vartheta_{t}^{\lambda,\rho}},u_{t}^{\lambda\rho})v_{t}d\lambda\\
\!+\!\int_{0}^{1}\!\widetilde{E}[\langle(\partial_{\mu}\sigma)(t,X^{\lambda,\rho}_{t},P_{\vartheta_{t}^{\lambda,\rho}};\widetilde{\vartheta}_{t}^{\lambda,\rho}) , (\frac{1}{\rho}(\widetilde{X}^{\rho}_{\cdot\wedge t}-\widetilde{X}_{\cdot\wedge t})),\widetilde{v} ) \rangle_{1}]d\lambda.
\end{aligned}
\end{equation}
Thus, for $p\geq 2$, we get from SDE \eqref{eq_5.23_1}
\begin{equation*}
\begin{aligned}
E\[\sup_{s\in[0,t]}\|\frac{1}{\rho}(X^{\rho}_{s}-X_{s})\|^{p}\]
%=&\ E[\sup_{t\in[0,T]}|\frac{1}{\rho}\int_{0}^{t}(\sigma(s,X^{\rho}_{s},P_{(X^{\rho}_{\cdot\wedge s},u^{\rho}) },u_{s}^{\rho} )-\sigma(t,X_{s},P_{(X_{\cdot\wedge s},u)},u_{s}))dB_{s}|^{p}]\\
\leq&\ C_{p}E\[\Bigr\{\int_{0}^{t}\|\frac{1}{\rho}(\sigma(s,X^{\rho}_{s},P_{(X^{\rho}_{\cdot\wedge s},u^{\rho}) } )-\sigma(s,X_{s},P_{(X_{\cdot\wedge s},u)}))\|^2ds\Bigr\}^{\frac{p}{2}}\]\\
%\leq&\ C_{p}E\[\Bigr\{\int_{0}^{t}\(\|\frac{1}{\rho}(X^{\rho}_{s}-X_{s})\|^{2}+\frac{1}{\rho^{2}}W^{2}_{2}(P_{(X^{\rho}_{\cdot\wedge s},u^{\rho})},P_{(X_{\cdot\wedge s},u)})\)ds\Bigr\}^{\frac{p}{2}}\]\\
\leq&\  C_{p}E\[\Bigr\{\int_{0}^{t}\(\|\frac{1}{\rho}(X^{\rho}_{s}-X_{s})\|^{2}
+E\[\sup_{s\in[0,t]}|\frac{1}{\rho}(X^{\rho}_{s}-X_{s})|^{2}+\int_{0}^{T}|v_{s}|^{2}ds\]\)ds\Bigr\}^{\frac{p}{2}}\].
\end{aligned}
\end{equation*}
Consequently, combining H\"{o}lder's inequality and Gronwall's inequality we obtain
\begin{equation*}
\begin{aligned}
E\[\sup_{t\in[0,T]}\|\frac{1}{\rho}(X^{\rho}_{t}-X_{t})\|^{p}\]
%\leq
%C\int_{0}^{T}E[\sup_{r\in[0,t]}|\frac{1}{\rho}(X^{\rho}_{r}-X_{r})|^{p}]dt+
\leq C_{p}E\[\(\int_{0}^{T}|v_{t}|^{2}dt\)^{\frac{p}{2}}\].
\end{aligned}
\end{equation*}
%We finally obtain,
%\begin{equation*}
%\begin{aligned}
%E[\sup_{t\in[0,T]}|\frac{1}{\rho}(X^{\rho}_{t}-X_{t})|^{p}]\leq C\in\mathbb{R}_{+}.
%\end{aligned}
%\end{equation*}

Next we prove inequality \eqref{Eq21-2-2_1}.
From the BSDE in \eqref{Eqcontrol1} we have
\begin{equation}\label{Eq_5.33.1}
\begin{aligned}
\frac{1}{\rho}(Y_{t}^{\rho}-Y_{t})=\frac{1}{\rho}(\Phi(X^{\rho}_T,P_{(X^{\rho},u^{\rho})})-\Phi(X_{T},P_{(X,u)}))
+\int_{t}^{T}\frac{1}{\rho}(f(\theta_{s}^{\rho})-f(\theta_{s}))ds-\int_{t}^{T}\frac{1}{\rho}(Z_{s}^{\rho}-Z_{s})dB_s,
\end{aligned}
\end{equation}
%In order to make the corresponding proof for the associated backward equation,
$t\in[0,T]$. Let us introduce the following notations:
\begin{equation*}
\begin{aligned}
\theta_{t}:=&(t,X_t,Y_t,Z_t,P_{(X,Y_{\cdot\vee t},u)}),\ \theta_{t}^{\rho}:=(t,X^{\rho}_t,Y^{\rho}_t,Z^{\rho}_t,P_{(X^{\rho},Y^{\rho}_{\cdot\vee t},u^{\rho})}),\\
Y_{t}^{\lambda,\rho}:=&Y_{t}+\lambda(Y_{t}^{\rho}-Y_{t}),\quad\qquad Z_{t}^{\lambda,\rho}:=Z_{t}+\lambda(Z_{t}^{\rho}-Z_{t}),\\
%u_{t}^{\lambda\rho}:=u_{t}+\lambda\rho v_{t}.\\
\theta_{t}^{\lambda,\rho}:=&(t,X_t^{\lambda,\rho},Y_t^{\lambda,\rho},
Z_t^{\lambda,\rho},P_{(X^{\lambda,\rho},Y_{\cdot\vee t}^{\lambda,\rho},u^{\lambda\rho})}).\\
\end{aligned}
\end{equation*}
Then, in a similar way to $\sigma$ in \eqref{eq_Chain1}, for the driving coefficient $f$ of the BSDE we have
\begin{equation}\label{eq_Chain2}
\begin{aligned}
&\frac{1}{\rho}[f(\theta^{\rho}_{t})-f(\theta_{t})]\\
=&\int_{0}^{1}f_{x}(\theta^{\lambda,\rho}_{t})\frac{1}{\rho}(X^{\rho}_{t}-X_{t})d\lambda
  +\int_{0}^{1}f_{y}(\theta^{\lambda,\rho}_{t})\frac{1}{\rho}(Y^{\rho}_{t}-Y_{t})d\lambda
+\int_{0}^{1}f_{z}(\theta^{\lambda,\rho}_{t})\frac{1}{\rho}(Z^{\rho}_{t}-Z_{t})d\lambda\\
  %+\int_{0}^{1}f_{v}(\theta^{\lambda,\rho}_{t},u_{t}^{\lambda\rho})v_{t}d\lambda\\
&+\int_{0}^{1}\widetilde{E}[\langle(\partial_{\mu}f)(\theta^{\lambda,\rho}_{t};
(\widetilde{X}^{\lambda,\rho},\widetilde{Y}_{\cdot\vee t}^{\lambda,\rho},\widetilde{u}^{\lambda,\rho})) , (\frac{1}{\rho}(\widetilde{X}^{\rho}-\widetilde{X}), \frac{1}{\rho}(\widetilde{Y}_{\cdot\vee t}^{\rho}-\widetilde{Y}_{\cdot\vee t}),\widetilde{v}) \rangle_{2}]d\lambda.
\end{aligned}
\end{equation}
Recall that all derivatives of $f$ are supposed to be bounded.
Applying a well-known BSDE standard estimate to \eqref{Eq_5.33.1}, we get that, for all $p\geq2$, there exists $C_p\in\mathbb{R}_+$ such that, for all $0\leq t\leq t'\leq T$,
\begin{equation}\label{Eq21-1-1}
\begin{aligned}
&\ E\[\sup_{s\in[t,t']}\Bigr|\frac{1}{\rho}\(Y_s^\rho-Y_s\)\Bigr|^p
+\(\int_t^{t'}\Bigr|\frac{1}{\rho}(Z_s^\rho-Z_s)\Bigr|^2ds\)^{\frac{p}{2}}\]
\leq C_p E\[\Bigr|\frac{1}{\rho}(Y_{t'}^\rho-Y_{t'})\Bigr|^p
+\(\int_t^{t'}\frac{1}{\rho}|f(\theta_s^\rho)-f(\theta_s)|ds\)^p\].
\end{aligned}
\end{equation}
We observe that, due to our assumptions (H5) and (H6)-(ii) on $f$,
\begin{equation}\notag
\begin{aligned}
I_s^{\lambda,\rho}:=  &\ \widetilde{E}\[\Bigr\langle(\partial_\mu f)
(\theta_s^{\lambda,\rho};(\widetilde{X}^{\lambda,\rho},\widetilde{Y}^{\lambda,\rho}_{\cdot\vee s},\widetilde{u}^{\lambda,\rho})),
(\frac{1}{\rho}(\widetilde{X}^\rho-\widetilde{X}),\frac{1}{\rho}(\widetilde{Y}^\rho_{\cdot\vee s}-\widetilde{Y}_{\cdot\vee s}),\widetilde{v})\Bigr\rangle_2\]\\
= &\ \widetilde{E}\[\int_0^T(\partial_\mu f)_1
(\theta_s^{\lambda,\rho};(\widetilde{X}^{\lambda,\rho},\widetilde{Y}^{\lambda,\rho}_{\cdot\vee s},\widetilde{u}^{\lambda,\rho}))(r)\frac{1}{\rho}(\widetilde{X}^\rho_r-\widetilde{X}_r)dr\]\\
&\ +\widetilde{E}\[\int_0^T(\partial_\mu f)_2
(\theta_s^{\lambda,\rho};(\widetilde{X}^{\lambda,\rho},\widetilde{Y}^{\lambda,\rho}_{\cdot\vee s},\widetilde{u}^{\lambda,\rho}))(r)\frac{1}{\rho}(\widetilde{Y}^\rho_{r\vee s}-\widetilde{Y}_{r\vee s})dr\]\\
&+\widetilde{E}\[\int_0^T(\partial_\mu f)_3
(\theta_s^{\lambda,\rho};(\widetilde{X}^{\lambda,\rho},\widetilde{Y}^{\lambda,\rho}_{\cdot\vee s},\widetilde{u}^{\lambda,\rho}))(r)\widetilde{v}_r dr\],
\end{aligned}
\end{equation}
and
\begin{equation}\notag
|I_s^{\lambda,\rho}|\leq C\widetilde{E}\[\int_0^T\(\|\frac{1}{\rho}(\widetilde{X}^\rho_r-\widetilde{X}_r)\|
+\|\frac{1}{\rho}(\widetilde{Y}^\rho_{r\vee s}-\widetilde{Y}_{r\vee s})\|+\|\widetilde{v}_r\|\)dr\],\ s\in[0,T],\ \lambda,\rho\in(0,1).
\end{equation}
Then, from \eqref{Eq21-1-1} combined with \eqref{eq_Chain2} we conclude
\begin{equation}\label{Eq21-1-2A}
\begin{aligned}
&\ E\[\sup_{s\in [t,t']}\|\frac{1}{\rho}(Y_s^\rho-Y_s\)\|^p
+\(\int_t^{t'}\Bigr|\frac{1}{\rho}(Z_s^\rho-Z_s)\Bigr|^2ds\)^{\frac{p}{2}}\]\\
\leq &\ C_p E\[\(\frac{1}{\rho}|Y_{t'}^\rho-Y_{t'}|\)^p\]
+C_p E\Bigr[\Bigr(\int_t^{t'}\Bigr\{\frac{1}{\rho}|X_{s}^\rho-X_{s}|
+\frac{1}{\rho}|Y_{s}^\rho-Y_{s}|+\frac{1}{\rho}|Z_{s}^\rho-Z_{s}|\\%+|\widetilde{\vartheta}_s|
&\ +\widetilde{E}\Bigr[\int_0^{T}\(|\frac{1}{\rho}(\widetilde{X}^\rho_r-\widetilde{X}_r)|
+|\frac{1}{\rho}(\widetilde{Y}^\rho_{r\vee s}-\widetilde{Y}_{r\vee s})|+|\widetilde{v}_r|\)dr\Bigr]\Bigr\}ds\Bigr)^p\Bigr]\\
\leq &\ C_p E\[\(\frac{1}{\rho}|Y_{t'}^\rho-Y_{t'}|\)^p\]
+C_p|t'-t|^{\frac{p}{2}}\Bigr\{E\[\sup_{s\in[0,T]}\bigr|\frac{1}{\rho}(X_s^\rho-X_s)\bigr|^p\]
+E\[\sup_{s\in[t,T]}\bigr|\frac{1}{\rho}(Y_s^\rho-Y_s)\bigr|^p\]\\
\end{aligned}
\end{equation}
\begin{equation}\notag
\begin{aligned}
&\ +E\Bigr[\Bigr(\int_t^{t'}\Bigr(\frac{1}{\rho}|Z_s^\rho-Z_s|\Bigr)^2ds\Bigr)^{\frac{p}{2}}\Bigr]
+E\Bigr[\Bigr(\int_0^{T}|v_r|^2dr\Bigr)^{\frac{p}{2}}\Bigr]\Bigr\}.
\end{aligned}
\end{equation}
Let now $\delta>0$ be such that $C_p\delta^{\frac{p}{2}}\leq \frac{1}{2}$. Then, for $t'=T$ and $t=T-\delta$,
\begin{equation}\label{Eq21-1-2}
\begin{aligned}
&E\[\sup_{s\in [T-\delta,T]}\|\frac{1}{\rho}\(Y_s^\rho-Y_s\)\|^p\]
+E\[\(\int_{T-\delta}^{T}\|\frac{1}{\rho}\(Z_s^\rho-Z_s\)\|^2ds\)^{\frac{p}{2}}\]\\
\leq\ &2C_p E\[\|\frac{1}{\rho}(Y_{T}^\rho-Y_{T})\|^p\]
+E\[\sup_{s\in[0,T]}\|\frac{1}{\rho}(X_s^\rho-X_s)\|^p\]
+E\[\(\int_0^{T}|v_r|^2dr\)^{\frac{p}{2}}\].
\end{aligned}
\end{equation}
But, as
\begin{equation}\notag
\begin{aligned}
E\[\|\frac{1}{\rho}(Y_T^\rho-Y_{T})\|^p\]
=&\ E\[\(\frac{1}{\rho}\bigr|\Phi(X_T^\rho,P_{(X^\rho,u^\rho)})-\Phi(X_T,P_{(X,u)})\bigr|\)^p\]\\
\leq&\ C_p E\[\sup_{s\in[0,T]}\|\frac{1}{\rho}(X_s^\rho-X_s)\bigr|^p\]
+C_p E\[\(\int_0^{T}|v_r|^2dr\)^{\frac{p}{2}}\],
\end{aligned}
\end{equation}
estimate \eqref{Eq21-1-2} yields that, for some constant $C_p\in\mathbb{R}_+$,
\begin{equation}\notag
\begin{aligned}
&E\[\sup_{s\in[T-\delta,T]}\|\frac{1}{\rho}(Y_s^\rho-Y_{s})\|^p\]
+E\[\(\int_{T-\delta}^T\|\frac{1}{\rho}(Z_s^\rho-Z_{s})\|^2ds\)^{\frac{p}{2}}\]\\
\leq&\ C_p \(E\[\sup_{s\in[0,T]}\|\frac{1}{\rho}(X_s^\rho-X_s)\|^p\]
+ E\[\(\int_0^{T}|v_r|^2dr\)^{\frac{p}{2}}\]\).
\end{aligned}
\end{equation}
Consequently, from \eqref{SDEestimates1_1} we obtain
\begin{equation}\label{Eq21-2-2D}
E\left[\sup_{s\in[T-\delta,T]}\left|\frac{1}{\rho}(Y_s^\rho-Y_{s})\right|^p\right]
+E\left[\left(\int_{T-\delta}^T\left|\frac{1}{\rho}(Z_s^\rho-Z_{s})\right|^2ds\right)^{\frac{p}{2}}\right]\leq C_p E\Bigr[\Bigr(\int_0^{T}|v_r|^2dr\Bigr)^{\frac{p}{2}}\Bigr],\ \rho\in(0,1).
\end{equation}
Next let us take $t'=T-\delta$ and $t=T-2\delta.$ Then, from (\ref{Eq21-1-2A})
\begin{equation}\label{Eq21-1-2B}
\begin{aligned}
&\ E\[\sup_{s\in [T-2\delta,T-\delta]}\|\frac{1}{\rho}(Y_s^\rho-Y_s)\|^p
+\(\int_{T-2\delta}^{T-\delta}\Bigr|\frac{1}{\rho}(Z_s^\rho-Z_s)\Bigr|^2ds\)^{\frac{p}{2}}\]\\
\leq &\ C_p E\[\|\frac{1}{\rho}(Y_{T-\delta}^\rho-Y_{T-\delta})\|^p\]
+C_p|t'-t|^{\frac{p}{2}}\Bigr\{E\[\sup_{s\in[0,T]}\bigr|\frac{1}{\rho}(X_s^\rho-X_s)\bigr|^p\]
		+E\[\sup_{s\in[T-2\delta,T-\delta]}\bigr|\frac{1}{\rho}(Y_s^\rho-Y_s)\bigr|^p\]\\
&\ +E\[\sup_{s\in[T-\delta,T]}\bigr|\frac{1}{\rho}(Y_s^\rho-Y_s)\bigr|^p\]+E\Bigr[\Bigr(\int_{T-2\delta}^{T-\delta}\Bigr|\frac{1}{\rho}(Z_s^\rho-Z_s)\Bigr|^2ds\Bigr)^{\frac{p}{2}}\Bigr]
		+E\Bigr[\Bigr(\int_0^{T}|v_r|^2dr\Bigr)^{\frac{p}{2}}\Bigr]\Bigr\}.
	\end{aligned}
\end{equation}
Then, as $C_p\delta^{p/2}\le 1/2$, (\ref{Eq21-2-2D}) and \eqref{SDEestimates1_1} yield
\begin{equation}\label{Eq21-1-2C}
	E\[\sup_{s\in[T-2\delta,T-\delta]}\|\frac{1}{\rho}(Y_s^\rho-Y_{s})\|^p\]
	+E\[\(\int_{T-2\delta}^{T-\delta}\|\frac{1}{\rho}(Z_s^\rho-Z_{s})\|^2ds\)^{\frac{p}{2}}\]\leq C_p E\[\Bigr(\int_0^{T}|v_r|^2dr\Bigr)^{\frac{p}{2}}\],\ \rho\in(0,1).
\end{equation}
In the same way the estimates (\ref{Eq21-1-2A}),\ (\ref{Eq21-2-2D}) and
(\ref{Eq21-1-2C}) allow to conclude that
\begin{equation}\label{Eq21-1-2E}
	E\[\sup_{s\in[T-3\delta,T-2\delta]}\|\frac{1}{\rho}(Y_s^\rho-Y_{s})\|^p\]
	+E\[\(\int_{T-3\delta}^{T-2\delta}\|\frac{1}{\rho}(Z_s^\rho-Z_{s})\|^2ds\)^{\frac{p}{2}}\]\leq C_p E\Bigr[\Bigr(\int_0^{T}|v_r|^2dr\Bigr)^{\frac{p}{2}}\Bigr],\ \rho\in(0,1).
\end{equation}
Finally, by iterating this argument and taking the sums over the estimates over all subintervals $[(T-\ell\delta)^+,(T-(\ell-1)\delta)^+],\ 1\le \ell\le [T/\delta]+1,$\ we get
\begin{equation}\notag
	E\[\sup_{s\in[0,T]}\|\frac{1}{\rho}(Y_s^\rho-Y_{s})\|^p\]
	+E\[\(\int_{0}^{T}\|\frac{1}{\rho}(Z_s^\rho-Z_{s})\|^2ds\)^{\frac{p}{2}}\]\leq C_p E\[\(\int_0^{T}|v_r|^2dr\)^{\frac{p}{2}}\],\ \rho\in(0,1).
\end{equation}
This proves (\ref{Eq21-2-2_1}).
\end{proof}

In what follows, we will use the notations introduced in the proof of Lemma \ref{lem5.1.1}.
\begin{lemma}\label{lem5.2}
We assume \textnormal{(H5)} and \textnormal{(H6)} hold. Then,
$\lim\limits_{\rho\rightarrow0}\frac{1}{\rho}(X^{\rho}_{t}-X_{t})=X^1_{t}, \  \lim\limits_{\rho\rightarrow0}\frac{1}{\rho}(Y^{\rho}_{t}-Y_{t})=Y^1_{t},$\
$ \lim\limits_{\rho\rightarrow0}\frac{1}{\rho}(Z^{\rho}_{t}-Z_{t})=Z^1_{t},$\ with convegence in $S_{\mathbb{F}}^2\times S_{\mathbb{F}}^2\times M_{\mathbb{F}}^2$.
\end{lemma}
\begin{proof}
Let $\widehat{X}_{t}:=X^{\rho}_{t}-X_{t}, \  \widehat{Y}_{t}:=Y^{\rho}_{t}-Y_{t}, \ \widehat{Z}_{t}:=Z^{\rho}_{t}-Z_{t}$, and recall that $(X,Y,Z)$ is the solution related with the optimal control $u$. Then,
\begin{equation}\label{Eq_control2}
\left\{
\begin{aligned}
\ d\widehat{X}_{t}=&\ [\sigma(t,X^{\rho}_t,P_{(X^{\rho}_{\cdot\wedge t},u^{\rho})})-\sigma(t,X_t,P_{(X_{\cdot\wedge t},u)})]dB_t,\\
-d\widehat{Y}_{t}=&\ [f(t,X^{\rho}_t,Y^{\rho}_t,Z^{\rho}_t,P_{(X^{\rho},Y^{\rho}_{\cdot\vee t},u^{\rho})})-f(t,X_t,Y_t,Z_t,P_{(X,Y_{\cdot\vee t},u)})]dt-\widehat{Z}_{t}dB_t,\\
\widehat{X}_{0}=&\ 0, \ \widehat{Y}_{T}=\Phi(X^{\rho}_T,P_{(X^{\rho},u^{\rho})})-\Phi(X_{T},P_{(X,u)})).
\end{aligned}
\right.
\end{equation}
From Lemma \ref{lem5.1.1} we know that $(\widehat{X},\ \widehat{Y},\ \widehat{Z})$\ converges to $0$ in $S_{\mathbb{F}}^2(0,T)\times S_{\mathbb{F}}^2(0,T)\times M_{\mathbb{F}}^2(0,T)$ as $\rho$ tends to $0$.
%From our Lipschitz condition on $\sigma$, the Burkholder-Davis-Gundy inequality as well as the fact that %$\lim\limits_{\rho\rightarrow0}W_{2}(\mu^{\rho},\mu)=0$ as $\mu^{\rho}\rightarrow\mu$,
%From Lemma \ref{lem5.1.1} we see easily that
%\begin{equation}\label{SDEestimates0}
%\begin{aligned}
%E[\sup_{t\in[0,T]}|X^{\rho}_{t}-X_{t}|^{p}]\rightarrow 0,\ \mbox{as}\ \rho\rightarrow 0,\ \forall\ p\geq 1.
%\end{aligned}
%\end{equation}
%Indeed, from \eqref{SDEestimates1_1} in Lemma \ref{lem5.1.1} that, for $p\geq2$, we have that there is a constant %$C_{p}\in\mathbb{R}_{+}$ such that
%\begin{equation}\label{SDEestimates1}
%E\[\sup_{t\in[0,T]}\Bigr|\frac{1}{\rho}(X^{\rho}_{t}-X_{t})\Bigr|^{p}\]\leq %C_{p}E\[\(\int_{0}^{T}|v|^{2}dt\)^{\frac{p}{2}}\],\ \rho\in(0,1).
%\end{equation}
Moreover, from \eqref{eq_5.23_1} and \eqref{eq_Chain1},
\begin{equation}\label{eq_5.25_1}
\begin{aligned}
\frac{1}{\rho}(X_{t}^{\rho}-X_{t})
=&R_{t}^{1}
%+R_{t}^{2}
+R_{t}^{2}
+\int_{0}^{t}\sigma_{x}(s,X_{s},P_{(X_{\cdot\wedge s},u)})\frac{1}{\rho}(X^{\rho}_{s}-X_{s})dB_{s}\\
%+\int_{0}^{t}\sigma_{v}(s,X_{s},P_{(X_{\cdot\wedge s},u)},u_{s})v(s)dB_{s}\\
&+\int_{0}^{t}\widetilde{E}[\langle
(\partial_{\mu}\sigma)(s,X_{s},P_{(X_{\cdot\wedge s},u)};(\widetilde{X}_{\cdot\wedge s},\widetilde{u})), (\frac{1}{\rho}(\widetilde{X}^{\rho}_{\cdot\wedge s}-\widetilde{X}_{\cdot\wedge s}),\widetilde{v})\rangle_{1}]dB_{s},
\end{aligned}
\end{equation}
where
\begin{equation*}
\begin{aligned}
R_{t}^{1}&=\int_{0}^{t}\Bigr\{\int_{0}^{1}(\sigma_{x}(s,X^{\lambda,\rho}_{s},P_{\vartheta_{s}^{\lambda,\rho}})-\sigma_{x}(s,X_{s},P_{(X_{\cdot\wedge s},u)}))\frac{1}{\rho}(X^{\rho}_{s}-X_{s})d\lambda\Bigr\}dB_{s},\\
%R_{t}^{2}&=\int_{0}^{t}\{\int_{0}^{1}(\sigma_{v}(s,X^{\lambda,\rho}_{s},P_{\vartheta_{s}^{\lambda,\rho}},u_{s}^{\lambda,\rho})-\sigma_{v}(s,X_{s},P_{(X_{\cdot\wedge s},u)},u_{s}))v_{s}d\lambda\}dB_{s},\\
R_{t}^{2}&=\int_{0}^{t}\Bigr\{\int_{0}^{1}\widetilde{E}[\langle(\partial_{\mu}\sigma)(s,X^{\lambda,\rho}_{s},P_{\vartheta_{s}^{\lambda,\rho}};\widetilde{\vartheta}_{s}^{\lambda,\rho})- (\partial_{\mu}\sigma)(s,X_{s},P_{(X_{\cdot\wedge s},u)};(\widetilde{X}_{\cdot\wedge s},\widetilde{u})),(\frac{1}{\rho}(\widetilde{X}^{\rho}_{\cdot\wedge s}-\widetilde{X}_{\cdot\wedge t}),\widetilde{v})\rangle_{1}]d\lambda\Bigr\}dB_{s}.
\end{aligned}
\end{equation*}
We begin with estimating $R_{t}^{1}$. % and $R_{t}^{2}$.
 From the Burkholder-Davis-Gundy inequality and thanks to our Lipschitz conditions we have
\begin{equation}
\begin{aligned}
&E[\sup_{t\in[0,T]}|R_{t}^{1}|^{2}] \\ %=E[\sup_{t\in[0,T]}|\int_{0}^{t}\{\int_{0}^{1}(\sigma_{x}(s,X^{\lambda,\rho}_{s},P_{\vartheta_{s}^{\lambda,\rho}},u_{s}^{\lambda,\rho})-\sigma_{x}(s,X_{s},P_{(X_{\cdot\wedge s},u)},u_{s}))\frac{1}{\rho}(X^{\rho}_{s}-X_{s})d\lambda\}dB_{s}|^{2}]\\
\leq& \ CE\[\int_{0}^{T}\|\int_{0}^{1}(\sigma_{x}(t,X^{\lambda,\rho}_{t},P_{\vartheta_{t}^{\lambda,\rho}})-\sigma_{x}(t,X_{t},P_{(X_{\cdot\wedge t},u)}))\frac{1}{\rho}(X^{\rho}_{t}-X_{t})d\lambda\|^{2}dt\], \\
\leq&\ CE\[\int_{0}^{T}\|(|X^{\rho}_{t}-X_{t}|+W_{2}(P_{\vartheta_{t}^{\lambda,\rho}},P_{(X_{\cdot\wedge t},u)}))\frac{1}{\rho}(X^{\rho}_{t}-X_{t})\|^{2}dt\]\\
\leq&\ C\rho^{2}E[\int_{0}^{T}|\frac{1}{\rho}(X^{\rho}_{t}-X_{t})|^{4}dt\]
+E\[\int_{0}^{T}W^{2}_{2}(P_{\vartheta_{t}^{\lambda,\rho}},P_{(X_{\cdot\wedge t},u)})\Bigr|\frac{1}{\rho}(X^{\rho}_{t}-X_{t})\|^{2}dt\],\ \rho\in(0,1),\ \lambda\in(0,1).\\
%&\ +CE\[\int_{0}^{T}|\rho v_{t}|^{2}|\frac{1}{\rho}(X^{\rho}_{t}-X_{t})|^{2}dt\]\\
%\leq&\ C\rho^{2}E[\int_{0}^{T}|\frac{1}{\rho}(X^{\rho}_{t}-X_{t})|^{4}dt]
%+E[\int_{0}^{T} |\frac{1}{\rho}(X^{\rho}_{s}-X_{t})|^{2}]\times E[|\vartheta_{t}^{\lambda,\rho}-(X_{\cdot\wedge t}, u)|^{2}_{\mathcal{C}_{T}\times\mathcal{U}_{T}}]\\
%&\ +CE\[\int_{0}^{t}|\rho v_{t}|^{2}\|\frac{1}{\rho}(X^{\rho}_{t}-X_{t})\|^{2}dt\],\ \rho>0,\ \lambda\in[0,1].
\end{aligned}
\end{equation}
%using estimation formula \eqref{SDEestimates1},
Thus,
\begin{equation}\label{W_2estimates1}
\begin{aligned}
E[\sup_{t\in[0,T]}|R_{t}^{1}|^{2}]\leq&\ C\rho^{2}E\[\sup_{t\in[0,T]}\|\frac{1}{\rho}(X^{\rho}_{t}-X_{t})\|^{4}\]\\
&\ +C\rho^{2}E\[\sup_{t\in[0,T]}|\frac{1}{\rho}(X^{\rho}_{t}-X_{t})\|^{2}\]\times E\[\sup_{t\in[0,T]}\|\frac{1}{\rho}(X^{\rho}_{ t}-X_{t})\|^{2}+\int_{0}^{T}|v_{t}|^{2}dt\],
\end{aligned}
\end{equation}
and using H\"{o}lder's inequality, combined with \eqref{SDEestimates1_1} we obtain
%\begin{equation*}
%\begin{aligned}
%&E[\sup_{t\in[0,T]}|\frac{1}{\rho}(X^{\rho}_{t}-X_{t})|^{2}\int_{0}^{T}|v_{t}|^{2}dt]\\
%\leq& \ %C(E[\sup_{t\in[0,T]}|\frac{1}{\rho}(X^{\rho}_{t}-X_{t})|^{4}])^{\frac{1}{2}}\times(E[(\int_{0}^{T}|v_{t}|^{2}dt)^{2}])^{\frac{1}{2}}\\
%\leq&\ C(E[\sup_{t\in[0,T]}|\frac{1}{\rho}(X^{\rho}_{t}-X_{t})|^{4}])^{\frac{1}{2}}\times(E[\int_{0}^{T}|v_{t}|^{4}dt])^{\frac{1}{2}},
%\end{aligned}
%\end{equation*}

\begin{equation}\label{eq_R_t^1}
\begin{aligned}
\Bigr(E\Bigr[\sup_{t\in[0,T]}|R_{t}^{1}|^{2}\Bigr]\Bigr)^{\frac{1}{2}}\leq C\rho.
\end{aligned}
\end{equation}
%In analogy to $R_{t}^{1}$, we have
%\begin{equation}
%\begin{aligned}
%E[\sup_{0\in[0,T]}|R_{t}^{2}|^{2}]
%=&\E[\sup_{0\in[0,T]}|\int_{0}^{t}\{\int_{0}^{1}(\sigma_{v}(s,X^{\lambda,\rho}_{s},P_{\vartheta_{s}^{\lambda,\rho}},u_{s}^{\lambda,\rho})-\sigma_{v}(s,X_{s},P_{(X_{\cdot\wedge s},u)},u_{s}))v_{s}d\lambda\}dB_{s}|^{2}]\\
%\leq& \ CE[\int_{0}^{T}|\int_{0}^{1}(\sigma_{v}(t,X^{\lambda,\rho}_{t},P_{\vartheta_{t}^{\lambda,\rho}},u_{t}^{\lambda\rho})-\sigma_{v}(t,X_{t},P_{(X_{\cdot\wedge t},u)},u_{t}))v_{t}d\lambda|^{2}dt],\\
%\leq&\ CE[\int_{0}^{T}(|X_{t}^{\rho}-X_{t}|^{2}+W^{2}_{2}(P_{\vartheta_{t}^{\lambda,\rho}},P_{(X_{\cdot\wedge t},u)})+|\rho v_{t}|^{2})+)|v_{t}|^{2}dt]\\
%\leq&\  C\rho^{2}E[\sup_{t\in[0,T]}|\frac{1}{\rho}(X^{\rho}_{t}-X_{t})|^{2}\int_{0}^{T}|v_{t}|^{2}dt]+C\rho^{2}E[\int_{0}^{T}|v_{t}|^{4}dt]\\
%&\ +C\rho^{2}E[\sup_{t\in[0,T]}|\frac{1}{\rho}(X^{\rho}_{ t}-X_{t})|^{2}+\int_{0}^{T}|v_{t}|^{2}dt]\times E[\int_{0}^{T}|v_{t}|^{2}dt]\\
%\leq&\ C\rho^{2}.
%\end{aligned}
%\end{equation}
Let us now estimate $R_{t}^{2}$. For this end we decompose the integrand of $R_{t}^{2}$,
\begin{equation*}
\begin{aligned}
\widetilde{E}[\langle(\partial_{\mu}\sigma)(s,X^{\lambda,\rho}_{s},P_{\vartheta_{s}^{\lambda,\rho}};\widetilde{\vartheta}_{s}^{\lambda,\rho})- (\partial_{\mu}\sigma)(s,X_{s},P_{(X_{\cdot\wedge s},u)};(\widetilde{X}_{\cdot\wedge s},\widetilde{u})),(\frac{1}{\rho}(\widetilde{X}^{\rho}_{\cdot\wedge s}-\widetilde{X}_{\cdot\wedge s}),\widetilde{v})\rangle_{1}]=R_{s}^{2,1}+R_{s}^{2,2},
\end{aligned}
\end{equation*}
where
\begin{equation*}
\begin{aligned}
R_{s}^{2,1}:=&\ \widetilde{E}[\langle(\partial_{\mu}\sigma)(s,X^{\lambda,\rho}_{s},P_{\vartheta_{s}^{\lambda,\rho}};\widetilde{\vartheta}_{s}^{\lambda,\rho})
-(\partial_{\mu}\sigma)(s,X_{s},P_{\vartheta_{s}^{\lambda,\rho}};\widetilde{\vartheta}_{s}^{\lambda,\rho}),(\frac{1}{\rho}(\widetilde{X}^{\rho}_{\cdot\wedge s}-\widetilde{X}_{\cdot\wedge s}),\widetilde{v})\rangle_{1}],\\
R_{s}^{2,2}:=&\ \widetilde{E}[\langle(\partial_{\mu}\sigma)(s,X_{s},P_{\vartheta_{s}^{\lambda,\rho}};\widetilde{\vartheta}_{s}^{\lambda,\rho})- (\partial_{\mu}\sigma)(s,X_{s},P_{(X_{\cdot\wedge s},u)};(\widetilde{X}_{\cdot\wedge s},\widetilde{u})),(\frac{1}{\rho}(\widetilde{X}^{\rho}_{\cdot\wedge s}-\widetilde{X}_{\cdot\wedge s}),\widetilde{v})\rangle_{1}].
\end{aligned}
\end{equation*}
As
\begin{equation*}
\begin{aligned}
R_{s}^{2,1}
=&\ \widetilde{E}\[\int_{0}^{T}((\partial_{\mu}\sigma)_{1}(s,X^{\lambda,\rho}_{s},P_{\vartheta_{s}^{\lambda,\rho}};\widetilde{\vartheta}_{s}^{\lambda,\rho})(r)
-(\partial_{\mu}\sigma)_{1}(s,X_{s},P_{\vartheta_{s}^{\lambda,\rho}};\widetilde{\vartheta}_{s}^{\lambda,\rho})(r))
\frac{1}{\rho}(\widetilde{X}^{\lambda,\rho}_{r\wedge s}-\widetilde{X}_{r\wedge s})dr\\
&\ +\int_{0}^{T}((\partial_{\mu}\sigma)_{2}(s,X^{\lambda,\rho}_{s},P_{\vartheta_{s}^{\lambda,\rho}};\widetilde{\vartheta}_{s}^{\lambda,\rho})(r)
-(\partial_{\mu}\sigma)_{2}(s,X_{s},P_{\vartheta_{s}^{\lambda,\rho}};\widetilde{\vartheta}_{s}^{\lambda,\rho})(r))
\widetilde{v}(r)dr\],
\end{aligned}
\end{equation*}
thanks to the Lipschitz assumption on $\partial_{\mu}\sigma$ and \eqref{SDEestimates1_1} we obtain
\begin{equation*}
\begin{aligned}
E\[\sup_{t\in[0,T]}\|\int_{0}^{t}R_{s}^{2,1}dB_{s}\|^{2}\]
&\leq\ C\rho^{2}E\[\sup_{t\in[0,T]}|\frac{1}{\rho}(X^{\rho}_{t}-X_{t})\|^{2}\]\times E\[\sup_{t\in[0,T]}\|\frac{1}{\rho}(X^{\rho}_{ t}-X_{t})\|^{2}+\int_{0}^{T}|v_{t}|^{2}dt\]\\
&\leq\ C\rho^{2},\ \rho\geq0,\ \lambda\in[0,1].
\end{aligned}
\end{equation*}

%Then, we have
%\begin{equation*}
%\begin{aligned}
%(E[\sup_{t\in[0,T]}|\int_{0}^{t}R_{s}^{3,1}ds|^{2}])^{\frac{1}{2}}\leq C\rho,\ \rho\geq0.
%\end{aligned}
%\end{equation*}
Similar to the above estimate, for $R_{s}^{2,2}$ we have
\begin{equation*}
\begin{aligned}
&|R_{s}^{2,2}|^{2}
=|\widetilde{E}[\langle(\partial_{\mu}\sigma)(s,X_{s},P_{\vartheta_{s}^{\lambda,\rho}};\widetilde{\vartheta}_{s}^{\lambda,\rho})- (\partial_{\mu}\sigma)(s,X_{s},P_{(X_{\cdot\wedge s},u)};(\widetilde{X}_{\cdot\wedge s},\widetilde{u})),(\frac{1}{\rho}(\widetilde{X}^{\rho}_{\cdot\wedge s}-\widetilde{X}_{\cdot\wedge s}),\widetilde{v})\rangle_{1}]|^{2}\\
%\leq& \ C|\widetilde{E}[|(\partial_{\mu}\sigma)(s,X_{s},P_{\vartheta_{s}^{\lambda,\rho}},u_{s};\widetilde{\vartheta}_{s}^{\lambda,\rho})- (\partial_{\mu}\sigma)(s,X_{s},P_{(X_{\cdot\wedge s},u)},u_{s};(\widetilde{X}_{\cdot\wedge s},\widetilde{u}))|^{2}_{BV_{T}\times\mathcal{U}_{T}}]\times \widetilde{E}[|(\frac{1}{\rho}(\widetilde{X}^{\rho}_{\cdot\wedge s}-\widetilde{X}_{\cdot\wedge s}),\widetilde{v})|^{2}_{\mathcal{C}_{T}\times\mathcal{U}_{T}}]\\
%\leq&\ C\widetilde{E}[W_{2}^{2}(P_{\vartheta_{s}^{\lambda,\rho}},P_{(X_{\cdot\wedge s},u)})+|\widetilde{\vartheta}_{s}^{\lambda,\rho}-(\widetilde{X}_{\cdot\wedge s}, \widetilde{u})|^{s}_{\mathcal{C}_{T}\times\mathcal{U}_{T}}]\\
%\leq&\ C\widetilde{E}[|\widetilde{\vartheta}_{s}^{\lambda,\rho}-(\widetilde{X}_{\cdot\wedge s}, \widetilde{u})|^{2}_{\mathcal{C}_{T}\times\mathcal{U}_{T}}]\\
%\leq&\ C\rho^{2}\widetilde{E}[\sup_{t\in[0,T]}|\frac{1}{\rho}(\widetilde{X}^{\rho}_{ t}-\widetilde{X}_{t})|^{2}+\int_{0}^{T}|v_{t}|^{2}dt]\ (\lambda\in[0,1])\\
\leq& C\(\widetilde{E}\[\(W_{2}(P_{\vartheta_{s}^{\lambda,\rho}},P_{(X_{\cdot\wedge s},u)})
\!+ \!\sup_{r\leq s}|\widetilde{X}^{\rho}_{r}-\widetilde{X}_{r}|\!+\! \rho\(\int_{0}^{T}|\widetilde{v}_{s}|^{2}ds\)^{\frac{1}{2}}\)
\(\!\int_{0}^{T}\!\|\frac{1}{\rho}(\widetilde{X}^{\rho}_{r\wedge s}\!-\!\widetilde{X}_{r\wedge s})\|dr\!+\!\int_{0}^{T}|\widetilde{v}_{s}|ds\)\]\)^{2}\\
\leq& C\rho^{2}\(\widetilde{E}\[\sup_{r\leq T}\|\frac{1}{\rho}(\widetilde{X}^{\rho}_{r}-\widetilde{X}_{r})\|^{2}+\int_{0}^{T}|\widetilde{v}_{s}|^{2}ds\]\)^{2}\\
\leq&\ C\rho^{2},\ \rho\geq0,
\end{aligned}
\end{equation*}
%C是随行变化的，依赖于T与L.
and, consequently,
\begin{equation}
\begin{aligned}
&(E[\sup_{t\in[0,T]}|R_{t}^{2}|^{2}])^{\frac{1}{2}}\leq C\rho,\ \rho >0.
\end{aligned}
\end{equation}
Let $X_{t}(\rho):=\frac{1}{\rho}(X^{\rho}_{t}-X_{t})-X_{t}^{1}$.
Then, from \eqref{eq_5.25_1}
\begin{equation*}
\begin{aligned}
X_{t}(\rho)=&\ R_{t}^{1}%+R_{t}^{2}
+R_{t}^{2}+\int_{0}^{t}\sigma_{x}(s,X_{s},P_{(X_{\cdot\wedge s},u)})X_{s}(\rho)dB_{s}\\
&\ +\int_{0}^{t}\widetilde{E}[\langle
(\partial_{\mu}\sigma)(s,X_{s},P_{(X_{\cdot\wedge s},u)};(\widetilde{X}_{\cdot\wedge s},\widetilde{u})), (\widetilde{X}_{\cdot\wedge s}(\rho),0)\rangle_{1}]dB_{s},
\end{aligned}
\end{equation*}
and as $\sigma_{x}$ and $\partial_\mu\sigma$ are bounded,
%in $L^{2}(\Omega, \mathcal{F},P; BV_{T}\times \mathcal{U}_{T})$,
it follows from Gronwall's inequality that
\begin{equation}\label{eq_X*0}
E[\sup_{t\in[0,T]}|X_{t}(\rho)|^{2}]\leq CE[\sup_{t\in[0,T]}|R_{t}^{1}|^{2}
%+\sup_{t\in[0,T]}|R_{t}^{2}|^{2}
+\sup_{t\in[0,T]}|R_{t}^{2}|^{2}]\leq C\rho^{2},\ \rho>0.
\end{equation}
Consequently, $\frac{1}{\rho}(X^{\rho}-X)\rightarrow X^{1}$ in ${S}^{2}_{\mathbb{F}}$ as $\rho\rightarrow0$ .

%%%%%%%%%%%%%%下面看BSDE
Using the notations introduced for the proof of Lemma \ref{lem5.1.1}, the BSDE in \eqref{Eq_control2} can be rewritten as follows:
\begin{equation}
\begin{aligned}
\widehat{Y}_{t}= \widehat{Y}_{T}+\int_{t}^{T}(f(\theta_{s}^{\rho})-f(\theta_{s}))ds-\int_{t}^{T}\widehat{Z}_{s}dB_s,\ t\in[0,T],
\end{aligned}
\end{equation}
where $\widehat{Y}_{T}=\Phi(X^{\rho}_T,P_{(X^{\rho},u^{\rho})})-\Phi(X_{T},P_{(X,u)})$.

%%%%%%%%%%%%%%%%%%%%%%%%%%%%%%%%%%%%%%%%%%%%%%%%%%%%%%%%%%%%%%%%%
Consequently, from equation \eqref{eq_Chain2} in the proof of Lemma \ref{lem5.1.1} we have
\begin{equation*}
\begin{aligned}
&\ \frac{1}{\rho}(Y_{t}^{\rho}-Y_{t})\\
=&\ \sum_{i=0}^{5}I_{t}^{i}+\int_{t}^{T}f_{x}(\theta_{s})\frac{1}{\rho}(X^{\rho}_{s}-X_{s})ds
  +\int_{t}^{T}f_{y}(\theta_{s})\frac{1}{\rho}(Y^{\rho}_{s}-Y_{s})ds+\int_{t}^{T}f_{z}(\theta_{s})\frac{1}{\rho}(Z^{\rho}_{s}-Z_{s})ds\\
 %+\int_{t}^{T}f_{v}(\theta_{s},u_{s})v_{s}ds\\
&\ +\int_{t}^{T}\widetilde{E}[\langle(\partial_{\mu}f)(\theta_{s};
(\widetilde{X},\widetilde{Y}_{\cdot\vee s},\widetilde{u})) , (\frac{1}{\rho}(\widetilde{X}^{\rho}-\widetilde{X}), \frac{1}{\rho}(\widetilde{Y}_{\cdot\vee s}^{\rho}-\widetilde{Y}_{\cdot\vee s}),\widetilde{v}) \rangle_{2}]ds-\int_{t}^{T}\frac{1}{\rho}(Z_{s}^{\rho}-Z_{s})dB_s\\
&\ +\Phi_{x}(X_{T},P_{(X,u)}))\frac{1}{\rho}(X_{T}^{\rho}-X_{T})
+\widetilde{E}[\langle(\partial_{\mu}\Phi)(X_{T},P_{(X,u)};(\widetilde{X},\widetilde{u})),(\frac{1}{\rho}(\widetilde{X}^{\rho}-\widetilde{X}),\widetilde{v})\rangle_{1}],\\
\end{aligned}
\end{equation*}
where
\begin{equation*}
\begin{aligned}
I^{0}:=&\int_{0}^{1}(\Phi_{x}(X^{\lambda\rho}_T,P_{(X^{\lambda\rho},u^{\lambda\rho})})-\Phi_{x}(X_{T},P_{(X,u)}))\frac{1}{\rho}(X^{\rho}_{T}-X_{T}))d\lambda,\\
I^{1}:=&\int_0^1\widetilde{E}[\langle(\partial_{\mu}\Phi)(X^{\lambda\rho}_T,P_{(X^{\lambda\rho},u^{\lambda\rho})};(\widetilde{X}^{\lambda\rho},\widetilde{u}^{\lambda\rho}))-(\partial_{\mu}\Phi)(X_{T},P_{(X,u)};(\widetilde{X},\widetilde{u})),(\frac{1}{\rho}(\widetilde{X}^{\rho}-\widetilde{X}),\widetilde{v})\rangle_{1}]d\lambda,\\
I_{t}^{2}:=&\int_{t}^{T}\int_{0}^{1}(f_{x}(\theta^{\lambda\rho}_{s})-f_{x}(\theta_{s}))\frac{1}{\rho}(X^{\rho}_{s}-X_{s})d\lambda ds,\\
I_{t}^{3}:=&\int_{t}^{T}\int_{0}^{1}(f_{y}(\theta^{\lambda\rho}_{s})-f_{y}(\theta_{s}))\frac{1}{\rho}(Y^{\rho}_{s}-Y_{s})d\lambda ds,\\
I_{t}^{4}:=&\int_{t}^{T}\int_{0}^{1}(f_{z}(\theta^{\lambda\rho}_{s})-f_{z}(\theta_{s}))\frac{1}{\rho}(Z^{\rho}_{s}-Z_{s})d\lambda ds,\\
%I_{t}^{5}:=&\int_{t}^{T}\int_{0}^{1}(f_{v}(\theta^{\lambda\rho}_{s})-f_{v}(\theta_{s}))v_{s}d\lambda ds,\\
I_{t}^{5}:=&\int_{t}^{T}\int_{0}^{1}\widetilde{E}[\langle((\partial_{\mu}f)(\theta^{\lambda\rho}_{s};
(\widetilde{X}^{\lambda\rho},\widetilde{Y}^{\lambda\rho}_{\cdot\vee s},\widetilde{u}^{\lambda\rho})) -(\partial_{\mu}f)(\theta_{s};
(\widetilde{X},\widetilde{Y}_{\cdot\vee s},\widetilde{u})) ), (\frac{1}{\rho}(\widetilde{X}^{\rho}-\widetilde{X}),\\ &\frac{1}{\rho}(\widetilde{Y}_{\cdot\vee s}^{\rho}-\widetilde{Y}_{\cdot\vee s}),\widetilde{v})) \rangle_{2}]d\lambda ds.
\end{aligned}
\end{equation*}
In analogy to the estimate of $R_{t}^{1}$, by using the Lipschitz continuity of $\Phi_{x}$, we see that
\begin{equation*}
\begin{aligned}
&E[|\Phi_{x}(X^{\lambda,\rho}_T,P_{(X^{\lambda,\rho},u^{\lambda\rho})})-\Phi_{x}(X_{T},P_{(X,u)}))\frac{1}{\rho}(X^{\rho}_{T}-X_{T})|^{2}]\\
\leq&\ CE[|(|X^{\lambda,\rho}_T-X_{T}|+ W_{2}(P_{(X^{\lambda,\rho},u^{\lambda\rho})},P_{(X,u)})) \frac{1}{\rho}(X^{\rho}_{T}-X_{T})|^{2}]
\leq\ C\rho^{2},\ \lambda\in[0,1],
\end{aligned}
\end{equation*}
and so $E[|I^{0}|^{2}]\leq C\rho^{2}$.
For $I_{t}^{3}$, using the Lipschitz continuity of $f_{y}$ and the definition of $\theta^{\lambda,\rho}$ and $\theta$ we have
\begin{equation*}
\begin{aligned}
E[\sup_{0\in[0,T]}|I_{t}^{3}|^{2}]
\leq &\ E\[\int_{0}^{T}\int_{0}^{1}\(\frac{1}{\rho}|Y^{\rho}_{t}-Y_{t}|\)^{2}\bigr\{(|X^{\rho}_{s}-X_{s}|+|Y^{\rho}_{s}-Y_{s}|+|Z^{\rho}_{s}-Z_{s}| \\ %+|v_{s}|
&\ + W_{2}(P_{(X^{\lambda,\rho},Y^{\lambda,\rho}_{\cdot\vee s},u^{\lambda\rho})}, P_{(X,Y_{\cdot\vee s},u)}))\bigr\}^{2}d\lambda ds \]\\
\leq &\ CE\[\sup_{t\in[0,T]}\frac{1}{\rho^{2}}|Y^{\rho}_{t}-Y_{t}|^{4}\]+\frac{1}{\rho^{2}}CE\[\(\int_{0}^{T}(|X^{\rho}_{s}-X_{s}|+|Y^{\rho}_{s}-Y_{s}|+|Z^{\rho}_{s}-Z_{s}|)^{2}ds\)^{2}\]\\
&\ +\frac{1}{\rho^{2}}C\(\int_{0}^{T}\int_{0}^{1}W_{2}(P_{(X^{\lambda,\rho},Y^{\lambda,\rho}_{\cdot\vee s},u^{\lambda\rho})}, P_{(X,Y_{\cdot\vee s},u)})^{2}d\lambda ds\)^{2}\\
\leq&\ CE\[\sup_{t\in[0,T]}\frac{1}{\rho^{2}} (|X^{\rho}_{s}-X_{s}|^{4}+|Y^{\rho}_{s}-Y_{s}|^{4})\]+CE\[\frac{1}{\rho^{2}}\(\int_{0}^{T}|Z^{\rho}_{s}-Z_{s}|^{2}ds\)^{2}\]\\
&\ +\rho^{2}\(E\[\int_{0}^{T}|v_{s}|^{2}ds\]\)^{2}.
\end{aligned}
\end{equation*}

%The similar as \eqref{W_2estimates1} and \eqref{eq_R_t^1}, we get,
Thus, using \eqref{Eq21-2-2_1} we obtain
\begin{equation}\label{eq_I_t^3}
\begin{aligned}
\Bigr(E\Bigr[\sup_{0\in[0,T]}|I_{t}^{3}|^{2}\Bigr]\Bigr)^{\frac{1}{2}}\leq C\rho.
\end{aligned}
\end{equation}
The terms $I_{t}^{2}$ and $I_{t}^{5}$ can be estimated with analogous arguments, and so
\begin{equation}\label{eq_I_t^256}
\begin{aligned}
\Bigr(E\Bigr[\sup_{0\in[0,T]}|I_{t}^{i}|^{2}\Bigr]\Bigr)^{\frac{1}{2}}\leq C\rho,\ i=2,\ 5.
\end{aligned}
\end{equation}
Similarly we also show that
\begin{equation}\label{eq_I_t111}
\begin{aligned}
\Bigr(E\Bigr[\sup_{0\in[0,T]}|I_{t}^{1}|^{2}\Bigr]\Bigr)^{\frac{1}{2}}\leq C\rho.
\end{aligned}
\end{equation}
%Now we left the estimation formulas of $I_{t}^{4}$ and $I_{t}^{6}$ to be proved. It's clearly that $(E[\sup_{0\in[0,T]}|I_{t}^{6}|^{2}])^{\frac{1}{2}}\leq C\rho$. For $I_{t}^{4}$,
It remains to estimate $I_{t}^{4}$. For this we remark that
\begin{equation*}
\begin{aligned}
E[\sup_{t\in[0,T]}|I_{t}^{4}|^{2}]
%=&E[\sup_{t\in[0,T]}|\int_{t}^{T}\int_{0}^{1}(f_{z}(\theta^{\lambda\rho}_{s},u^{\lambda\rho}_{s})-f_{z}(\theta_{s},u_{s}))\frac{1}{\rho}(Z^{\rho}_{s}-Z_{s})d\lambda ds|^{2}],\ \lambda\in[0,1],\\
%\leq& CE[\sup_{t\in[0,T]}|\int_{t}^{T}(|\theta^{\rho}_{t}-\theta_{t}|+\rho|v_{t}|)\frac{1}{\rho}(Z^{\rho}_{t}-Z_{t})dt|^{2}]\\
\leq&\ C\rho^{2}E\[\(\sup_{t\in[0,T]}\|\frac{1}{\rho}(X^{\rho}_{t}-X_{t})\|^{2}+\sup_{t\in[0,T]}\|\frac{1}{\rho}(Y^{\rho}_{t}-Y_{t})\|^{2}
\)\int_{0}^{T}\|\frac{1}{\rho}(Z^{\rho}_{t}-Z_{t})\|^{2}dt\]\\
&\ +\!C\rho^{2}E\[\int_{0}^{T}\|\frac{1}{\rho}(Z^{\rho}_{t}-Z_{t})\|^{2}dt\]E\[\sup_{t\in[0,T]}\|\frac{1}{\rho}(X^{\rho}_{ t}-X_{t})\|^{2}\!+\!\sup_{t\in[0,T]}\|\frac{1}{\rho}(Y^{\rho}_{ t}\!-\!Y_{t})\|^{2}\!+\!\int_{0}^{T}\!|v_{s}|^{2}ds\]\\
\end{aligned}
\end{equation*}
\begin{equation*}
\begin{aligned}
&\ + C\rho^{2}E\[(\int_{0}^{T}\|\frac{1}{\rho}(Z^{\rho}_{t}-Z_{t})\|^{2}dt)^{2}\],
\end{aligned}
\end{equation*}
and thanks to \eqref{SDEestimates1_1} and \eqref{Eq21-2-2_1} we obtain
\begin{equation}\label{eq_I_t^4}
\begin{aligned}
\Bigr(E\Bigr[\sup_{t\in[0,T]}|I_{t}^{4}|^{2}\Bigr]\Bigr)^{\frac{1}{2}}\leq C\rho.
\end{aligned}
\end{equation}
%%%%%%%%%%%%%%%%后面Y*_{t}

Setting now $Y_{t}(\rho):=\frac{1}{\rho}(Y^{\rho}_{t}-Y_{t})-Y_{t}^{1},\ Z_{t}(\rho):=\frac{1}{\rho}(Z^{\rho}_{t}-Z_{t})-Z_{t}^{1}$ and recalling that $X_{t}(\rho):=\frac{1}{\rho}(X^{\rho}_{t}-X_{t})-X_{t}^{1}$, we have
\begin{equation*}
\begin{aligned}
Y_{t}(\rho)=&\ \sum_{i=0}^{5}I_{t}^{i}+\Phi_{x}(X_{T},P_{(X,u)}))\widetilde{X}_{T}(\rho)
+\widetilde{E}[\langle(\partial_{\mu}\Phi)(X_{T},P_{(X,u)};(\widetilde{X},\widetilde{u})),\widetilde{X}(\rho),0)\rangle_{1}]\\
&\ +\int_{t}^{T}f_{x}(\theta_{s})X_{s}(\rho)ds
  +\int_{t}^{T}f_{y}(\theta_{s})Y_{s}(\rho)ds
  +\int_{t}^{T}f_{z}(\theta_{s})Z_{s}(\rho)ds\\
&\ +\int_{t}^{T}\widetilde{E}[\langle(\partial_{\mu}f)(\theta_{s};
(\widetilde{X},\widetilde{Y}_{\cdot\vee s},\widetilde{u})) , (\widetilde{X}(\rho), \widetilde{Y}_{\cdot\vee s}(\rho)),0)) \rangle_{2}]ds-\int_{t}^{T}Z_{s}(\rho)dB_s,\ t\in[0,T].
\end{aligned}
\end{equation*}
%and as \eqref{eq_X*0}, $f_{x},f_{y},f_{z}$ be bounded, $\partial_\mu f$ is bounded in $L^{2}(\Omega, \mathcal{F},P; BV_{T}\times BV_{T}\times \mathcal{U}_{T})$ and by using Gronwall's inequality, we finally obtain
The solution $(Y(\rho),Z(\rho))$ of the above equation can be estimated in the same way as \eqref{Eq_5.33.1}. This yields
\begin{equation*}
\begin{aligned}
E\Bigr[\sup_{t\in[0,T]}|Y_{t}(\rho)|^{2}+\int_{0}^{T}|Z_{t}(\rho)|^{2}dt\Bigr]\leq CE\Bigr[\sup_{t\in[0,T]}|\sum_{i=0}^{5}I_{t}^{i}|^{2}\Bigr]+ CE\Bigr[\sup_{t\in[0,T]}|X_{t}(\rho)|^{2}\Bigr].
\end{aligned}
\end{equation*}
Finally, from our estimates of $I^{i}_{t},\ 1\leq i \leq 5$, and from \eqref{eq_X*0} it follows that $$E\Bigr[\sup_{t\in[0,T]}|Y_{t}(\rho)|^{2}+\int_{0}^{T}|Z_{t}(\rho)|^{2}dt\Bigr]\leq C\rho^{2},\ \rho\in(0,1).$$
Consequently, $(\frac{1}{\rho}(Y^{\rho}-Y)\rightarrow X^{1},\frac{1}{\rho}(Z^{\rho}-Z)\rightarrow Z^{1})$ in ${S}^{2}_{\mathbb{F}}\times {M}^{2}_{\mathbb{F}}$ as $\rho\rightarrow0$ .
\end{proof}

%From above on, we have proved equation \eqref{eq_vari} has a unique solution $(X^1,\ Y^1,\ Z^1)$. And the solution $(\widehat{X},\widehat{Y},\widehat{Z})$ converges to $(X^1,\ Y^1\ Z^1$ in $\mathcal{S}^{2}_{\mathbb{F}}\times \mathcal{S}^{2}_{\mathbb{F}}\times \mathcal{M}^{2}_{\mathbb{F}}$\ as $\rho$ tends to $0$.
Let us now study the so-called variational inequality. For this note that, because $u(\cdot)$\ is an optimal control, it holds
\begin{equation}\label{eq_J_v}
\rho^{-1}[J(u(\cdot)+\rho v(\cdot))-J(u(\cdot))]\geq 0.
\end{equation}
Thus thanks to the Lemmas \ref{lem5.1} and \ref{lem5.2} we have
\begin{theorem}\label{theo_variin}
We suppose \textnormal{(H5)} and \textnormal{(H6)} are satisfied. Then, the following variational inequality holds true:
\begin{equation}\label{eq_variin}
\begin{aligned}
0\leq & E\Bigr[\int_{0}^{T}(L_{x}(\theta_{t})X^{1}_{t}+L_{y}(\theta_{t})Y^{1}_{t}+L_{z}(\theta_{t})Z^{1}_{t}%+L_{v}(\theta_{t},u_{t})v{t}\\
+\widetilde{E}[\langle(\partial_{\mu}L)(\theta_{t};(\widetilde{X},\widetilde{Y},\widetilde{u})), (\widetilde{X}^{1},\widetilde{Y}^{1},\widetilde{v})\rangle_{2}] )dt\\
&\quad +\varphi_{x}(X_{T},P_{(X,Y,u)})X_{T}^{1}+\widetilde{E}[\langle(\partial_{\mu}\varphi)(X_{T},P_{(X,Y,u)};(\widetilde{X},\widetilde{Y},\widetilde{u})), (\widetilde{X}^{1},\widetilde{Y}^{1},\widetilde{v})\rangle_{2}]\Bigr],
\end{aligned}
\end{equation}
where $\theta_{t}=(t,X_{t},Y_{t},Z_{t},P_{(X,Y_{\cdot\vee t},v)})$.
\end{theorem}
\begin{proof}
Letting $\rho\rightarrow0$ in \eqref{eq_J_v}, and using Lemma \ref{lem5.2}, we obtain \eqref{eq_variin}.
\end{proof}

In order to derive the maximum principle, we need to introduce the adjoint equation. For this, we make use of the following
notations:
\begin{equation}\label{adj_notation}
\begin{aligned}
&\sigma_{x}(t):=\sigma_{x}(t,X_{t},P_{(X_{\cdot\wedge t},u)});\\
&(\partial_{\mu}\sigma)^{*}_{1}(t)[k]
:=E\Bigr[\widetilde{E}\Bigr[\int_{t}^{T}\{(\partial_{\mu}\sigma)_{1}(r,\widetilde{X}_{r},P_{(X_{\cdot\wedge r},u)};(X_{\cdot\wedge r},u))(t)\widetilde{k}_{r}
+(\partial_{\mu}\sigma)_{1}(t,\widetilde{X}_{t},P_{(X_{\cdot\wedge t},u)};(X_{\cdot\wedge t},u))(r)\widetilde{k}_{r}\}dr
\Bigr]\Bigr|\mathcal{F}_{t}\Bigr];\\
&(\partial_{\mu}\sigma)^{*}_{2}(t)[k]
:=E\Bigr[\widetilde{E}\Bigr[\int_{0}^{T}
(\partial_{\mu}\sigma)_{2}(r,\widetilde{X}_{r},P_{(X_{\cdot\wedge r},u)};(X_{\cdot\wedge r},u))(t)\widetilde{k}_{r}dr
\Bigr]\Bigr|\mathcal{F}_{t}\Bigr];\\
&(\partial_{\mu}f)_{j}^{*}(t)[p]
:=E\Bigr[\widetilde{E}\Bigr[\int_{0}^{T}
(\partial_{\mu}f)_{j}(\widetilde{\theta}_{r};(X,Y_{\cdot\vee r},u))(t)\widetilde{p}(r)dr\Bigr]\Bigr|\mathcal{F}_{t}\Bigr],\ j=1,3;\\
&(\partial_{\mu}f)_{2}^{*}(t)[p]
:=E\Bigr[\widetilde{E}\Bigr[\int_{0}^{t}
\Bigr\{(\partial_{\mu}f)_{2}(\widetilde{\theta}_{r};(X,Y_{\cdot\vee r},u))(t)\widetilde{p}(r)
+(\partial_{\mu}f)_{2}(\widetilde{\theta}_{t};(X,Y_{\cdot\vee t},u))(r)\widetilde{p}(t)\Bigr\}dr\Bigr]\Bigr|\mathcal{F}_{t}\Bigr];\\
&(\partial_{\mu}L)_{j}^{*}(t)
:=E\Bigr[\widetilde{E}\Bigr[\int_{0}^{T}
(\partial_{\mu}L)_{j}(\widetilde{\theta}_{r};(X,Y,u))(t)dr\Bigr]\Bigr|\mathcal{F}_{t}\Bigr],\ j=1,2,3;\\
&(\partial_{\mu}\varphi)_{j}^{*}(t)
:=E\Bigr[\widetilde{E}\Bigr[
(\partial_{\mu}\varphi)_{j}(\widetilde{X}_{T},P_{(X,Y,u)};(X,Y,u))(t)\Bigr|\mathcal{F}_{t}],\ j=1,2,3;\\
&(\partial_{\mu}\Phi)_{j}^{*}(t)[p(T)]
:=E[\widetilde{E}[
(\partial_{\mu}\Phi)_{j}(\widetilde{X}_{T},P_{(X,u)};(X,u))(t)\widetilde{p}(T)|\mathcal{F}_{t}],\ j=1,2.
\end{aligned}
\end{equation}
Observe that the above processes are bounded and $\mathbb{F}$-adapted. With the notations introduced above we consider the following adjoint FBSDE,
\begin{equation}\label{eq_adjSDE}
\left\{
\begin{aligned}
dp(t)=&\ \{f_{y}(\theta_{t})p(t)+(\partial_{\mu}f)^{*}_{2}(t)[p]
-L_{y}(\theta_{t})-(\partial_{\mu}L)^{*}_{2}(t)-(\partial_{\mu}\varphi)^{*}_{2}(t)\}dt\\
&\ +\{f_{z}(\theta_{t})p(t)-L_{z}(\theta_{t})\Bigr\}dB_{t},\ t\in [0,T],\\
p(0)=&0,
\end{aligned}
\right.
\end{equation}
\begin{equation}\label{eq_adjBSDE}
\left\{
\begin{aligned}
dq(t)=&\ -\Bigr\{\sigma_{x}(t)k(t)+(\partial_{\mu}\sigma)^{*}_{1}(t)[k]-f_{x}(\theta_{t})p(t)\\
&\ -(\partial_{\mu}f)^{*}_{1}(t)[p]+L_{x}(\theta_{t})+(\partial_{\mu}L)^{*}_{1}(t)
+(\partial_{\mu}\varphi)^{*}_{1}(t)-(\partial_{\mu}\Phi)^{*}(t)[p(T)] \}dt\\
&\ +k(t)dB_{t},\ t\in [0,T],\\
q(T)=&\ \varphi_{x}(X_{T},P_{(X,Y,u)})-\Phi_{x}(X_{T},P_{(X,u)})p(T).
\end{aligned}
\right.
\end{equation}
%With \textnormal{(H5)} and \rm{(H6)},  all derivatives bounded. So also all terms in \eqref{adj_notation} are bounded.
\begin{lemma}
Under the assumptions \textnormal{(H5)} and \rm{(H6)} the adjoint equation \eqref{eq_adjSDE}-\eqref{eq_adjBSDE} has a unique adapted solution $(p,q,k)\in {S}^{2}_{\mathbb{F}}\times {S}^{2}_{\mathbb{F}}\times {M}^{2}_{\mathbb{F}}$.
\end{lemma}
\begin{proof}
Equation \eqref{eq_adjSDE} is an affine mean-field forward stochastic equation with delay (Recall the definition of $(\partial_{\mu}f)_{2}^{*}(t)[p]$), thus the existence and the uniqueness is obvious. After getting the solution of \eqref{eq_adjSDE}, equation \eqref{eq_adjBSDE} is a mean-field BSDE with anticipation like \eqref{eq_variBSDE} (Recall the definition of $(\partial_{\mu}\sigma)_{1}^{*}(t)[k]$). From Lemma \ref{lem5.1}, we see that there exists a unique solution $(q,k)\in {S}^{2}_{\mathbb{F}}\times {M}^{2}_{\mathbb{F}}$. The proof is complete.
\end{proof}

Before giving the maximum principle, we provide the following lemma which is important.
\begin{lemma}\label{Lemma5.5}
Let $p$ be the solution to the adjoint SDE \eqref{eq_adjBSDE}, $(q,k)$ the solution to the adjoint BSDE \eqref{eq_adjBSDE}, and $(X^{1},Y^{1},Z^{1})$ the solution to \eqref{eq_vari}. Then we have
\begin{equation}\label{eq_XqYp}
\begin{aligned}
 E[X_{T}^{1}q(T)+Y_{T}^{1}p(T)]
=&\ E\[\int_{0}^{T}X_{t}^{1}
\Bigr\{(\partial_{\mu}\Phi)^{*}_{1}(t)[p(T)]-L_{x}(\theta_{t})
-(\partial_{\mu}L)^{*}_{1}(t)-(\partial_{\mu}\varphi)^{*}_{1}(t)\}dt\]\\
&\ -E\[\int_{0}^{T}Y_{t}^{1}\big\{L_{y}(\theta_{t})+(\partial_{\mu}L)_{2}^{*}(t)
+(\partial_{\mu}\varphi)_{2}^{*}(t)\big\}dt\]-E\[\int_{0}^{T}Z_{t}^{1} L_{z}(\theta(t))dt\]\\
&\ +E\[\int_{0}^{T}v_{t}\{(\partial_{\mu}\sigma)^{*}_{2}(t)[k]- (\partial_{\mu}f)^{*}_{3}(t)[p] \}dt\]
\end{aligned}
\end{equation}
\end{lemma}
\begin{proof}
The proof of this lemma follows immediately from It\^{o}'s formula and Fubini's theorem, so we omit it.
\end{proof}

Set $\nu=P_{(X,u)}, \mu= P_{(X,Y,u)}$. The terminal value $Y_{T}=\Phi(X_{T},P_{(X,u)})$ and also $\varphi(X_{T},P_{(X,Y,u)})$ depend on the law of the whole path $(X,u)$ and $(X,Y,u)$, respectively. This has as consequence that they produce their own coefficients
$(\partial_{\mu}\Phi)_{j}^{*}(t)[p(T)]$ and $(\partial_{\mu}\varphi)^{*}_{i}(t),\ j=1,2,\ i=1,2,3$, which we have to
take into account in the definition of the Hamiltonian. It adds that the derivatives with respect to the measure $(\partial_\mu\sigma)_1^*(t)[k],\, (\partial_\mu f)_i^*(t)[p],\, i=1,2,$ depend on the whole solution process of the adjoint FBSDE.
This makes that our Hamiltonian cannot be defined in the classical way. We define the Hamiltonian just as follows:
\begin{equation}\label{Eq_Hamiltonian5.50}
\begin{aligned}
H(t,x,y,z,\nu,\mu):=(-f(t,x,y,z,\mu),\ \sigma(t,x,\nu),\ L(t,x,y,z,\mu),-\Phi(\cdot,\nu),\varphi(\cdot,\mu)).
 \end{aligned}
\end{equation}
For the derivatives of the Hamiltonian, we introduce the following notations. For $\nu=P_{(X,u)}, \mu= P_{(X,Y,u)}$,
$\sigma(t)=\sigma(t,X_{t},P_{(X_{\cdot\wedge t},u)}),$ $\sigma_{x}(t)=\sigma_{x}(t,X_{t},P_{(X_{\cdot\wedge t},u)}),$ $
f_{x}(\theta_{t})=f_{x}(t,X_{t},Y_{t},Z_{t},\mu),$ and $L_{x}(\theta_{t})=f_{x}(t,X_{t},Y_{t},Z_{t},\mu)$, we put
$(\partial_{x}H)(t)=(-f_{x}(\theta_{t}), \sigma_{x}(t),L_{x}(\theta_{t}),0,0)$ and $H_{x}(t)=((\partial_{x}H)(t),(p(t),k(t),1,0,0))_{\mathbb{R}^{5}}$,
where $(\cdot,\cdot)_{\mathbb{R}^{5}}$ denotes the inner product in $\mathbb{R}^{5}$.

In the same sense we define
\begin{equation*}
\begin{aligned}
(\partial_{y}H)(t)=(-f_{y}(\theta_{t}),0,L_{y}(\theta_{t}),0,0),\
(\partial_{z}H)(t)=(-f_{z}(\theta_{t}),0,L_{z}(\theta_{t}),0,0),
\end{aligned}
\end{equation*}
and
\begin{equation*}
\begin{aligned}
H_{y}(t)=((\partial_{y}H)(t),(p(t),k(t),1,0,0))_{\mathbb{R}^{5}},\
H_{z}(t)=((\partial_{z}H)(t),(p(t),k(t),1,0,0))_{\mathbb{R}^{5}}.
\end{aligned}
\end{equation*}
Concerning the derivatives with respect to the measure, we write
\begin{equation*}
\begin{aligned}
(\partial_{\mu}H)_{x}(t)=&-(\partial_{\mu}f)^{*}_{1}(t)[p]+(\partial_{\mu}\sigma)^{*}_{1}(t)[k]+(\partial_{\mu}L)^{*}_{1}(t)
-(\partial_{\mu}\Phi)^{*}_{1}(t)[p(T)]+(\partial_{\mu}\varphi)^{*}_{1}(t),\\
(\partial_{\mu}H)_{y}(t)=&-(\partial_{\mu}f)^{*}_{2}(t)[p]+(\partial_{\mu}L)^{*}_{2}(t)+(\partial_{\mu}\varphi)^{*}_{2}(t),\\
(\partial_{\mu}H)_{v}(t)=&-(\partial_{\mu}f)^{*}_{3}(t)[p]+(\partial_{\mu}\sigma)^{*}_{2}(t)[k]+(\partial_{\mu}L)^{*}_{3}(t)
-(\partial_{\mu}\Phi)^{*}_{2}(t)[p(T)]+(\partial_{\mu}\varphi)^{*}_{3}(t).
\end{aligned}
\end{equation*}
Putting
\begin{equation}\label{notation0708}
\begin{aligned}
\mathcal{D}_{x}H(t):=&H_{x}(t)+(\partial_{\mu}H)_{x}(t),\
\mathcal{D}_{y}H(t):=H_{y}(t)+(\partial_{\mu}H)_{y}(t),\\
\mathcal{D}_{z}H(t):=&H_{z}(t),\
\mathcal{D}_{v}H(t):=(\partial_{\mu}H)_{v}(t).
\end{aligned}
\end{equation}
We can write now the adjoint FBSDE \eqref{eq_adjSDE}-\eqref{eq_adjBSDE} in the following way:
\begin{equation}\label{Eq_adjointFBSDE2}
\left\{
\begin{aligned}
dp(t)=& -\mathcal{D}_{y}H(t)dt-\mathcal{D}_{z}H(t)dB_{t},\ t\in[0,T],\\
dq(t)=& -\mathcal{D}_{x}H(t)dt+k(t)dB_{t},\ t\in [0,T],\\
p(0)=&0,\ q(T)=\ \varphi_{x}(X_{T},P_{(X,Y,u)})-\Phi_{x}(X_{T},P_{(X,u)})p(T).
\end{aligned}
\right.
\end{equation}

Then we will get the following stochastic maximum principle:
\begin{theorem}\label{Thm5.1}
Let $u$ be an optimal control of the mean-field FBSDE control problem. Then, recalling the definition of $\mathcal{D}_{v}H(t)$, we have the maximum principle:
\begin{equation}\label{eqSMP}
\begin{aligned}
\mathcal{D}_{v}H(t)(v-u(t))\geq 0,\ \text{for all}\ v\in U,\ dtdP\mbox{-}a.e.
 \end{aligned}
\end{equation}
where $(p,q,k)$ is the solution of the adjoint equation \eqref{eq_adjSDE} and \eqref{eq_adjBSDE}.
\end{theorem}
\begin{proof}
Applying It\^{o}'s formula to $(X^{1}_{t}q(t)+Y^{1}_{t}p(t))$, we get \eqref{eq_XqYp}. Recall that $Y_{T}=\Phi(X_{T},P_{(X,u)})$, and so
\begin{equation*}
\begin{aligned}
Y^{1}_{T}=&\ \Phi_{x}(X_{T},P_{(X,u)})X^{1}_{T}+\widetilde{E}[\langle(\partial_{\mu}\Phi(X_{T},P_{(X,u)};(\widetilde{X},\widetilde{u})), (\widetilde{X}^{1},\widetilde{v}))\rangle_{1}]\\
=&\ \Phi_{x}(X_{T},P_{(X,u)})X^{1}_{T}+\widetilde{E}\[\int_{0}^{T}(\partial_{\mu}\Phi)_{1}(X_{T},P_{(X,u)};(\widetilde{X},\widetilde{u}))(r)\widetilde{X}^{1}_{r}dr\]\\
&\ +\widetilde{E}\Bigr[\int_{0}^{T}(\partial_{\mu}\Phi)_{2}(X_{T},P_{(X,u)};(\widetilde{X},\widetilde{u}))(r)\widetilde{v}_{r}dr\Bigr].
\end{aligned}
\end{equation*}
Then, we have
\begin{equation}\label{eq_YT1p}
\begin{aligned}
&\ E[X_{T}^{1}q(T)+Y_{T}^{1}p(T)]\\
=&\ E[X^{1}_{T}(q(T)+\Phi_{x}(X_{T},P_{(X,u)})p(T))]
+\int_{0}^{T}X_{T}^{1}E[\widetilde{E}[(\partial_{\mu}\Phi)_{1}(\widetilde{X}_{T},P_{(X,u)};(X,u))(t)\widetilde{p}(T)]|\mathcal{F}_{t}]dt]\\
&\ +E\Bigr[\int_{0}^{T}v_{t}E[\widetilde{E}[(\partial_{\mu}\Phi)_{2}(\widetilde{X}_{T},P_{(X,u)};(X,u))(t)\widetilde{p}(T)]|\mathcal{F}_{t}]dt\Bigr]\\
=&\ E[X_{T}^{1}(q(T)+\Phi_{x}(X_{T},P_{(X,u)})p(T))]
+E\Bigr[\int_{0}^{T}X_{T}^{1}(\partial_{\mu}\Phi)_{1}^{*}(t)[p(T)]dt\Bigr]
+E\Bigr[\int_{0}^{T}v_{t}(\partial_{\mu}\Phi)_{2}^{*}(t)[p(T)]dt\Bigr].
\end{aligned}
\end{equation}
Recall $q(T)=\varphi_{x}(X_{T},P_{(X,Y,u)})-\Phi_{x}(X_{T},P_{(X,u)})p(T)$ and the definition of $(\partial_{\mu}\Phi)_{j}^{*}(t)[p(T)]$. So
\begin{equation*}
\begin{aligned}
&\ E[X^{1}_{T}\varphi_{x}(X_{T},P_{(X,Y,u)})]=E[X^{1}_{T}(q(T)+\Phi_{x}(X_{T},P_{(X,u)})p(T))]\\
=&\ E[X^{1}_{T}q(T)+Y^{1}_{T}p(T)]-E\Bigr[\int_{0}^{T}X^{1}_{T}(\partial_{\mu}\Phi)_{1}^{*}(t)[p(T)]dt\Bigr]
-E\[\int_{0}^{T}v_{t}(\partial_{\mu}\Phi)_{2}^{*}(t)[p(T)]dt\].
\end{aligned}
\end{equation*}
We substitute equation \eqref{eq_XqYp} into the above relation. The new formula we get for $E[X^{1}_{T}\varphi_{x}(X_{T},P_{(X,Y,u)})]$ is now combined with the variational inequality \eqref{eq_variin}. For $v\in \mathcal{U}_{ad}$  such that $(u+ v)\in \mathcal{U}_{ad}$, we obtain
\begin{equation*}
\begin{aligned}
0\leq \partial_{v}J(u)
=&\ E\[\int_{0}^{T}v_{t}\Bigr\{-(\partial_{\mu}f)^{*}_{3}(t)[p]+(\partial_{\mu}\sigma)^{*}_{2}(t)[k]+(\partial_{\mu}L)^{*}_{3}(t)
-(\partial_{\mu}\Phi)^{*}_{2}(t)[p(T)]+(\partial_{\mu}\varphi)^{*}_{3}(t)
\Bigr\}dt\]\\
&\ +E\Bigr[\int_{0}^{T}(X_{t}^{1}-X_{T}^{1})(\partial_{\mu}\Phi)_{1}^{*}(t)[p(T)]]dt\Bigr].
\\
\end{aligned}
\end{equation*}
However, as $(\partial_{\mu}\Phi)_{1}^{*}(t)[p(T)]$ is $\mathcal{F}_{t}$-measurable and
$X^{1}$ is a martingale, this latter expression disappears. Consequently, recalling the notation $\mathcal{D}_{v}H(t)$, we see that the preceding formula just states
\begin{equation}
E\[\int_{0}^{T} \mathcal{D}_{v}H(t)\cdot v_{t} dt\]\geq 0,\ dtdP\mbox{-}a.e.,\ \mbox{for all}\ v\in\mathcal{U}_{ad}\ \mbox{such that}\ u+v\in\mathcal{U}_{ad}.
\end{equation}
Then we conclude by using a standard argument that
\begin{equation*}
\mathcal{D}_{v}H(t)(v-u(t))\geq 0,\ \text{for all}\ v\in U,\ dtdP\mbox{-}a.e.
\end{equation*}
\end{proof}

\section{A sufficient condition for optimality}

Last but not least, in this section we show that the optimality condition given by
 Pontryagin's SMP \eqref{eqSMP} is not only necessary, but, combined with a suitable convexity assumption
 for the Hamiltonian $H$, it is also sufficient. More precisely, we have the following.

\begin{theorem}\label{Thmsufficient}
Let us suppose the convexity of the Hamiltonian $(-f(t,x,y,z,\mu)p(t),\sigma(t,x,\nu)k(t),L(t,x,y,z,\mu),$
$-\Phi(x',\nu)p(T),\varphi(x',\mu))$  in $(x,x',y,z,\nu,\mu)\in \mathbb{R}^4\times\mathcal{P}_2(\mathcal{C}_{T}\times\mathcal{U}_T)\times\mathcal{P}_2(\mathcal{C}_{T}^2\times\mathcal{U}_T)$,
where $(p,q,k)$ is the solution of the adjoint equation \eqref{eq_adjSDE}-\eqref{eq_adjBSDE}. Furthermore, we continue to suppose the standard assumptions (H5)-(H6)
of the preceding Section 5. Then, if an admissible control $u\in\mathcal{U}_{ad}$ satisfies \eqref{eqSMP},
it is optimal:
\begin{equation}\nonumber
\begin{aligned}
J(w)\geq J(u),\ \text{for all}\ w\in\mathcal{U}_{ad} \ \text{such that}\ v:=w-u\in\mathcal{U}_{ad}.
\end{aligned}
\end{equation}
\end{theorem}

\begin{remark}
The convexity of $(x,x',y,z,\mu,\nu)\rightarrow(-f(t,x,y,z,\mu)p(t),-\Phi(x',\nu)p(T))$ can be got by supposing,
for instance, the convexity of $-f(t,\cdot,\cdot,\cdot,\cdot)$ and $-\Phi$, and by taking assumptions on the coefficients of
the adjoint forward SDE \eqref{eq_adjSDE} which guarantee that $p(t)\geq 0$, $t\in[0,T]$.
\end{remark}

\begin{proof}[Proof (of Theorem \ref{Thmsufficient})]
Let $w\in\mathcal{U}_{ad}$ and $v:=w-u\in\mathcal{U}_{ad}$. From the
convexity of $(x,x',y,z,\mu)\rightarrow \varphi(x',\mu)+L(t,x,y,z,\mu)$, $t\in[0,T]$,
with using the notations of  the preceding section, we have
\begin{equation}\label{eq070801}
\begin{aligned}
&\ J(w)-J(u)\\
=&\ E\big[\varphi(X^{w}_{T},P_{(X^{w},w)})-\varphi(X_{T},P_{(X,u)})+\int_{0}^{T}\big(L(\theta^{w}_t)-L(\theta_t)\big)dt\big]\\
\geq&\ E\Bigr[\varphi_{x}(X_{T},P_{(X,u)})(X^{w}_{T}-X_{T})+\widetilde{E}[\langle(\partial_{\mu}\varphi)(X_{T},P_{(X,Y,u)};(\widetilde{X},\widetilde{Y},\widetilde{u}))
,(\widetilde{X}^{w}-\widetilde{X},\widetilde{Y}^{w}-\widetilde{Y},\widetilde{v})\rangle_{2}]\Bigr]\\
&\ +E\Bigr[\int_{0}^{T}\Big(L_{x}(\theta_{t})(X^{w}_{t}-X_{t})+L_{y}(\theta_{t})(Y^{w}_{t}-Y_{t})+L_{z}(\theta_{t})(Z^{w}_{t}-Z_{t})\\
&\ +\widetilde{E}[\langle(\partial_{\mu}L)(\theta_{t};(\widetilde{X},\widetilde{Y},\widetilde{u})), (\widetilde{X}^{w}-\widetilde{X},\widetilde{Y}^{w}-\widetilde{Y},\widetilde{v})\rangle_{2}] \Big)dt\Bigr].
\end{aligned}
\end{equation}
Here $({X}^{w},{Y}^{w},{Z}^{w})$ denotes the solution \eqref{Eqcontrol} associated
with the control $w\in\mathcal{U}_{ad}$, $\theta^{w}_t=(X^{w}_t,Y^{w}_t,Z^{w}_t,P_{(X^{w},Y^{w},w)})$,
and $(X,Y,Z)=(X^u,Y^u,Z^u)$, $\theta=\theta^u$. On the other hand, from the convexity of $-\Phi(x',\nu)p(T)$,
and using the notations introduced in \eqref{adj_notation}, as well as equation \eqref{Eqcontrol}, we get
\begin{equation}
\begin{aligned}
&\ E[\Phi_{x}(X_{T},P_{(X,u)})(X^{w}_{T}-X_{T})p(T)]\\
\geq&\ E[(\Phi(X^{w}_{T},P_{(X^{w},w)})-\Phi(X_{T},P_{(X,u)}))p(T)]
-E[\widetilde{E}[\langle(\partial_{\mu}\Phi)(X_{T},P_{(X,u)};(\widetilde{X},\widetilde{u}))
,(\widetilde{X}^{w}-\widetilde{X},\widetilde{v})\rangle_{1}]p(T)]\\
=&\ E[p(T)(Y^{w}_{T}-Y_{T})]
 -E\Bigr[\int_{0}^{T}(\partial_{\mu}\Phi)_{1}^{*}(t)[p(T)](X^{w}_{t}-X_{t})dt\Bigr]
-E\[\int_{0}^{T}(\partial_{\mu}\Phi)_{2}^{*}(t)[p(T)]v_tdt\Bigr].
\end{aligned}
\end{equation}
Thus, recalling from \eqref{eq_adjBSDE} that $\varphi_{x}(X_{T},P_{(X,Y,u)})=q(T)+\Phi_{x}(X_{T},P_{(X,u)})p(T)$, we obtain
\begin{equation}\label{eq070802}
\begin{aligned}
&\ E[\varphi_{x}(X_{T},P_{(X,Y,u)})(X^{w}_{T}-X_{T})]\\
\geq&\ E[q(T)(X^{w}_{T}-X_{T})+p(T)(Y^{w}_{T}-Y_{T})]
 -E\Bigr[\int_{0}^{T}(\partial_{\mu}\Phi)_{1}^{*}(t)[p(T)](X^{w}_{t}-X_{t})dt\Bigr]\\
&\ -E\[\int_{0}^{T}(\partial_{\mu}\Phi)_{2}^{*}(t)[p(T)]v_tdt\Bigr].
\end{aligned}
\end{equation}
Applying It\^{o}'s formula to $q(t)(X^{w}_{t}-X_{t})+p(t)(Y^{w}_{t}-Y_{t})$,
a straight-forward computations yields
\begin{equation}\label{eq070803}
\begin{aligned}
&\ E[q(T)(X^{w}_{T}-X_{T})+p(T)(Y^{w}_{T}-Y_{T})]\\
=&\ E\Bigr[\int_{0}^{T} (-\mathcal{D}_{x}H(t)(X^{w}_t-X_t)-\mathcal{D}_{y}H(t)(Y^{w}_t-Y_t)-\mathcal{D}_{z}H(t)(Z^{w}_t-Z_t) )dt\Bigr]\\
&\ + E\Bigr[\int_{0}^{T}(-(f(\theta^{w}_{t})-f(\theta_t))p(t)+ (\sigma(t,X^{w}_t,P_{(X^{w}_{\cdot\wedge t},w)})-\sigma(t,X_t,P_{(X_{\cdot\wedge t},u)}))k(t)    )dt\Bigr].
\end{aligned}
\end{equation}
Recall that the notations of $\mathcal{D}_{x}H(t), \mathcal{D}_{y}H(t), \mathcal{D}_{z}H(t)$ have
been introduced in \eqref{notation0708}. Consequently,  substituting \eqref{eq070803} in \eqref{eq070802}
and the estimate thus obtained in \eqref{eq070801}, we obtain by using that
\begin{equation}\nonumber
\begin{aligned}
&\ E[\widetilde{E}[\langle(\partial_{\mu}\varphi)(X_{T},P_{(X,Y,u)};(\widetilde{X},\widetilde{Y},\widetilde{u}))
,(\widetilde{X}^{w}-\widetilde{X},\widetilde{Y}^{w}-\widetilde{Y},\widetilde{v})\rangle_{2}]]\\
=&\ E\Bigr[\int_{0}^{T} \big((\partial_{\mu}\varphi)_{1}^{*}(t)(X^{w}_{t}-X_{t})+(\partial_{\mu}\varphi)_{2}^{*}(t)(Y^{w}_{t}-Y_{t})
+(\partial_{\mu}\varphi)_{3}^{*}(t)v_{t}\big)dt\Bigr]
\end{aligned}
\end{equation}
and the similar formula for $\partial_{\mu}L$ (see \eqref{adj_notation}), the following
\begin{equation}\label{AAA}
\begin{aligned}
&\ J(w)-J(u)\\
\ge&\ E\Bigr[\int_{0}^{T}\bigr\{
(X^{w}_{t}-X_{t})\bigr(-\mathcal{D}_{x}H(t)-(\partial_{\mu}\Phi)^{*}_{1}(t)[p(T)]+(\partial_{\mu}\varphi)^{*}_{1}(t)
+L_{x}(\theta_{t})+(\partial_{\mu}L)^{*}_{1}(t)\bigr)\\
&\ +\!(Y^{w}_{t}\!-\!Y_{t})\bigr( -\mathcal{D}_{y}H(t)\!+\!(\partial_{\mu}\varphi)^{*}_{2}(t) \!+\!L_{y}(\theta_{t})\!+\!(\partial_{\mu}L)^{*}_{2}(t)\bigr)
 \!+\!(Z^{w}_{t}\!-\!Z_{t})\bigr( -\!\mathcal{D}_{z}H(t)\!+\! L_{z}(\theta_{t}) \bigr)
   \bigr\}dt\Bigr]\\
&\ + E\Bigr[\int_{0}^{T}\bigr( -(f(\theta^{w}_{t})-f(\theta_t))p(t)+ (\sigma(t,X^{w}_t,P_{(X^{w}_{\cdot\wedge t},w)})-\sigma(t,X_t,P_{(X_{\cdot\wedge t},u)}))k(t) \bigr)dt\Bigr]\\
&\ +E\Bigr[\int_{0}^{T}v_{t}\bigr(-(\partial_{\mu}\Phi)^{*}_{2}(t)[p(T)]+(\partial_{\mu}\varphi)^{*}_{3}(t)+(\partial_{\mu}L)^{*}_{3}(t)
\bigr)dt\Bigr]
\end{aligned}
\end{equation}
Recalling the definition of $\mathcal{D}_{x}H(t), \mathcal{D}_{y}H(t), \mathcal{D}_{z}H(t)$, the right-hand side of \eqref{AAA} can be writen as follows:
 \begin{equation}%\label{EqconvexityL}
\begin{aligned}
J(w)-J(u)
\geq I_{1}+I_{2} + E\Bigr[\int_{0}^{T}\mathcal{D}_{v}H(t)v_{t} dt\Bigr],
\end{aligned}
\end{equation}
where
 \begin{equation}%\label{EqconvexityL}
\begin{aligned}
I_{1}=&\ E\Bigr[\int_{0}^{T}\big\{
-\Bigr[f(\theta^{w}_{t})-f(\theta_t)-f_{x}(\theta_{t})(X^{w}_{t}-X_{t})-f_{y}(\theta_{t})(Y^{w}_{t}-Y_{t})
-f_{z}(\theta_{t})(Z^{w}_{t}-Z_{t})\Bigr]p(t)\\
&+(\partial_{\mu}f)_{1}^{*}(t)[p](X^{w}_{t}-X_{t})+(\partial_{\mu}f)_{2}^{*}(t)[p](Y^{w}_{t}-Y_{t})
+(\partial_{\mu}f)_{3}^{*}(t)[p]v_{t}\big\}dt\Bigr]\\
\end{aligned}
\end{equation}
and
 \begin{equation}%\label{EqconvexityL}
\begin{aligned}
I_2=&\ E\Bigr[\int_{0}^{T}\big\{\big[(\sigma(t,X^{w}_t,P_{(X^{w}_{\cdot\wedge t},w)}))-\sigma(t,X_t,P_{(X_{\cdot\wedge t},u)})- \sigma_{x}(t)(X^{w}_{t}-X_{t}) \big]k(t)\\
&\ -(\partial_{\mu}\sigma)^{*}_{1}(t)[k](X^{w}_{t}-X_{t})-(\partial_{\mu}\sigma)^{*}_{2}(t)[k]v_{t}
   \big\}dt\Bigr].\\
\end{aligned}
\end{equation}
But, as
\begin{equation}%\label{EqconvexityL}
\begin{aligned}
&E\Bigr[\int_{0}^{T}\big\{(\partial_{\mu}f)_{1}^{*}(t)[p](X^{w}_{t}-X_{t})+(\partial_{\mu}f)_{2}^{*}(t)[p](Y^{w}_{t}-Y_{t})
+(\partial_{\mu}f)_{3}^{*}(t)[p]v_{t} \big\}dt\Bigr]\\
=&\ E\Bigr[\int_{0}^{T} \widetilde{E}[\langle(\partial_{\mu}f)(\theta_{t};(\widetilde{X},\widetilde{Y},\widetilde{u})), (\widetilde{X}^{w}-\widetilde{X},\widetilde{Y}^{w}-\widetilde{Y},\widetilde{v})\rangle_{2}]p(t)
dt\Bigr]
\end{aligned}
\end{equation}
and
\begin{equation}%\label{E}
\begin{aligned}
 &E\Bigr[\int_{0}^{T}
\big\{(\partial_{\mu}\sigma)^{*}_{1}(t)[k](X^{w}_{t}-X_{t})+(\partial_{\mu}\sigma)^{*}_{2}(t)[k]v_{t}\big\}dt\Bigr]\\
=&\ E\Bigr[\int_{0}^{T} \widetilde{E}[\langle
(\partial_{\mu}\sigma)(t,X_{t},P_{(X_{\cdot\wedge t},u)};(\widetilde{X}_{\cdot\wedge t},\widetilde{u})), (\widetilde{X}^{w}_{\cdot\wedge t}-\widetilde{X}_{\cdot\wedge t},\widetilde{v})\rangle_{1}]k(t)dt\Bigr],
\end{aligned}
\end{equation}
it follows from the convexity of $(-f(x,y,z,\mu)p(t),\sigma(t,x,\nu)k(t))$ that $I_1\geq0$ and $I_2\geq0$.
Consequently, thanks to our assumption $\mathcal{D}_{v}H(t)v_t\geq0$, we obtain $J(w)-J(u)\geq0$.
\end{proof}

%\clearpage


\begin{thebibliography}{0}

%\bibitem{AC2021}
%N. Agrama and S.E. Choutrib. Mean-field FBSDE and optimal control. \emph{Stochastic Analysis and Applications},39 (2), 235-251, 2021.

%\bibitem{AD}
%D.~Andersson and B.~Djehiche. A maximum principle for SDEs of mean-field type. \emph{Appl. Math. Optim.}, 63, 341-356, 2011.

\bibitem{A}
F.~Antonelli. Backward-forward stochastic differential equations. \emph{Ann. Appl.} 3, 777-793, 1993.


%\bibitem{B1}
%R.~Bellman. Dynamic programming. \emph{Princeton Univ.Press.}, 1957.

%\bibitem{B2}
%A.~Bensoussan. Lectures on stochastic control. In: Lecture Notes in Mathematics. \emph{Springer}, Berlin, 972, 1-62, 1981.

%\bibitem{BDL}
%R.~Buckdahn, B.~Djehiche and J.~Li. A general stochastic maximum principle for SDEs of mean-field type. \emph{Appl. Math. Optim.}, 64, 197-216, 2011.

\bibitem{BDLP2009}
R.~Buckdahn, B.~Djehiche, J.~Li and S.~Peng. Mean-field backward stochastic differential equations: A limit approach. \emph{Ann. Probab.}, 37, 1524-1565, 2009.

\bibitem{BLM}
R. Buckdahn, J. Li, J. Ma, A stochastic maximum principle for general mean-field systems. \emph{Applied Mathematics \& Optimization}, 74 (3), 507-534, 2016.

\bibitem{BLP1}
R.~Buckdahn, J.~Li, and S.~Peng. Mean-field backward stochastic differential equations and related patial differential equations. \emph{Stochastic Processes and their Applications}, 119, 3133-3154, 2009.

\bibitem{CP2013}
P. Cardaliaguet. Notes on mean field games. From P.-L. Lions Lectures at Coll\`{e}ge de France.(2013).
Available at http://www.college-de-france.fr.

\bibitem{CD2018I}
R. Carmona and F. Delarue. Probabilistic Theory of Mean Field Games with Applications I. Springer, 2018.

%\bibitem{CD2018II}
%R. Carmona and F. Delarue. Probabilistic Theory of Mean Field Games with Applications II. Springer, 2018.


%\bibitem{CDLL2019} P.~Cardaliaguet, F.~Delarue, J. M.~Lasry and P.L.~Lions. The master equation and the convergence
   %problem in mean field games:(ams-201). \emph{Princeton University Press}. 2019.

\bibitem{DEdwards2011}
D.A. Edwards. On the Kantorovich-Rubinstein theorem. \emph{Expositiones Mathematicae}, 29, 387-398, (2011).

\bibitem{D}
F.~Delarue. A forward-backward stochastic algorithm for quasi-linear PDEs. \emph{Annals of Applied Probability}, 16, 140-184, 2006.

%\bibitem{EKHM2008}
%N. El Karoui, S. Hamad\`{e}ne, and A. Matoussi. Backward stochastic differential equations and applications. \emph{Indifference pricing: theory and applications}. 267-320, 2008.

%\bibitem{FT2018}{M. Fujii, A. Takahashi. Quadratic-exponential growth BSDEs with jumps and their Malliavin's differentiability. \emph{Stochastic Processes and their Applications}, 128 (6), 2083-2130, 2018.}

%\bibitem{HL2014}{T. Hao, J. Li. Backward stochastic differential equations coupled with value function and related optimal control problems. \emph{In Abstract and Applied Analysis}, 2014, 1-17, 2014. }

%\bibitem{H}
%U.G.~Haussmann. A stochastic maximum principle for optimal control of diffusions.  \emph{Longman Scientific and Technical, Essex}, 1986.


\bibitem{Hu2017}
M. Hu. Stochastic global maximum principle for optimization with recursive utilities. \emph{Probability, Uncertainty and Quantitative Risk}, 2 (1), 1-20, 2017.

%\bibitem{HJX2018}
%{M. Hu, S. Ji, X. Xue. A global stochastic maximum principle for fully coupled forward-backward stochastic systems. \emph{SIAM Journal on Control and Optimization}, 56 (6), 4309-4335, 2018.}


\bibitem{HP}
Y.~Hu, and S.~Peng. Solution of forward-backward stochastic differential equations. \emph{Probab. Theory Related Fields}, 103, 273-283, 1995.

%\bibitem{KH1982}
% H.~Kunita. Stochastic differential equations and stochastic flows of diffeomorphisms.Ecole d'\'{e}t\'{e}
%de Probabilit\'{e} de Saint-Flour-1982.\emph{ Lect. Notes Math.}, 1097, 143-303, 1982.

%\bibitem{K1}
%H.J.~Kushner. On the stochastic maximum principle: fixed time of control. \emph{J. Math. Anal. Appl.}, 11, 78-92, 1965

%\bibitem{K2}
%H.J.~Kushner. Necessary conditions for continuous parameter stochastic optimization problems. \emph{SIAM J. Control Optim.}, 10, 550-565, 1972.

%\bibitem{KS1991}
%I. Karatzas, S.E. Shreve. Brownian motion and stochastic calculus[M]. \emph{Springer Science \& Business Media}, 1991.

%\bibitem{L}
%J.~Li. Stochastic maximum principle in mean-field controls. \emph{Automatica}, 366-373, 2012.

%\bibitem{L2018}
%{J. Li. Mean-field forward and backward SDEs with jumps and associated nonlocal quasi-linear integral-PDEs. \emph{Stochastic Processes and their Applications}, 128 (9), 3118-3180, 2018}

\bibitem{LH2016}
J. Li and M. Hui. Weak solutions of mean-field stochastic differential equations and application to zero-sum stochastic differential games. \emph{SIAM Journal on Control and Optimization}, 54 (3), 1826-1858, 2016.



\bibitem{LW2015}
{J. Li, Q.M. Wei. Stochastic differential games for fully coupled FBSDEs with jumps. Applied Mathematics \& Optimization, 71 (3), 411-448, 2015.}

%\bibitem{MMS2019}
%{A. Matoussi, A. Manai, R. Salhi. Mean-Field Backward-Forward SDE with Jumps and Storage problem in Smart Grids. arXiv preprint arXiv:1906.08525. 2019.}

\bibitem{MPY}
J.~Ma, P.~Protter, J.~Yong. Solving forward-backward stochastic differential equations explicitly-a four step scheme. \emph{Probab Related Fields}, 98, 339-359, 1994.

%\bibitem{MY}
%J.~Ma, J.M.~Yong. Forward-backward stochastic differential equations and their applications. \emph{Springer},
%Berlin, 1999.

%\bibitem{MWZZ}
%J.~Ma, Z.~Wu, D.Z. Zhang, J.F. Zhang. On wellposedness of forward-backward SDEs---A unified approach. \emph{Annals of Applied Probability}, 25(4), 2168-2214, 2015. http://arxiv.org/abs/1110.4658, 2011.
 %Hindawi.2014.

\bibitem{N2006}
D. Nualart. The Malliavin Calculus and Related Topics. \emph{Berlin, Heidelberg,
New York: Springer}, 1995.

\bibitem{PP1}
 E.~Pardoux, S.~Peng. Adapted solution of a backward stochastic differential equation. \emph{Systems Control Lett}, 14, 55-61, 1990.

\bibitem{PT}
E.~Pardoux, S.~Tang. Forward-backward stochastic differential equations and quasilinear parabolic PDEs. \emph{Probab. Theory Relat. Fields}, 114, 123-150, 1999.

%\bibitem{P1990}
%S.~Peng. A general stochastic maximum principle for optimal control promblems. \emph{SIAM J. Contr. Optim.}, 28(4), 966-979, 1990.


%\bibitem{Peng1993}
%S. Peng. Backward stochastic difffferential equations and applications to the optimal control. Appl. Math. Optim. 27 ,1251-144, 1993.

%\bibitem{PBG}
%L.S.~Pontryagin, V.G.~Boltanskii, and R.V.~Gamkrelidze. The mathematical theory of optimal processes. \emph{Inter-scien. N.Y.}, 1962.

%\bibitem{PW}
%S.~Peng, Z.~Wu. Fully coupled forward-backward stochastic differential equations and applications to optimal control. \emph{SIAM J. Control Optim.}, 37, 825-843, 1999.

%\bibitem{Q}
%Y.L.~Qin. Mean-field forward-backward stochastic differential equations and related questions. \emph{Shandong University (Weihai) Master Thesis}, April, 2012.

\bibitem{V2003}
C. Villani. Topics in optimal transportations. \emph{Graduate Studies in Mathematics,
58, AMS, Providence, RI}, 2003.

 \bibitem{W3}
Z.~Wu. Maximum principle for optimal control problem of fully coupled forward-backward stochastic systems. \emph{Systems Sci. math. Sci.}, 11(3), 249-259, 1998.


%\bibitem{WZ2020}
%C. Wu and J. Zhang. Viscosity solutions to parabolic master equations and McKean-VlasovSDEs with closed-loop controls, \emph{Ann. Appl. Probab.}, 30, 936-986, 2020.

%\bibitem{Y}
%J.M.~Yong. Forward-backward stochastic differential equations with mixed initial-terminal conditions. \emph{Tranctions of the American Mathematical Society}, 362 (2), 1047-1096, 2010.

 \bibitem{Z}
D.~Zhang. Forward-backward stochastic differential equations and backward linear quadratic stochastic optimal control problem. \emph{Communications in mathematical research}, 25(5), 402-410, 2009.

\bibitem{ZhangLX2020}
S. Zhang, X. Li, and J.Xiong. A stochastic maximum principle for partially observed stochastic control systems with delay. \emph{Systems and Control Letters}, 146 (2020), 104812, 2020.
%\bibitem{Za1985}
%M. Zakai . The Malliavin calculus.\emph{ Acta Applicandae Mathematica},  3(2), 175-207, 1985.

\end{thebibliography}
\end{document}